\documentclass[a4paper, 12pt]{article}
\usepackage[russian]{babel}
\usepackage{amsmath}
\usepackage{amsfonts}
\newcommand{\mx}{\operatorname{\mathfrak M}}
\newcommand{\mn}{\operatorname{\mathfrak m}}
\renewcommand{\Im}{\operatorname{Im}}
\newcommand{\Ker}{\operatorname{Ker}}
\newcommand{\loc}{\operatorname{loc}}
\newcommand{\e}{\operatorname{\mathfrak e}}
\renewcommand{\span}{\operatorname{span}}
\newcommand{\close}{\operatorname{close}}
\newcommand{\s}{\operatorname{\mathcal s}}
\newcommand{\card}{\operatorname{card}}
\newcommand{\diam}{\operatorname{diam}}
\newcommand{\supp}{\operatorname{supp}}
\newcommand{\supvrai}{\operatornamewithlimits{sup\,vrai}}
\newcommand{\N}{\mathbb N}
\newcommand{\Z}{\mathbb Z}
\newcommand{\R}{\mathbb R}
\renewcommand{\C}{\mathbb C}
\newcommand{\Nu}{\mathcal N}

\newcommand{\mes}{\operatorname{mes}}

\begin{document}

\author{ С. Н. Кудрявцев }
\title{Теорема типа Литтлвуда-Пэли для ортопроекторов
на подпространства всплесков}
\date{}\maketitle
\begin{abstract}
В статье доказано утверждение, представляющее собой
аналог теоремы Литтлвуда-Пэли для ортопроекторов на
подпространства всплесков, соответствующие многомерному
кратно-масштабному анализу, порожденному тензорным произведением
гладких финитных масштабирующих функций одной переменной.
\end{abstract}

Ключевые слова: ортопроектор, подпространства всплесков,
масштабирующая функция, кратно-масштабный анализ,
теорема Литтлвуда-Пэли
\bigskip

\centerline{Введение}

Как известно (см., например, [1], [2] и др. работы), важное значение для
вывода порядковых оценок
точности приближения в $L_p$ способами,
основанными на применении кратных тригонометрических рядов, классов
периодических функций нескольких переменных с условиями на смешанные
производные (разности) имеют теорема Литтлвуда-Пэли для кратных рядов
Фурье и следствия из нее.
По поводу теоремы Литтлвуда-Пэли для кратных рядов Фурье
см., например, [3, п. 1.5.2], [4] и приведенную там литературу.
Для средств приближения классов (непериодических) функций, Заданных
на кубе $ I^d, $ подчиненных условиям на смешанные разности, аналог теоремы
Литтлвуда-Пэли установлен автором в [5].
Для получения соответствующих оценок точности приближения классов
непериодических функций смешанной гладкости, заданных на всем
пространстве $ \R^d, $ полезно иметь аналоги этих утверждений для средств
приближения таких классов функций.
С этой целью в работе доказан аналог теоремы Литтлвуда-Пэли для
семейства ортопроекторов $ \{\mathcal E_\kappa, \kappa \in \Z_+^d \}, $
(см. п. 2.2.) на подпространства всплесков, соответствующие
кратно-масштабному анализу (КМА), порожденному тензорным произведением
гладких финитных функций одной переменной,
а именно, показано, что при $ 1 < p < \infty $ для любой функции
$ f \in L_p(\R^d) $ выполняются неравенства
\begin{equation*} \tag{1}
c_{1} \| f \|_{L_p(\R^d)} \le
( \int_{\R^d} (\sum_{\kappa \in \Z_+^d} |(\mathcal E_\kappa f)(x)|^2 )^{p/2}
dx)^{1/p} \le c_{2} \| f \|_{L_p(\R^d)}
\end{equation*}
с некоторыми константами $ c_1, c_2 > 0, $ зависящими от $ d,p $ и функций,
порождающих КМА.
Задача построения масштабирующих функций и всплесков не является предметом
рассмотрения настоящей работы, и по этому поводу см., например, [6], [7].

Работа состоит из введения и двух параграфов. В \S 1 приведены
предварительные сведения об ортопроекторах, рассматриваемых в
работе. В п. 2.2. \S 2 на основании сведений из \S 1  и п. 2.1.
устанавливается справедливость (1). Перейдем к точным формулировкам

и доказательствам.
\bigskip

\centerline{\S 1. Операторы проектирования на подпространства, }
\centerline{порожденные ограниченными функциями}
\centerline{с компактными носителями}
\bigskip

 1.1. В этом пункте вводятся обозначения, используемые
в настоящей работе, а также приводятся некоторые факты, необходимые
в дальнейшем.

Для $ d \in \N $ через $ \Z_+^d $ обозначим множество
$$
\Z_+^d =\{\lambda =(\lambda_1,  \ldots, \lambda_d) \in \Z^d:
\lambda_j \ge0, j=1, \ldots, d\}.
$$
 Обозначим также при  $ d \in \N $ для $ l \in \Z_+^d $ через
$ \Z_+^d(l) $ множество
 $$
 \Z_+^d(l) =\{ \lambda  \in \Z_+^d: \lambda_j \le l_j, j=1, \ldots,  d\}.
 $$

Для комплексного числа $ z $ через $ \overline z $ обозначается
комплексно-сопряженное число.

Через $ C_1(D) $ обозначается пространство комплекснозначных непрерывно
дифференцируемых функций в области $ D \subset \R^d. $

Для измеримого по Лебегу множества $ D \subset \R^d$
при $ 1 \le p \le \infty $ через  $  L_p(D),$  как обычно,
обозначается пространство  всех   комплекснозначных измеримых на $ D $
функций $f,$ для которых определена норма
\begin{equation*}
\|f\|_{L_p(D)}= \begin{cases} (\int_D |f(x)|^p
dx)^{1/p}, 1 \le p<\infty;\\
\supvrai_{x \in D}|f(x)|, p=\infty. \end{cases}
\end{equation*}
$ L_2(D) $ является гильбертовым пространством со скалярным произведением
$ \langle f, g \rangle = \int_D f \overline g dx, f,g \in L_2(D). $

Обозначим еще через $ L_1^{\loc}(\R^d) $ линейное пространство всех локально
суммируемых в $ \R^d $ функций, принимающих комплексные значения.

Введем еще следующие обозначения.

Для $ x,y \in \R^d $ положим $ xy =x \cdot y = (x_1 y_1, \ldots,
x_d y_d), $ а для $ x \in \R^d $ и $ A \subset \R^d $ определим
$$
x A = x \cdot A = \{xy: y \in A\}.
$$

Будем обозначать
$$
(x,y) = \sum_{j=1}^d x_j y_j, x, y \in \R^d.
$$

Обозначим через $ \R_+^d $ множество $ x \in \R^d, $ для которых $
x_j >0 $ при $ j=1,\ldots,d,$ и для $ a \in \R_+^d, x \in \R^d $
положим $ a^x = a_1^{x_1} \ldots a_d^{x_d}.$

При $ d \in \N $ определим множества
\begin{multline*}
I^d  =  \{x \in \R^d: 0 < x_j < 1,j=1,\ldots,d\},\\
\overline I^d  =  \{x \in \R^d: 0 \le x_j \le 1,j=1,\ldots,d\},\\
B^d  =  \{x \in \R^d: -1 \le x_j \le 1,j=1,\ldots,d\}.
\end{multline*}

Через $ \e $ будем обозначать вектор в $ \R^d, $ задаваемый
равенством $ \e =(1,\ldots,1). $

При $ d \in \N $ для $ x \in \R^d $ и $ J = \{j_1,\ldots,j_k\}
\subset \N: 1 \le j_1 < j_2 < \ldots < j_k \le d, $ через $ x^J $
обозначим вектор $ x^J = (x_{j_1},\ldots,x_{j_k}) \in \R^k, $ а
для множества $ A \subset \R^d $ положим $ A^J = \{x^J: x \in A\}. $

При $ d \in \N $ для $ x \in \R^d $ через $\s(x) $ обозначим
множество $ \s(x) = \{j =1,\ldots,d: x_j \ne 0\}, $ а для
множества $ J \subset \{1,\ldots,d\} $ через $ \chi_J $ обозначим
вектор из $ \R^d $ с компонентами
\begin{equation*}
(\chi_J)_j = \begin{cases} 1,  \text{для} j \in J; \\ 0,
 \text{для} j \in (\{1,\ldots,d\} \setminus J). \end{cases}
\end{equation*}
Будем также обозначать через $ \chi_A $ характеристическую функцию
множества $ A \subset \R^d. $

Для $ d \in \N, x,y \in \R^d $ будем писать $ x \le y (x < y), $
если для каждого $ j=1,\ldots,d $ выполняется неравенство $ x_j
\le y_j (x_j < y_j). $

Для банахова пространства $ X $ (над $ \C$) обозначим $ B(X) = \{x
\in X: \|x\|_X \le 1\}. $

Для банаховых пространств $ X,Y $ через $ \mathcal B(X,Y) $
обозначим банахово пространство, состоящее из непрерывных линейных
операторов $ T: X \mapsto Y, $ с нормой
$$
\|T\|_{\mathcal B(X,Y)} = \sup_{x \in B(X)} \|Tx\|_Y.
$$
Отметим, что если $ X=Y,$ то $ \mathcal B(X,Y) $ является
банаховой алгеброй. Отметим еще, что при $ Y =  \C $ пространство
$ \mathcal B(X, \C) $ обозначается также $ X^*. $

В завершении этого пункта приведем сведения о кратных
рядах, которыми будем пользоваться в дальнейшем.

При $ d \in \N $ для $ y \in \R^d $ обозначим
\begin{eqnarray*}
\mx(y)   =  \max_{j=1,\ldots,d} y_j,
\end{eqnarray*}
$$
\mn(y) = \min_{j=1,\ldots,d} y_j
$$
и для банахова пространства $ X, $ вектора $ x \in X $ и семейства
$ \{x_\kappa \in X, \kappa \in \Z_+^d\} $ будем писать $ x =
\lim_{ \mn(\kappa) \to \infty} x_\kappa, $ если для любого $
\epsilon >0 $ существует $ n_0 \in \N $ такое, что для любого $
\kappa \in \Z_+^d, $ для которого $ \mn(\kappa) > n_0, $
справедливо неравенство $ \|x -x_\kappa\|_X  < \epsilon. $

Пусть $ X $ -- банахово пространство (над $ \C $), $ d \in \N $ и
$ \{ x_\kappa \in X: \kappa \in \Z_+^d\} $ -- семейство векторов.
Тогда под суммой ряда $ \sum_{\kappa \in \Z_+^d} x_\kappa $ будем
понимать вектор $ x \in X, $ для которого выполняется равенство $
x = \lim_{\mn(k) \to \infty} \sum_{\kappa \in \Z_+^d(k)} x_\kappa. $

При $ d \in \N $ через $ \Upsilon^d $ обозначим множество
$$
\Upsilon^d = \{ \epsilon \in \Z^d: \epsilon_j \in \{0,1\},
j=1,\ldots,d\}.
$$

Имеет место

   Предложение  1.1.1

Пусть $ X $ -- банахово пространство, а вектор $ x \in X $ и
семейство $ \{x_\kappa \in X: \kappa \in \Z_+^d\} $ таковы, что $
x = \lim_{ \mn(\kappa) \to \infty} x_\kappa, $ Тогда для семейства $
\{ \mathcal X_\kappa \in X, \kappa \in \Z_+^d \}, $ определяемого
равенством
$$
\mathcal X_\kappa = \sum_{\epsilon \in \Upsilon^d: \s(\epsilon)
\subset \s(\kappa)} (-\e)^\epsilon x_{\kappa -\epsilon}, \kappa
\in \Z_+^d,
$$
справедливо равенство
$$
x = \sum_{\kappa \in \Z_+^d} \mathcal X_\kappa.
$$

Предложение является следствием того, что при $ k \in \Z_+^d $
выполняется равенство
\begin{equation*} \tag{1.1.1}
\sum|{\kappa \in \Z_+^d(k)}
\mathcal X_\kappa = x_k \text{(см. [8])}.
\end{equation*}

Замечание.

Легко заметить, что для любого семейства чисел $ \{x_\kappa \in
\R: x_\kappa \ge 0, \kappa \in \Z_+^d\} , $ если ряд $
\sum_{\kappa \in \Z_+^d} x_\kappa $ сходится, т.е. существует
предел $ \lim\limits_{\mn(k) \to \infty} \sum_{\kappa \in
\Z_+^d(k)} x_\kappa, $ то для любой последовательности подмножеств
$ \{Z_n \subset \Z_+^d, n \in \Z_+\}, $
 таких, что $ \card Z_n < \infty,
Z_n \subset Z_{n+1},  n \in \Z_+, $
и $ \cup_{ n \in \Z_+} Z_n = \Z_+^d, $
справедливо равенство
$ \sum_{\kappa \in \Z_+^d} x_\kappa =
\lim_{ n \to \infty} \sum_{\kappa \in Z_n} x_\kappa. $
Отсюда несложно понять, что если для семейства векторов
$ \{x_\kappa \in X, \kappa \in \Z_+^d\}  $ банахова пространства $ X $
ряд $ \sum_{\kappa \in \Z_+^d} \| x_\kappa \|_X $ сходится,
то для любой последовательности подмножеств
$ \{Z_n \subset \Z_+^d, n \in \Z_+\}, $
 таких, что $ \card Z_n < \infty,
Z_n \subset Z_{n+1},  n \in \Z_+, $
и $ \cup_{ n \in \Z_+} Z_n = \Z_+^d, $
в $ X $ соблюдается равенство
$ \sum_{\kappa \in \Z_+^d} x_\kappa =
\lim_{ n \to \infty} \sum_{\kappa \in Z_n} x_\kappa. $

\bigskip

1.2. В этом пункте изучаются некоторые свойства операторов проектирования,
рассматриваемых в работе.

Прежде всего отметим, что если $ H:  X \mapsto Y $ -- изоморфизм линейного
пространства $ X $ на линейное пространство $ Y $ и
$ P: Y \mapsto Y $ -- проекционный оператор, то $ H^{-1} P H $ --
проектор, действующий в $ X $ и
\begin{equation*} \tag{1.2.1}
\Im H^{-1} P H = H^{-1} \Im P,
\end{equation*}
\begin{equation*} \tag{1.2.2}
\Ker H^{-1} P H = H^{-1} \Ker P.
\end{equation*}

Лемма 1.2.1

Пусть $ P_0, P_1 $ -- операторы проектирования, действующеие в
линейном пространстве $ X, $ такие, что
\begin{equation*} \tag{1.2.3}
\Im P_0 \subset \Im P_1,
\Ker P_1 \subset \Ker P_0.
\end{equation*}
Тогда $ P_1 -P_0 $ является оператором
проектирования, причем,
\begin{equation*} \tag{1.2.4}
\Im (P_1 -P_0) = \Im P_1 \cap \Ker P_0;
\Ker (P_1 -P_0) = \Im P_0 +\Ker P_1.
\end{equation*}

Доказательство.

Действительно, с учетом (1.2.3), а также того факта, что
$ \Im (E -P_1) = \Ker P_1, $ имеем
\begin{multline*}
(P_1 -P_0)^2 = P_1^2 -P_1 P_0 -P_0 P_1 +P_0^2 = \\
P_1 -P_0 -P_0 (P_1 -E +E) +P_0 = P_1 -P_0.
\end{multline*}

Далее, ввиду (1.2.3)
$$
\Im (P_1 -P_0) \subset \Im P_1
$$
и
$$
P_0 (P_1 -P_0) =  P_0 (P_1 -E +E) -P_0^2 = 0 +P_0 -P_0 =0,
$$
т.е. $ \Im (P_1 -P_0) \subset \Ker P_0, $
а, следовательно, $ \Im (P_1 -P_0) \subset \Im P_1 \cap \Ker P_0. $
С другой стороны для $ x \in \Im P_1 \cap \Ker P_0 $ выполняется соотношение
$ x = P_1 x = P_1x -0 = P_1 x -P_0 x = (P_1 -P_0) x \in \Im (P_1 -P_0). $
Тем самым установлена справедливость первого равенства в (1.2.4).

Для проверки второго равенства в (1.2.4), учитывая, что $ P_1 -P_0 $ --
проектор, имеем
\begin{multline*}
\Ker (P_1 -P_0) = \Im (E -(P_1 -P_0)) = \Im((E -P_1) +P_0)\\
 \subset \Im (E -P_1) +\Im P_0 = \Ker P_1 +\Im P_0.
\end{multline*}
Наконец, для $ x = y +z, $ где $ y \in \Ker P_1, z \in \Im P_0, $
в силу (1.2.3) соблюдается равенство
$ (P_1 -P_0) x = P_1 y -P_0 y +P_1 z -P_0 z = 0 -0 +z -z = 0, $
и, значит, $ \Ker P_1 +\Im P_0 \subset \Ker (P_1 -P_0), $ что завершает
вывод второго равенства в (1.2.4). $ \square $

При $ \delta \in \R_+ $ обозначим через $ h_\delta $ отображение, которое каждой
(комплекснозначной) функции $ f  $, заданной на некотором множестве
$ S  \subset \R, $ ставит в соответствие функцию $ h_\delta f, $ определяемую
на множестве $ \{ x \in \R: \delta x \in S\} $ равенством
$ (h_\delta f)(x) = f(\delta x). $
Ясно, что отображение $ h_\delta $
является биекцией на себя  множества всех функций с  областью  определения
в $ \R. $ При этом,  для $ \delta, \sigma \in \R_+ $ имеет место равенство
$ h_{\delta \sigma} = h_\delta h_\sigma $ и
обратное отображение $ h_\delta^{-1} $  задается
равенством
$ h_\delta^{-1}  = h_{\delta^{-1}}.$

Отметим, что при $ 1 \le p  \le \infty $ отображение $ h_\delta, $
суженное на $ L_p(\R), $ является линейным гомеоморфизмом
пространства $ L_p(\R) $ на себя и для $ f \in L_p(\R) $ выполняется равенство
\begin{equation*}
\| h_\delta f\|_{L_p(\R)} =
\delta^{-p^{-1}} \|f\|_{L_p(\R)},
\end{equation*}
а, следовательно, для $ f \in L_p(\R) $ соблюдается равенство
\begin{equation*}
\|h_\delta^{-1}f\|_{L_p(\R)} =
\delta^{p^{-1}} \|f\|_{L_p(\R)}.
\end{equation*}

Кроме того, при $ \delta \in \R_+, 1 \le p \le \infty $ для $ f \in L_p(\R),
g \in L_{p^\prime}(\R) (p^\prime = p/(p -1)) $  справедливо соотношение
$$
\int_\R (h_\delta f) \overline {(h_\delta g)} dx = \int_\R f(\delta x)
\overline {g(\delta x)} dx =
\delta^{-1} \int_\R f(x) \overline {g(x)} dx
$$
а, следовательно,
\begin{equation*} \tag{1.2.5}
\int_\R (h_\delta f) \overline {(h_\delta g)} dx = 0 \iff
\int_\R f \overline g dx =0.
\end{equation*}

Для семейства $ \{ x_\nu \in X, \nu \in \Z \} $ векторов банахова пространства
$ X $ под суммой $ \sum_{\nu \in \Z} x_\nu $ будем понимать
$ \lim_{m, n \in \Z_+, \min(m, n) \to \infty} \sum_{\nu =-m}^n x_\nu. $

Через $ l_2, $ как обычно, обозначается гильбертово пространство векторов
$ x = \{x_\nu \in \C, \nu \in \Z\}, $ для которых ряд
$ \sum_{\nu \in \Z} | x_\nu |^2 $ сходится, со скалярным произведением
$$
\langle x, y \rangle = \sum_{\nu \in \Z } x_\nu \overline {y_\nu}, x, y \in l_2.
$$

Пусть $ X_0 $ -- замкнутое линейное подпространство в $ L_2(\R), $
для которого существует функция $ \phi \in X_0 $
такая, что система функций
$ \{ \phi(\cdot -\nu), \nu \in \Z \} $ является базисом Рисса в $ X_0, $
т.е. замыкание в $ L_2(\R) $ линейной оболочки этой системы
\begin{equation*} \tag{1.2.6}
\close_{L_2(\R)} (\span \{\phi(\cdot -\nu), \nu \in \Z \}) = X_0
\end{equation*}
и существуют константы $ A, B \in \R_+, $
обладающие тем свойством, что для любого семейства
$ \{ c_\nu, \nu \in \Z \} \in l_2 $ соблюдаются неравенства
\begin{equation*} \tag{1.2.7}
A \| \{c_\nu \} \|_{l_2} \le
\| \sum_{\nu \in \Z} c_\nu \phi(\cdot -\nu) \|_{L_2(\R)} \le
B \| \{c_\nu \} \|_{l_2} \text{(см. [6], [7]).}
\end{equation*}
Условия (1.2.6), (1.2.7) эквивалентны тому, что отображение
$ A­_\phi: l_2 \mapsto X_0, $ которое каждому
вектору $ \{c_\nu, \nu \in \Z \} \in l_2 $ ставит в соответствие
функцию $ \sum_{\nu \in \Z } c_\nu \phi(\cdot -\nu) \in X_0 \cap L_2(\R) $
является линейным гомеоморфизмом пространства $ l_2 $ на
подпространство $ X_0 \cap L_2(\R). $

Для функции $ \phi \in X_0, $ удовлетворяющей условиям (1.2.6), (1.2.7),
найдем двойственную ей функцию $ \phi^* \in X_0$ (см. гл. 1 из [6]),
т.е. такую, что при $ \nu, \mu \in \Z $ соблюдаются равенства
\begin{equation*} \tag{1.2.8}
\int_\R \phi(x -\nu) \overline {\phi^*(x -\mu)} dx = \begin{cases} 0,
\text{ при } \nu \in \Z \setminus \{\mu\}; \\
1, \text{при } \nu = \mu.
\end{cases}
\end{equation*}

При этом семейство функций $ \{\phi^*(\cdot -\nu), \nu \in \Z \} $ также является базисом Рисса в $ X_0 $ и
для каждого $ f \in X_0 $ имеет место представление
\begin{equation*} \tag{1.2.9}
f = \sum_{\nu \in \Z}  c_\nu \phi(\cdot -\nu),
\end{equation*}
где семейство $ \{ c_\nu, \nu \in \Z \} \in l_2, $ а для каждого
$ \nu \in \Z $ выполняется равенство
\begin{equation*} \tag{1.2.10}
c_\nu = \int_\R f(x) \overline {\phi^*(x -\nu)} dx.
\end{equation*}

Обозначим через $ U_0 $ оператор ортогонального проектирования
пространства $ L_2(\R) $ на подпространство $ X_0, $
т.е. $ U_0: L_2(\R) \mapsto L_2(\R) $ -- непрерывный проектор, у которого
\begin{equation*} \tag{1.2.11}
\Im U_0 = X_0,
\end{equation*}
\begin{equation*} \tag{1.2.12}
\Ker U_0 = \Im (E -U_0) = \{ f \in L_2(\R): \int_\R f \overline g dx =0
\forall g \in X_0 \},
\end{equation*}
где $ E $ -- тождественный оператор в $ L_2(\R). $

Из (1.2.9), (1.2.10), (1.2.11), (1.2.12) вытекает, что для $ f \in L_2(\R) $
справедливо равенство
$$
U_0 f = \sum_{ \nu \in \Z} (\int_\R (U_0 f)(x) \overline {\phi^*(x -\nu)} dx)
\phi(\cdot -\nu),
$$
причем,
\begin{multline*}
\int_\R (U_0 f)(x) \overline {\phi^*(x -\nu)} dx =
\int_\R ((U_0 f)(x) -f(x) +f(x)) \overline {\phi^*(x -\nu)} dx = \\
\int_\R ((U_0 f)(x) -f(x)) \overline {\phi^*(x -\nu)} dx +
\int_\R f(x) \overline {\phi^*(x -\nu)} dx =
\int_\R f(x) \overline {\phi^*(x -\nu)} dx,
\end{multline*}
и, значит,
\begin{equation*} \tag{1.2.13}
U_0 f = \sum_{ \nu \in \Z} (\int_\R f(x) \overline {\phi^*(x -\nu)} dx)
\phi(\cdot -\nu).
\end{equation*}

Теперь при $ k \in \N $ положим
\begin{equation*} \tag{1.2.14}
X_k = h_2^k ( X_0 ) = h_2 (h_2^{k -1} X_0) = h_2 ( X_{k -1}),
\end{equation*}
\begin{equation*} \tag{1.2.15}
U_k = h_2^k U_0 (h_2^k)^{-1} =
h_{2^k} U_0 (h_{2^k})^{-1} = h_{2^k} U_0 h_{2^{-k}}.
\end{equation*}

Сопоставляя (1.2.1), (1.2.2), (1.2.14), (1.2.15), (1.2.11), (1.2.12), (1.2.5), (1.2.6),
заключаем, что $ U_k: L_2(\R) \mapsto L_2(\R) $ является проектором,
у которого
\begin{multline*} \tag{1.2.16}
\Im U_k = X_k = h_{2^k} (\close_{L_2(\R)} (\span \{\phi(\cdot -\nu),
\nu \in \Z \})) = \\
\close_{L_2(\R)} (h_{2^k} (\span \{\phi(\cdot -\nu),
\nu \in \Z \})) = \\
\close_{L_2(\R)} (\span \{ h_{2^k} (\phi(\cdot -\nu)),
\nu \in \Z \}) = \\
\close_{L_2(\R)} (\span \{ \phi(2^k \cdot -\nu)),
\nu \in \Z \}),
\end{multline*}
\begin{multline*} \tag{1.2.17}
\Ker U_k = h_{2^k} (\{ f \in L_2(\R): \int_\R f \overline g dx =0
\forall g \in X_0 \}) = \\
\{ h_{2^k}(f): f \in L_2(\R), \int_\R h_{2^k}(f)
\overline { h_{2^k}(g)} dx =0
\forall g \in X_0 \}) = \\
\{ F \in L_2(\R): \int_\R F \overline G dx =0 \forall G \in X_k \}.
\end{multline*}

Принимая во внимание (1.2.16), (1.2.17) и учитывая, что для $ f \in L_2(\R) $
справедливо соотношение
\begin{multline*}
\int_\R f(x) \overline {g(x)} dx =0 \forall g \in
\close_{L_2(\R)} (\span \{ \phi(2^k \cdot -\nu), \nu \in \Z \}) \iff \\
\int_\R f(x) \overline {g(x)} dx =0 \forall g \in
\span \{ \phi(2^k \cdot -\nu), \nu \in \Z \} \iff \\
\int_\R f(x) \overline {\phi(2^k x -\nu)} dx =0 \forall \nu \in \Z,
\end{multline*}
находим, что при  $ k \in \Z_+ $ верно равенство
\begin{equation*}
\Ker U_k = \{ f \in L_2(\R): \int_\R f(x) \overline {\phi(2^k x -\nu)} dx =0
\forall \nu \in \Z \}.
\end{equation*}

Отметим, что ввиду (1.2.15), (1.2.13) для $ f \in L_2(\R) $
при $ k \in \Z_+ $ выполняется равенство
\begin{multline*} \tag{1.2.18}
U_k f = h_{2^k} (\sum_{ \nu \in \Z} (\int_\R ( h_{2^{-k}} f)(x)
\overline {\phi^*(x -\nu)} dx) \phi(\cdot -\nu)) = \\
\sum_{ \nu \in \Z} (\int_\R f(2^{-k} x)
\overline {\phi^*(x -\nu)} dx) \phi(2^k \cdot -\nu) = \\
\sum_{ \nu \in \Z} 2^k (\int_\R f(x)
\overline {\phi^*(2^k x -\nu)} dx) \phi(2^k \cdot -\nu).
\end{multline*}

При $ d \in \N $ для $ t \in \R^d $ через $ 2^t $ будем обозначать
вектор $ 2^t = (2^{t_1}, \ldots, 2^{t_d}). $

Для $ d \in \N, \kappa \in \Z^d, \nu \in \Z^d $ обозначим через
$ Q_{\kappa, \nu}^d = 2^{-\kappa} \nu +2^{-\kappa} I^d. $

Предложение 1.2.2

Пусть $ d \in \N, $ функции $ \phi, \phi^* \in L_\infty(\R^d) $ имеют
компактные носители и при $ \kappa \in \Z_+^d $ линейный оператор
$ E_\kappa = E_\kappa^{\phi,\phi^*}: L_1^{\loc(\R^d)} \mapsto L_1^{\loc(\R^d)} $
определяется равенством
\begin{equation*} \tag{1.2.19}
(E_\kappa f) (x) = \sum_{\nu \in \Z^d} 2^{(\kappa, \e)}
(\int_{\R^d} f(y) \overline {\phi^*(2^\kappa y -\nu)} dy)
\phi(2^\kappa x -\nu), x \in \R^d, f \in L_1^{\loc(\R^d)},
\end{equation*}
где суммирование производится поточечно.
Тогда существует константа $ c_1(d, \phi, \phi^*) >0 $ такая, что
при $ 1 \le p \le \infty $ для $ f \in L_p(\R^d), \kappa \in \Z_+^d $
справедлива оценка
\begin{equation*} \tag{1.2.20}
\| E_\kappa f \|_{L_p(\R^d)} \le c_1 \| f \|_{L_p(\R^d)},
\end{equation*}
и для любой функции $ f \in L_1^{\loc(\R^d)}, $ имеющей компактный носитель,
существует компакт $ K \subset \R^d $ такой, что при любом
$ \kappa \in \Z_+^d $ носитель
\begin{equation*} \tag{1.2.21}
\supp (E_\kappa f) \subset K.
\end{equation*}

Доказательство.

При $ \kappa, \nu \in \Z^d $ обозначим
$ \phi_{\kappa, \nu}(\cdot) = \phi(2^\kappa \cdot -\nu),
\phi^*_{\kappa, \nu}(\cdot) = \phi^*(2^\kappa \cdot -\nu). $
Заметим, что существует константа $ c_2(d, \phi) > 0 $ такая, что
при $ \kappa, n \in \Z^d $ справедлива оценка
\begin{equation*} \tag{1.2.22}
\card \{\nu \in \Z^d: \supp \phi_{\kappa, \nu} \cap Q_{\kappa, n}^d \ne
\emptyset \} \le c_2.
\end{equation*}

В самом деле, учитывая, что $ \supp \phi_{\kappa, \nu} =
2^{-\kappa} \nu +2^{-\kappa} \supp \phi, $ имеем
\begin{multline*}
\{\nu \in \Z^d: \supp \phi_{\kappa, \nu} \cap Q_{\kappa, n}^d \ne \emptyset \} =\\
\{ \nu \in \Z^d: (2^{-\kappa} \nu +2^{-\kappa} \supp \phi) \cap
(2^{-\kappa} n  +2^{-\kappa} I^d) \ne \emptyset \} = \\
\{\nu \in \Z^d: (\nu + \supp \phi) \cap (n +I^d) \ne \emptyset \} = \\
n +\{\nu^\prime \in \Z^d: (\nu^\prime +n + \supp \phi) \cap (n +I^d)
\ne \emptyset \} = \\
n +\{\nu^\prime \in \Z^d: (\nu^\prime +\supp \phi) \cap I^d \ne \emptyset \},
\end{multline*}
а, следовательно,
$$
\card \{\nu \in \Z^d: \supp \phi_{\kappa, \nu} \cap Q_{\kappa, n}^d
\ne \emptyset \} =
\card \{\nu^\prime \in \Z^d: (\nu^\prime +\supp \phi) \cap I^d \ne \emptyset \}.
$$
Отсюда ввиду ограниченности $ \supp \phi $ вытекает (1.2.22).

Отметим, что (1.2.22) сохраняет свою силу, если $ \phi $ заменить
на $ \phi^*. $

Из (1.2.22) следует, что на каждом кубе $ Q_{\kappa,n}^d, n \in \Z^d, $
число отличных от 0 слагаемых в правой части (1.2.19) конечно.

Пусть $ 1 \le p < \infty, f \in L_p(\R^d). $ Тогда при $ \kappa \in \Z_+^d $
выводим
\begin{multline*} \tag{1.2.23}
\| E_\kappa f \|_{L_p(\R^d)}^p =
\|\sum_{\nu \in \Z^d}2^{(\kappa, \e)} (\int_{\R^d} f(y) \overline {\phi^*_{\kappa, \nu}(y)} dy)
\phi_{\kappa, \nu} \|_{L_p(\R^d)}^p = \\
\int_{\R^d} | \sum_{\nu \in \Z^d}
2^{(\kappa, \e)} (\int_{\R^d} f(y) \overline {\phi^*_{\kappa, \nu}(y)} dy)
\phi_{\kappa, \nu}|^p dx = \\
\sum_{n \in \Z^d} \int_{Q_{\kappa, n}^d} | \sum_{\nu \in \Z^d}
2^{(\kappa, \e)} (\int_{\R^d} f(y) \overline {\phi^*_{\kappa, \nu}(y)} dy)
\phi_{\kappa, \nu}|^p dx = \\
\sum_{n \in \Z^d} \int_{Q_{\kappa, n}^d} | \sum_{\nu \in \Z^d:
\supp \phi_{\kappa, \nu} \cap Q_{\kappa, n}^d \ne \emptyset}
2^{(\kappa, \e)} (\int_{\R^d} f(y) \overline {\phi^*_{\kappa, \nu}(y)} dy)
\phi_{\kappa, \nu}|^p dx \le \\
\sum_{n \in \Z^d} \int_{Q_{\kappa, n}^d} (\sum_{\nu \in \Z^d:
\supp \phi_{\kappa, \nu} \cap Q_{\kappa, n}^d \ne \emptyset}
2^{(\kappa, \e)}
|\int_{\R^d} f(y) \overline {\phi^*_{\kappa, \nu}(y)} dy|
\|\phi\|_{L_\infty(\R^d)})^p dx = \\
\sum_{n \in \Z^d} (\mes Q_{\kappa, n}^d) 2^{(\kappa, \e) p}
\|\phi\|_{L_\infty(\R^d)}^p
(\sum_{\nu \in \Z^d: \supp \phi_{\kappa, \nu} \cap Q_{\kappa, n}^d
\ne \emptyset}
|\int_{\R^d} f(y) \overline {\phi^*_{\kappa, \nu}(y)} dy| )^p \le \\
2^{(\kappa, \e) (p -1)} \|\phi\|_{L_\infty(\R^d)}^p
\sum_{n \in \Z^d}
(\card \{\nu \in \Z^d: \supp \phi_{\kappa, \nu} \cap Q_{\kappa, n}^d \ne
\emptyset\})^{p/p^\prime} \cdot \\
(\sum_{\nu \in \Z^d: \supp \phi_{\kappa, \nu} \cap Q_{\kappa, n}^d \ne \emptyset}
|\int_{\R^d} f(y) \overline {\phi^*_{\kappa, \nu}(y)} dy|^p) \le \\
2^{(\kappa, \e) (p -1)} \|\phi\|_{L_\infty(\R^d)}^p
\sum_{n \in \Z^d} c_2^{p/p^\prime}
(\sum_{\nu \in \Z^d: \supp \phi_{\kappa, \nu} \cap Q_{\kappa, n}^d \ne \emptyset}
|\int_{\supp \phi^*_{\kappa,\nu}} f(y)
\overline {\phi^*_{\kappa, \nu}(y)} dy|^p) \le \\
c_2^{p/p^\prime} 2^{(\kappa, \e) (p -1)} \|\phi\|_{L_\infty(\R^d)}^p
\sum_{n \in \Z^d}
(\sum_{\nu \in \Z^d: \supp \phi_{\kappa, \nu} \cap Q_{\kappa, n}^d \ne \emptyset}
(\int_{\supp \phi^*_{\kappa,\nu}} |f|^p dy) \cdot \\
(\int_{2^{-\kappa \nu} +2^{-\kappa} \supp \phi^*} |\phi^*_{\kappa, \nu}|^{p^\prime} dy)^{p/p^\prime}) =
c_2^{p/p^\prime} 2^{(\kappa, \e) (p -1)} \|\phi\|_{L_\infty(\R^d)}^p  \cdot \\
\sum_{n \in \Z^d}
(\sum_{\nu \in \Z^d: \supp \phi_{\kappa, \nu} \cap Q_{\kappa, n}^d \ne \emptyset}
(\int_{\R^d} \chi_{\supp \phi^*_{\kappa,\nu}} |f|^p dy)
2^{-(\kappa, \e) p/p^\prime}
(\int_{\supp \phi^*} |\phi^*|^{p^\prime} dz)^{p/p^\prime}) \le\\
c_2^{p/p^\prime} \|\phi\|_{L_\infty(\R^d)}^{p}
\|\phi^*\|_{L_\infty(\R^d)}^{p}
(\mes \supp \phi^* )^{p/p^\prime}
\sum_{n \in \Z^d}
\sum_{\nu \in \Z^d: \supp \phi_{\kappa, \nu} \cap Q_{\kappa, n}^d \ne \emptyset}
\int_{\R^d} \chi_{\supp \phi^*_{\kappa,\nu}} |f|^p dy \le \\
(c_3(d, \phi, \phi^*))^p \int_{\R^d} (\sum_{n \in \Z^d}
\sum_{\nu \in \Z^d: \supp \phi_{\kappa, \nu} \cap Q_{\kappa, n}^d \ne \emptyset}
\chi_{\supp \phi^*_{\kappa,\nu}}) |f|^p dy = \\
c_3^p \sum_{n^\prime \in \Z^d}
\int_{ Q_{\kappa, n^\prime}^d} (\sum_{n \in \Z^d}
\sum_{\nu \in \Z^d: \supp \phi_{\kappa, \nu} \cap Q_{\kappa, n}^d \ne \emptyset}
\chi_{\supp \phi^*_{\kappa,\nu}}) |f|^p dy = \\
c_3^p \sum_{n^\prime \in \Z^d}\int_{ Q_{\kappa, n^\prime}^d}
(\sum_{\substack{n \in \Z^d: \exists \nu \in \Z^d: \supp \phi_{\kappa, \nu} \cap Q_{\kappa, n}^d \ne \emptyset,\\
\supp \phi^*_{\kappa, \nu} \cap Q_{\kappa, n^\prime}^d \ne \emptyset}}
\sum_{\nu \in \Z^d: \supp \phi_{\kappa, \nu} \cap Q_{\kappa, n}^d \ne \emptyset}
\chi_{\supp \phi^*_{\kappa,\nu}}) |f|^p dy.
\end{multline*}

Заметим, что если для $ n, n^\prime \in \Z^d $ существует $ \nu \in \Z^d, $
для которого $ \supp \phi_{\kappa, \nu} \cap Q_{\kappa, n}^d \ne \emptyset $ и,
$ \supp \phi^*_{\kappa, \nu} \cap Q_{\kappa, n^\prime}^d \ne \emptyset, $
то, как видно из вывода (1.2.22), имеют место соотношения
$ (\nu -n +\supp \phi) \cap I^d \ne \emptyset $ и
$ (\nu -n^\prime +\supp \phi^*) \cap I^d \ne \emptyset, $
а, следовательно, существуют точки $ \xi \in \supp \phi, \eta \in \supp \phi^*  $
такие, что $ \nu -n +\xi \in I^d, \nu -n^\prime +\eta \in I^d. $
Поэтому для $ j =1, \ldots, d $ получаем
\begin{multline*}
n_j^\prime -n_j = \nu_j -n_j +\xi_j -\xi_j -(\nu_j -n_j^\prime) -\eta_j +\eta_j = \\
(\nu_j -n_j +\xi_j) -(\nu_j -n_j^\prime +\eta_j) +(\eta_j -\xi_j) > \\
0 -1 -|\eta_j -\xi_j| \ge -1 -\sup_{\eta \in \supp \phi^*, \xi \in \supp \phi}
|\eta_j -\xi_j|.
\end{multline*}
Точно так же для $ n, n^\prime \in \Z^d, $ для которых существует
$ \nu \in \Z^d: \supp \phi_{\kappa, \nu} \cap Q_{\kappa, n}^d \ne
\emptyset, \supp \phi_{\kappa,\nu}^* \cap Q_{\kappa, n^\prime}^d
\ne \emptyset, $ находим, что
$$
n_j^\prime -n_j \le 1 +\sup_{\eta \in \supp \phi^*, \xi \in \supp
\phi} |\eta_j -\xi_j|,
$$
и, значит,
$$
| n_j -n_j^\prime | \le 1 +\sup_{\eta \in \supp
\phi^*, \xi \in \supp \phi} | \eta_j -\xi_j |, j =1, \ldots, d.
$$

Учитывая это замечание, видим, что существует константа $ c_4(d,
\phi, \phi^*)
> 0 $ такая, что для $ n^\prime \in \Z^d $ справедлива оценка
\begin{equation*} \tag{1.2.24}
\card \{n \in \Z^d: \exists \nu \in \Z^d:
\supp \phi_{\kappa, \nu} \cap Q_{\kappa, n}^d \ne \emptyset,
\supp \phi^*_{\kappa, \nu} \cap Q_{\kappa, n^\prime}^d \ne \emptyset \}
\le c_4.
\end{equation*}
Принимая во внимание (1.2.24) и (1.2.22), заключаем, что при
$ n^\prime \in \Z^d $ выполняется неравенство
\begin{multline*} \tag{1.2.25}
\sum_{\substack{n \in \Z^d: \exists \nu \in \Z^d: \supp \phi_{\kappa, \nu} \cap
Q_{\kappa, n}^d \ne \emptyset,\\
\supp \phi^*_{\kappa, \nu} \cap Q_{\kappa, n^\prime}^d \ne \emptyset}}
\sum_{\nu \in \Z^d: \supp \phi_{\kappa, \nu} \cap Q_{\kappa, n}^d \ne \emptyset}
\chi_{\supp \phi^*_{\kappa,\nu}} \le \\
\sum_{\substack{n \in \Z^d: \exists \nu \in \Z^d: \supp \phi_{\kappa, \nu} \cap
Q_{\kappa, n}^d \ne \emptyset,\\
\supp \phi^*_{\kappa, \nu} \cap Q_{\kappa, n^\prime}^d \ne \emptyset}} c_2 \le
c_2 c_4.
\end{multline*}

Подставляя (1.2.25) в (1.2.23), приходим к неравенству
$$
\| E_\kappa f \|_{L_p(\R^d)}^p \le
c_3^p \sum_{n^\prime \in \Z^d}
\int_{ Q_{\kappa, n^\prime}^d} c_5 |f|^p dy \le
(c_1(d, \phi, \phi^*))^p \int_{\R^d} |f|^p dy,
$$
что влечет (1.2.20).

Аналогично проводится вывод (1.2.20) при  $ p = \infty. $

Для проверки (1.2.21) отметим, что при $ \kappa \in \Z_+^d $ для $ \nu \in \Z^d $
имеет место равенство
\begin{multline*}
\int_{\R^d} f(y) \overline {\phi^*(2^\kappa y -\nu)} dy =
\int_{\supp f \cap \supp \phi^*(2^\kappa \cdot -\nu)} f(y)
\overline {\phi^*(2^\kappa y -\nu)} dy = \\
\int_{\supp f \cap (2^{-\kappa} \nu +2^{-\kappa} \supp \phi^*)} f(y)
\overline {\phi^*(2^\kappa y -\nu)} dy,
\end{multline*}
и, следовательно,
\begin{multline*}
(E_\kappa f) (\cdot) = \sum_{\nu \in \Z^d: \supp f \cap (2^{-\kappa} \nu
+2^{-\kappa} \supp \phi^*) \ne \emptyset } 2^{(\kappa, \e)} \cdot \\
(\int_{\R^d} f(y) \overline {\phi^*(2^\kappa y -\nu)} dy) \phi(2^\kappa \cdot -\nu),
\end{multline*}
а, значит,
\begin{multline*}
\supp (E_\kappa f) \subset \cup_{\nu \in \Z^d: \supp f \cap (2^{-\kappa} \nu
+2^{-\kappa} \supp \phi^*) \ne \emptyset } \supp \phi(2^\kappa \cdot -\nu) = \\
\cup_{\nu \in \Z^d: \supp f \cap (2^{-\kappa} \nu +2^{-\kappa} \supp \phi^*)
\ne \emptyset } (2^{-\kappa} \nu +2^{-\kappa} \supp \phi).
\end{multline*}
Отсюда видим, что для $ x \in \supp (E_\kappa f) $ существует
$ \nu \in \Z^d, $ для которого $ x = 2^{-\kappa} \nu +2^{-\kappa} \xi, $
где $ \xi \in \supp \phi, $
а для некоторого $ \eta \in \supp \phi^* $ точка $ y = (2^{-\kappa} \nu
+2^{-\kappa} \eta) \in \supp f, $ что влечет неравенства
\begin{multline*}
| x_j -y_j | = | 2^{-\kappa_j} \nu_j +2^{-\kappa_j} \xi_j -
(2^{-\kappa_j} \nu_j +2^{-\kappa_j} \eta_j) | = \\
2^{-\kappa_j} | \xi_j -\eta_j| \le
\sup_{\xi \in \supp \phi, \eta \in \supp \phi^*} | \xi_j -\eta_j |, j =1, \ldots,d,
\end{multline*}
откуда следует (1.2.21).$ \square $

Предложение 1.2.3

Пусть функция $ \phi \in L_\infty(\R) $ имеет компактный носитель
и выполяются условия (1.2.6), (1.2.7), а функция
$ \phi^* \in L_\infty(\R) \cap X_0 $ также имеет компактный носитель и
соблюдаются равенства (1.2.8). Тогда для операторов $ E_\kappa, U_\kappa,
\kappa \in \Z_+ $ (см. (1.2.19), (1.2.15)) справедливы соотношения:

1) \begin{equation*} \tag{1.2.26}
E_\kappa f = U_\kappa f, f \in L_2(\R);
\end{equation*}

2) при $ 1 \le p < \infty $ оператор $ E_\kappa^p = E_\kappa \mid_{L_p(\R)} $
является оператором проектирования в $ L_p(\R) $, для которого
\begin{multline*} \tag{1.2.27}
\Im E_\kappa^p = \close_{L_p(\R)} (\span \{ \phi(2^\kappa \cdot -\nu),
\nu \in \Z \}) = \\
h_{2^\kappa} (\close_{L_p(\R)} (\span
\{ \phi(\cdot -\nu), \nu \in \Z \}));
\end{multline*}
\begin{multline*} \tag{1.2.28}
\Ker E_\kappa^p = \{ f \in L_p(\R): \int_\R f(x)
\overline {\phi^*(2^\kappa x -\nu)} dx =0 \forall \nu \in \Z \} = \\
\{ f \in L_p(\R): \int_\R f(x)
\overline {\phi(2^\kappa x -\nu)} dx =0 \forall \nu \in \Z \} = \\
\{ f \in L_p(\R): \int_\R f(x)
\overline {g(x)} dx =0 \forall g \in
\span \{\phi(2^\kappa \cdot -\nu),  \nu \in \Z\} \};
\end{multline*}

3) при $ 1 < P < \infty $ для $ f \in L_p(\R),
g \in L_{p^\prime}(\R) $ имеют место
равенства
\begin{equation*} \tag{1.2.29}
\int_{\R} (E_\kappa f) \cdot \overline g dx = \int_{\R} f \cdot
\overline {(E_\kappa g) } dx
\end{equation*}
и
\begin{equation*} \tag{1.2.30}
\int_{\R} (\mathcal E_\kappa f) \cdot \overline g dx = \int_{\R} f \cdot
\overline {(\mathcal E_\kappa g) } dx,
\end{equation*}
где оператор $ \mathcal E_\kappa = \mathcal E_\kappa^{\phi, \phi^*}:
L_1^{\loc(\R)} \mapsto L_1^{\loc(\R)}, $ определяется равенством
$$
\mathcal E_\kappa = E_\kappa -E_{\kappa -1}, \kappa \in \Z_+,
E_{-1} =0.
$$

Доказательство.

Чтобы убедиться в справедливости (1.2.26), заметим, что для каждой
функции $ f \in L_2(\R), $ носитель которой компактен, в силу (1.2.19),
(1.2.18) выполняется равенство
\begin{multline*}
( E_\kappa f) (\cdot) = \sum_{\nu \in \Z: \supp f \cap
\supp \phi^* (2^{\kappa} \cdot -\nu) \ne \emptyset } 2^\kappa \cdot \\
(\int_\R f(y) \overline {\phi^*(2^\kappa y -\nu)} dy) \phi(2^\kappa \cdot -\nu) =
(U_\kappa f) (\cdot),
\end{multline*}
(ибо множество $ \{ \nu \in \Z: \supp f \cap
\supp \phi^* (2^{\kappa} \cdot -\nu) \ne \emptyset \} $ -- конечно).
Учитывая, что $ E_\kappa $ и $ U_\kappa $ -- непрерывные операторы
в $ L_2(\R), $ а множество функций $ f \in L_2(\R), $
имеющих компактные носители, плотно в $ L_2(\R), $
отсюда заключаем, что (1.2.26) справедливо для всех $ f \in L_2(\R). $

Проверим, что $ E_\kappa^p $ является проектором. Для этого, учитывая, что для
$ f \in L_2(\R) \cap L_p(\R), $ благодаря (1.2.20), имеет место
соотношение $ E_\kappa^p f = E_\kappa f \in L_2(\R) \cap L_p(\R), $ используя
(1.2.26) и тот факт, что $ U_\kappa $ -- оператор проектирования, выводим
\begin{equation*} \tag{1.2.31}
E_\kappa^p E_\kappa^p f =
E_\kappa E_\kappa f =
U_\kappa U_\kappa f = U_\kappa f = E_\kappa f = E_\kappa^p f.
\end{equation*}
Отсюда в силу непрерывности оператора $ E_\kappa^p: L_p(\R) \mapsto L_p(\R) $
(см. (1.2.20)) и плотности $ L_2(\R) \cap L_p(\R) $ в $ L_p(\R) $
при $ 1 \le p < \infty $ вытекает соблюдение (1.2.31) для $ f \in L_p(\R), $
т.е. $ E_\kappa^p $ -- проектор в $ L_p(\R). $

Установим справедливость (1.2.27). Выбирая для $ f \in L_p(\R) $
последовательность функций $ \{ f_n \in L_\infty(\R), n \in \N \}, $ имеющих
компактные носители, сходящуюся к $ f $  в $ L_p(\R), $ вследствие (1.2.19)
и непрерывности оператора $ E_\kappa^p: L_p(\R) \mapsto L_p(\R) $
получаем, что
$ E_\kappa^p f_n \in \span \{ \phi(2^\kappa \cdot -\nu), \nu \in \Z \},
n \in \N, $ и $ E_\kappa^p f_n \to E_\kappa^p f $ в $ L_p(\R) $
при $ n \to \infty, $ т.е. $ \Im E_\kappa^p \subset \close_{L_p(\R)}
(\span \{ \phi(2^\kappa \cdot -\nu), \nu \in \Z \}). $
Обратно, для $ f \in \close_{L_p(\R)} (\span \{ \phi(2^\kappa \cdot -\nu),
\nu \in \Z \}), $ беря послдовательность $ \{ f_n \in
\span \{ \phi(2^\kappa \cdot -\nu), \nu \in \Z \}, n \in \N \}, $
сходящуюся к $ f $ в $ L_p(\R), $ используя непрерывность
$ E_\kappa \mid_{L_p(\R)} $ в $ L_p(\R), $ (1.2.26), (1.2.16) и тот факт,
что $ U_\kappa $ суть проектор, получаем,
$$
E_\kappa^p f = \lim_{n \to \infty} E_\kappa^p f_n =
\lim_{n \to \infty} E_\kappa f_n =
\lim_{n \to \infty} U_\kappa f_n =
\lim_{n \to \infty} f_n = f,
$$
следовательно,
$$
\close_{L_p(\R)} (\span \{ \phi(2^\kappa \cdot -\nu),
\nu \in \Z \}) \subset \Im E_\kappa^p.
$$
Тем самым, установлена справедливость (1.2.27).

Перейдем к выводу (1.2.28). Сначала получим (1.2.28) при $ \kappa =0. $
Пусть $ f \in \Ker E_0^p. $ Тогда при $ \mu \in \Z $ ввиду (1.2.19), (1.2.8)
соблюдаются равенства
\begin{multline*}
0 = \int_\R (E_0^p f)(x) \overline {\phi^*(x -\mu)} dx =
\int_\R (\sum_{\nu \in \Z} (\int_\R f(y) \overline {\phi^*(y -\nu)} dy)
\phi(x -\nu)) \overline {\phi^*(x -\mu)} dx = \\
\sum_{\nu \in \Z} \int_\R (\int_\R f(y) \overline {\phi^*(y -\nu)} dy)
\phi(x -\nu) \overline {\phi^*(x -\mu)} dx = \\
\sum_{\nu \in \Z} (\int_\R f(y) \overline {\phi^*(y -\nu)} dy) \int_\R
\phi(x -\nu) \overline {\phi^*(x -\mu)} dx =
\int_\R f(y) \overline {\phi^*(y -\mu)} dy).
\end{multline*}

Если же для $ f \in L_p(\R) $ выполняются равенства
$ \int_\R f(y) \overline {\phi^*(y -\nu)} dy =0 , \nu \in \Z, $
то согласно (1.2.19) $ E_0^p f =0. $
Таким образом, первое равенство в (1.2.28) при $ \kappa = 0 $ установлено.

Теперь заметим, что в силу (1.2.9), (1.2.10) имеет место разложение
\begin{multline*} \tag{1.2.32}
\phi^* (\cdot -\mu) = \sum_{ \nu \in \Z} (\int_\R \phi^*(y -\mu)
\overline {\phi^*(y -\nu)} dy) \phi(\cdot -\nu) = \\
\sum_{ \nu \in \Z: \supp \phi^*( \cdot -\mu \cap \supp \phi^*(\cdot -\nu)
\ne \emptyset} (\int_\R \phi^*(y -\mu)
\overline {\phi^*(y -\nu)} dy)
\phi(\cdot -\nu) \in \\
\span \{ \phi(\cdot -\nu), \nu \in \Z \}.
\end{multline*}

Поскольку система функций $ \{ \phi^*(\cdot -\nu) , \nu \in \Z \} $ образует
двойственный системе $ \{ \phi(\cdot -\nu) , \nu \in \Z \} $
базис Рисса в $ X_0 $ (см. выше), то меняя ролями $ \phi $ и $ \phi^* $
в (1.2.32), находим, что
$$
\phi (\cdot -\mu) \in \span \{ \phi^*(\cdot -\nu), \nu \in \Z \}, \mu \in \Z,
$$
а, следовательно,
$$
\span \{ \phi(\cdot -\nu), \nu \in \Z \} =
\span \{ \phi^*(\cdot -\nu), \nu \in \Z \}.
$$

Объединяя сказанное, заключаем, что
\begin{multline*}
\Ker E_0^p = \{ f \in L_p(\R): \int_\R f(x)
\overline {\phi^*( x -\nu)} dx =0 \forall \nu \in \Z \} = \\
\{ f \in L_p(\R): \int_\R f(x)
\overline {g(x)} dx =0 \forall g \in \span \{ \phi^*(\cdot -\nu), \nu \in \Z \} = \\
\{ f \in L_p(\R): \int_\R f(x)
\overline {g(x)} dx =0 \forall g \in \span \{ \phi(\cdot -\nu), \nu \in \Z \} = \\
\{ f \in L_p(\R): \int_\R f(x)
\overline {\phi( x -\nu)} dx =0 \forall \nu \in \Z \},
\end{multline*}
что завершает вывод (1.2.28) при $ \kappa =0. $

Наконец, убедимся в спраедливости (1.2.28) при $ \kappa \in \N. $
Принимая во внимание, что $ E_\kappa^p =
h_{2^\kappa} E_0^p (h_{2^\kappa})^{-1}, $ на основании (1.2.2), (1.2.5)
и (1.2.28) при $ \kappa =0 $ получаем
\begin{multline*}
\Ker E_\kappa^p =
h_{2^\kappa} (\{ f \in L_p(\R): \int_\R f(x)
\overline {\phi( x -\nu)} dx =0 \forall \nu \in \Z \}) = \\
\{ h_{2^\kappa} f: f \in L_p(\R): \int_\R f(x)
\overline {\phi( x -\nu)} dx =0 \forall \nu \in \Z \} = \\
\{ h_{2^\kappa} f: f \in L_p(\R): \int_\R (h_{2^\kappa} f)(x)
\overline {(h_{2^\kappa} \phi(\cdot -\nu))(x)} dx =0 \forall \nu \in \Z \} = \\
\{ f \in L_p(\R): \int_\R f(x)
\overline { \phi(2^\kappa x -\nu)} dx =0 \forall \nu \in \Z \}.
\end{multline*}
Ввиду последнего равенства ясно, что и оставшиеся равенства в (1.2.28) также
имеют место.

Для получения (1.2.29), принимая во внимание, (1.2.26), (1.2.16), (1.2.17),
а также то обстоятельство, что $ \Im (E -U_\kappa) = \Ker U_\kappa, $
для $ f \in L_p(\R) \cap L_2(\R), g \in L_{p^\prime}(\R) \cap L_2(\R) $ имеем
\begin{multline*}
\int_{\R} (E_\kappa f) \cdot \overline g dx =
\int_{\R} (U_\kappa f) \cdot \overline g dx =
\int_{\R} (U_\kappa f) \cdot \overline {( U_\kappa g +g -U_\kappa g)} dx = \\
\int_{\R} (U_\kappa f) \cdot \overline { U_\kappa g } dx +
\int_{\R} (U_\kappa f) \cdot \overline {( g -U_\kappa g)} dx =
\int_{\R} (U_\kappa f) \cdot \overline { (U_\kappa g )} dx,
\end{multline*}
и
\begin{multline*}
\int_{\R} f \cdot \overline {(E_\kappa g) } dx =
\int_{\R} f \cdot \overline {(U_\kappa g) } dx =
\int_{\R} (U_\kappa f +f -U_\kappa f) \cdot \overline {(U_\kappa g) } dx = \\
\int_{\R} (U_\kappa f ) \cdot \overline {(U_\kappa g) } dx +
\int_{\R} (f -U_\kappa f) \cdot \overline {(U_\kappa g) } dx =
\int_{\R} (U_\kappa f ) \cdot \overline {(U_\kappa g) } dx.
\end{multline*}
Сопоставляя эти равенства, заключаем, что (1.2.29) выполняется для
$ f \in L_p(\R) \cap L_2(\R),
g \in L_{p^\prime}(\R) \cap L_2(\R). $

Для вывода (1.2.29) для $ f \in L_p(\R), g \in L_{p^\prime}(\R) $
достаточно выбрать последовательность $ \{ f_n \in L_p(\R) \cap L_2(\R),
n \in \N \} $, сходящуюся к $ f $ в $ L_p(\R), $
и последовательность $ \{g_n \in L_{p^\prime}(\R) \cap L_2(\R), n \in \N \}, $
сходящуюся к $ g $ в $ L_{p^\prime}(\R), $ применить к ним на основании
доказанного (1.2.29), получив равенство
$$
\int_{\R} (E_\kappa f_n) \cdot \overline {g_n} dx = \int_{\R} f_n \cdot
\overline {(E_\kappa g_n) } dx,
$$

а затем, учитывая (1.2.20) и неравенство Гельдера, перейти к пределу в этом
равенстве при $ n \to \infty. $

Для вывода (1.2.30), используя (1.2.29), получаем
\begin{multline*}
\int_{\R} (\mathcal E_\kappa f) \cdot \overline g dx =
\int_{\R} ((E_\kappa f)
-(E_{\kappa -1} f))
\cdot \overline g dx =
\int_{\R} (E_\kappa f) \cdot \overline g dx -
\int_{\R} (E_{\kappa -1} f) \cdot \overline g dx = \\
\int_{\R} f \cdot \overline {(E_\kappa g) } dx -
\int_{\R} f \cdot \overline {(E_{\kappa -1} g) } dx =
\int_{\R} f \cdot (\overline {(E_\kappa g) } -
\overline {(E_{\kappa -1} g) }) dx =
\int_{\R} f \cdot \overline {(\mathcal E_\kappa g) } dx. \square
\end{multline*}

Предложение 1.2.4

Пусть выполнены условия предложения 1.2.3, а также $ X_0 \subset X_1 $
и $ \cup_{\kappa \in \Z_+} \span \{ \phi(2^\kappa \cdot -\nu),
\nu \in \Z \} $ плотно в $ L_2(\R). $ Тогда имеют место соотношения:

1) \begin{equation*} \tag{1.2.33}
\span \{ \phi(2^\kappa \cdot -\nu), \nu \in \Z\} \subset
\span \{ \phi(2^{\kappa +1} \cdot -\nu), \nu \in \Z \}, \kappa \in \Z_+;
\end{equation*}

2) при $ 1 \le p < \infty $ для $ \kappa, \kappa^\prime \in \Z_+:
\kappa^\prime \le  \kappa, $ выполняются равенства
\begin{equation*} \tag{1.2.34}
E_{\kappa^\prime}^p E_\kappa^p = E_\kappa^p
E_{\kappa^\prime}^p =
E_{\kappa^\prime}^p;
\end{equation*}

3) при $ 1 \le p < \infty $ для $ \kappa, \kappa^\prime \in \Z_+ $
соблюдаются равенства
\begin{equation*} \tag{1.2.35}
\mathcal E_\kappa^p \mathcal E_{\kappa^\prime}^p =
\begin{cases} \mathcal E_\kappa^p,
  \text{при $ \kappa = \kappa^\prime $};
       0,  \text{при $ \kappa \ne \kappa^\prime $};
\end{cases}
\end{equation*}
где $ \mathcal E_\kappa^p = \mathcal E_\kappa \mid_{L_p(\R)}, $
и при $ 1 < p < \infty, \kappa, \kappa^\prime \in \Z_+: \kappa \ne \kappa^\prime, $
для $ f \in L_p(\R), g \in L_{p^\prime}(\R) $ имеет место равенство
\begin{equation*} \tag{1.2.36}
\int_{\R} (\mathcal E_\kappa f) \cdot \overline (\mathcal
E_{\kappa^\prime} g) dx = 0;
\end{equation*}

4) при $ 1 \le p < \infty $ для $ \kappa \in \N $ соблюдаются равенства
$$
\Im \mathcal E_\kappa^p = \Im E_\kappa^p \cap \Ker E_{\kappa -1}^p;
\Ker \mathcal E_\kappa^p = \Im E_{\kappa -1}^p +\Ker E_\kappa^p;
$$

5) при $ 1 \le p < \infty $ в $ L_p(\R) $ справедливо равенство
\begin{equation*} \tag{1.2.37}
f = \sum_{ \kappa \in \Z_+} \mathcal E_\kappa f, f \in l_p(\R),
\end{equation*}
а для любой функции $ f \in
L_2(\R) $ имеет место равенство
\begin{equation*} \tag{1.2.38}
\| f \|_{L_2(\R)} = (\sum_{ \kappa \in \Z_+} \| \mathcal
E_\kappa f \|_{L_2(\R)}^2)^{1/2}.
\end{equation*}

Доказательство.

Прежде всего в условиях предложения с учетом (1.2.14) имеем
\begin{multline*} \tag{1.2.39}
X_\kappa = h_{2^\kappa} (X_0) \subset h_{2^\kappa} (X_1) =
h_{2^\kappa} (h_2 (X_0)) = (h_{2^\kappa} h_2) (X_0) = \\
h_{2^{\kappa +1}} (X_0) = X_{\kappa +1}.
\end{multline*}

Далее, при $ \mu \in \Z $ в силу  (1.2.39), (1.2.16), (1.2.26), (1.2.19),
а также компактности $ \supp \phi $ и $ \supp \phi^*, $ получаем
\begin{multline*}
\phi(2^\kappa \cdot -\mu) = U_{\kappa +1} \phi(2^\kappa \cdot -\mu) =
E_{\kappa +1} \phi(2^\kappa \cdot -\mu) = \\
\sum_{\substack{\nu \in \Z: \supp \phi(2^\kappa \cdot -\mu) \cap \\ (2^{-(\kappa +1)} \nu +
2^{-(\kappa +1)} \supp \phi^*) \ne \emptyset }} 2^{\kappa +1}
(\int_\R \phi(2^\kappa y -\mu) \overline {\phi^*(2^{\kappa +1} y -\nu)} dy)
\phi(2^{\kappa +1} \cdot -\nu) \in \\
\span \{ \phi(2^{\kappa +1} \cdot -\nu), \nu \in \Z\}, \kappa \in \Z_+,
\end{multline*}
что влечет (1.2.33).

Убедимся в справедливости (1.2.34).
Поскольку при $ \kappa^\prime, \kappa \in \Z_+: \kappa^\prime \le \kappa, $
вследствие (1.2.33), (1.2.27) имеет место включение
$ \Im E_{\kappa^\prime}^p \subset \Im E_\kappa^p, $
то для проектора $ E_\kappa^p $ вытекает, что
$$
E_\kappa^p (E_{\kappa^\prime}^p f) =
E_{\kappa^\prime}^p f,  f \in L_p(\R).
$$

Далее, в условиях п. 2 предложения для $ f \in L_p(\R) $ находим, что
$$
E_{\kappa^\prime}^p f = E_{\kappa^\prime}^p (E_\kappa^p
f) +E_{\kappa^\prime}^p (f -E_\kappa^p f) =
E_{\kappa^\prime}^p (E_\kappa^p f),
$$
ибо вследствие п. 2 предложения 1.2.3, (1.2.28), (1.2.33) справедливо
соотношение
\begin{multline*}
\Im (E -E_\kappa^p ) = \Ker E_\kappa^p = \\
\{ f \in L_p(\R): \int_\R f(x)
\overline {g(x)} dx =0 \forall g \in
\span \{\phi(2^\kappa \cdot -\nu),  \nu \in \Z\} \}
\subset \\
\{ f \in L_p(\R): \int_\R f(x)
\overline {g(x)} dx =0 \forall g \in
\span \{\phi(2^{\kappa^\prime} \cdot -\nu),  \nu \in \Z\} \} =
Ker E_{\kappa^\prime}^p,
\end{multline*}
и, значит, выполняется равенство
$$
E_{\kappa^\prime}^p (f -E_\kappa^p f) =0.
$$

Перейдем к проверке (1.2.35).
Пусть $ \kappa, \kappa^\prime \in \Z_+ $ и $ \kappa = \kappa^\prime. $
Тогда в силу (1.2.34) получаем
\begin{multline*}
(\mathcal E_\kappa^p)^2 = (E_\kappa^p)^2 - E_{\kappa
-1}^p E_\kappa^p -E_\kappa^p E_{\kappa -1}^p
+(E_{\kappa -1}^p)^2 =\\
= E_\kappa^p - E_{\kappa -1}^p -E_{\kappa -1}^p
+E_{\kappa -1}^p = E_\kappa^p - E_{\kappa -1}^p =
\mathcal E_\kappa^p.
\end{multline*}

Пусть $ \kappa \ne \kappa^\prime. $ Предположим, что $
\kappa^\prime < \kappa. $ Тогда, снова используя (1.2.34), выводим
\begin{multline*}
\mathcal E_{\kappa^\prime}^p \mathcal E_\kappa^p =
E_{\kappa^\prime}^p E_\kappa^p - E_{\kappa^\prime
-1}^p E_\kappa^p - E_{\kappa^\prime}^p E_{\kappa
-1}^p +
E_{\kappa^\prime -1}^p E_{\kappa -1}^p =\\
= E_{\kappa^\prime}^p - E_{\kappa^\prime -1}^p -
E_{\kappa^\prime}^p + E_{\kappa^\prime -1}^p =0.
\end{multline*}

Аналогично проверяется (1.2.35) при $
\kappa^\prime > \kappa. $

Для проверки (1.2.36) при $ f \in L_p(\R), g \in L_{p^\prime}(\R), $
применяя (1.2.30), (1.2.35), находим
\begin{multline*}
\int_{\R} (\mathcal E_\kappa f) \cdot \overline {(\mathcal
E_{\kappa^\prime} g)} dx =
\int_{\R} ( \mathcal E_{\kappa^\prime} \mathcal E_\kappa f) \cdot
\overline g dx = \\
\int_{\R} ( \mathcal E_{\kappa^\prime}^p \mathcal E_\kappa^p f) \cdot
\overline g dx =
\int_\R 0 \cdot \overline g dx =0, \kappa, \kappa^\prime \in \Z_+:
\kappa \ne \kappa^\prime.
\end{multline*}

Далее, сопоставляя (1.2.27), (1.2.28), (1.2.33) с (1.2.3), в соответствии
с (1.2.4) получаем равенства п. 4).

Для доказательства (1.2.37) с учетом предложения 1.1.1 покажем, что
при $ 1 \le p < \infty $ для $ f \in L_p(\R) $ имеет место соотношение
\begin{equation*} \tag{1.2.40}
\| f -E_\kappa f \|_{L_p(\R)} \to 0 \text{ при } \kappa \to \infty.
\end{equation*}

Сначала установим справедливость (1.2.40) при  $ p =2. $
Для $ f \in L_2(\R) $ и произвольного $ \epsilon >0, $ пользуясь плотностью
в $ L_2(\R) $ множества
$$
\cup_{\kappa \in \Z_+} \span \{ \phi(2^\kappa \cdot -\nu),
\nu \in \Z \},
$$
выберем $ \kappa_0 \in \Z_+ $ и $ g \in
\span \{ \phi(2^{\kappa_0} \cdot -\nu), \nu \in \Z \}, $ для которых
$$
\| f -g \|_{L_2(\R)} < \epsilon.
$$
Тогда при $ \kappa > \kappa_0 $ ввиду (1.2.26), (1.2.16), (1.2.39),
(1.2.20) находим, что
\begin{multline*}
\| f -E_\kappa f \|_{L_2(\R)} = \| f -g +g -U_\kappa f \|_{L_2(\R)} = \\
\| f -g +U_\kappa g -U_\kappa f \|_{L_2(\R)} =
\| f -g +E_\kappa (g -f) \|_{L_2(\R)} \le \\
\| f -g \|_{L_2(\R)} +
\| E_\kappa (g -f) \|_{L_2(\R)} \le
c_6 \| f -g \|_{L_2(\R)} < c_6 \epsilon.
\end{multline*}

Теперь убедимся в справедливости (1.2.40) при $ 1 \le p < \infty. $
Покажем сначала, что (1.2.40) имеет место для любой функции
$ f \in L_\infty(\R), $ имеющей компактный носитель. В этой ситуации,
выбирая в соответствии с (1.2.21) компакт $ K \subset \R $ такой, что
$ \supp f \cup \supp E_\kappa f \subset K, \kappa \in \Z_+, $
и используя (1.2.20) при $ p = \infty, $ а также (1.2.40)
при $ p =2, $ для произвольного $ \epsilon >0 $ получаем
\begin{multline*}
\| f - E_\kappa f \|_{L_p(\R)}^p =
\int_\R | f -(E_\kappa f) |^p dx =
\int_K | f -(E_\kappa f) |^p dx = \\
\int_{x \in K: | f(x) -(E_\kappa f)(x) | \le \epsilon}
| f -(E_\kappa f) |^p dx +\int_{x \in K: | f(x) -(E_\kappa f)(x) | >
\epsilon} | f -(E_\kappa f) |^p dx \le \\
\epsilon^p \mes K +\int_{x \in K: | f(x) -(E_\kappa f)(x) | >
\epsilon} | f -(E_\kappa f) |^p dx \le \\
\epsilon^p \mes K +\| f -(E_\kappa f) \|_{L_\infty(\R)}^p
\mes \{x \in K: | f(x) -(E_\kappa f)(x) | > \epsilon \} \le \\
\epsilon^p \mes K +(\| f \|_{L_\infty(\R)}
+\| E_\kappa f \|_{L_\infty(\R)})^p
\mes \{x \in K: | f(x) -(E_\kappa f)(x) | > \epsilon \} \le \\
\epsilon^p \mes K +(\| f \|_{L_\infty(\R)}
+\| E_\kappa f \|_{L_\infty(\R)})^p
\epsilon^{-2} \int_{ \{x \in K: | f(x) -(E_\kappa f)(x) | > \epsilon \}}
| f(x) -(E_\kappa f)(x) |^2 dx \le \\
\epsilon^p \mes K +(c_7 \| f \|_{L_\infty(\R)})^p
\epsilon^{-2} \int_\R | f(x) -(E_\kappa f)(x) |^2 dx = \\
\epsilon^p \mes K +(c_7 \| f \|_{L_\infty(\R)})^p
\epsilon^{-2} \| f -(E_\kappa f) \|_{L_2(\R)}^2 <
\epsilon^p (\mes K +(c_7 \| f \|_{L_\infty(\R)})^p )
\end{multline*}
для $ \kappa \in \Z_+ $ таких, что
$$
\| f -(E_\kappa f) \|_{L_2(\R)}^2 <
\epsilon^{p +2},
$$
т.е. (1.2.40) выполняется для любой функции $ f \in L_\infty(\R), $
имеющей компактный носитель. А поскльку множество таких функций
плотно в $ L_p(\R) $ при $ 1 \le p < \infty, $  то отсюда в силу (1.2.20)
заключаем, что (1.2.40) верно для любой $ f \in L_p(\R). $

Равенство (1.2.37) ввиду (1.1.1) является следствием
соотношения (1.2.40).

Наконец, обращаясь к выводу (1.2.38), заметим, что для любого конечного
множества $ S \subset \Z_+ $ и любого набора функций $ \{ g_\kappa \in
\Im \mathcal E_\kappa^2, \kappa \in S \} $ верно равенство
\begin{equation*} \tag{1.2.41}
\| \sum_{ \kappa \in S} g_\kappa \|_{L_2(\R)} = (\sum_{ \kappa
\in S} \| g_\kappa \|_{L_2(\R)}^2)^{1/2}.
\end{equation*}

В самом деле, ввиду (1.2.36) находим, что
\begin{multline*}
\| \sum_{ \kappa \in S} g_\kappa \|_{L_2(\R)}^2 = \int_{\R} |
\sum_{ \kappa \in S} g_\kappa|^2 dx =
\int_{\R} ( \sum_{ \kappa \in S} g_\kappa) ( \sum_{ \kappa^\prime \in S}
\overline { g_{\kappa^\prime}}) dx =\\
\int_{\R} \sum_{ \kappa \in S} \sum_{ \kappa^\prime \in S}
g_\kappa \overline {g_{\kappa^\prime}} dx =
\sum_{ \kappa \in S}
 \sum_{ \kappa^\prime \in S}
\int_{\R} g_\kappa \overline {g_{\kappa^\prime}} dx = \sum_{ \kappa \in S}
\int_{\R} | g_\kappa |^2 dx =
\sum_{ \kappa \in S} \| g_\kappa
\|_{L_2(\R)}^2,
\end{multline*}

откуда выводим (1.2.41).

Для получения равенства (1.2.38), в условиях теоремы благодаря (1.2.37)
и (1.2.41), выводим
\begin{multline*}
\| f \|_{L_2(\R)} = \| \lim_{ k \to \infty} \sum_{\kappa =0}^k
\mathcal E_\kappa f \|_{L_2(\R)} =
= \lim_{ k \to \infty} \|\sum_{\kappa =0}^k \mathcal E_\kappa f \|_{L_2(\R)}
= \\
\lim_{ k \to \infty} ( \sum_{\kappa =0}^k
\| \mathcal E_\kappa f \|_{L_2(\R)}^2)^{1/2} =
= ( \lim_{ k \to \infty} \sum_{\kappa =0}^k
\| \mathcal E_\kappa f \|_{L_2(\R)}^2)^{1/2} =
(\sum_{\kappa =0}^\infty
\| \mathcal E_\kappa f \|_{L_2(\R)}^2)^{1/2}. \square
\end{multline*}

Отметим, что функция $ \phi $ из предложения 1.2.4 является масштабирующей
для соответствующего кратно-масштабного анализа (см., например, [6], [7]).
\bigskip

1.3. В этом пункте приведем некоторые вспомогательные утверждения,
которые используются в следующем пункте и далее.

При доказательстве леммы 1.3.2 используется

Лемма 1.3.1

Пусть $ d \in \N, 1 \le p < \infty. $ Тогда если для функции $ f \in
L_p(\R^d) $ последовательность $ \{ f_n \in L_p(\R^d), n \in \N \} $ сходится
в $ L_p(\R^d) $ к $ f, $ то для каждого $ j=1,\ldots,d $ существует
последовательность $ \{ n_k \in \N: n_k < n_{k+1},k \in \N\} $ такая, что
почти для всех $ (x_1,\ldots,x_{j-1},x_{j+1},\ldots,x_d) \in \R^{d-1} $
подпоследовательность $ \{ f_{n_k}(x_1,\ldots,x_{j-1},\cdot,x_{j+1},\ldots,x_d),
k \in \N\} $ сходится в $ L_p(\R) $ к функции $ f(x_1,\ldots,x_{j-1},\cdot,
x_{j+1},\ldots,x_d). $

Доказательство.

Фиксировав $ j=1,\ldots,d, $ в условиях леммы, вследствие теоремы
Фубини, получаем, что почти для всех $
(x_1,\ldots,x_{j-1},x_{j+1},\ldots,x_d) \in \R^{d-1} $ каждая функция
$ f(x_1,\ldots,x_{j-1},\cdot,x_{j+1},\ldots,x_d),
f_n(x_1,\ldots,x_{j-1},\cdot,x_{j+1},\ldots,x_d), n \in \N, $
принадлежит $ L_p(\R) $ и
\begin{multline*}
\|f -f_n\|_{L_p(\R^d)}^p =  \int_{ \R^d} |f(x) -f_n(x)|^p dx =
\int_{ \R^{d-1}} \int_\R |f(x_1,\ldots,x_{j-1},x_j,x_{j+1},\ldots,x_d)-\\
- f_n(x_1,\ldots,x_{j-1},x_j,x_{j+1},\ldots,x_d)|^p dx_j dx_1\ldots
dx_{j-1} dx_{j+1} \ldots dx_d \to 0 \text{ при } n \to \infty,
\end{multline*}
т.е. последовательность функций $ \{ \int_\R |f(x) -f_n(x)|^p dx_j,
n \in \N \} $ сходится к $ 0 $ в $ L_1(\R^{d-1}) $. Поэтому
последовательность $ \{ \int_\R |f(x) -f_n(x)|^p dx_j, n \in \N \}
$ сходится к $ 0 $ по мере, и, значит, существует
последовательность $ \{ n_k \in \N: n_k < n_{k+1}, k \in \N\} $
такая, что почти для всех $ (x_1,\ldots,x_{j-1}, x_{j+1},\ldots,x_d)
\in \R^{d-1} $ подпоследовательность  \linebreak $ \{ \int_\R
|f(x_1,\ldots,x_{j-1},x_j,x_{j+1},\ldots,x_d)
-f_{n_k}(x_1,\ldots,x_{j-1},x_j,x_{j+1},\ldots,x_d)|^p dx_j =
\|f(x_1,\ldots,x_{j-1},\cdot,x_{j+1},\ldots,x_d)
-f_{n_k}(x_1,\ldots,x_{j-1},\cdot,x_{j+1},\ldots,x_d)\|_{L_p(\R)}^p,
k \in \N\} $ сходится к $ 0, $ что и завершает доказательство
леммы. $ \square $

Доказательство леммы 1.3.2 вполне аналогично доказательству
соответствующего утверждения из [9], [10].

Лемма 1.3.2

    Пусть $ d \in \N, 1 \le p < \infty. $ Тогда

1) при  $ j=1,\ldots,d$ для любого непрерывного линейного оператора $ T: L_p(\R)
\mapsto L_p(\R) $ существует единственный непрерывный линейный оператор
$ \mathcal T^j: L_p(\R^d) \mapsto L_p(\R^d), $ для которого для любой функции
$ f \in L_p(\R^d) $ почти для всех $ (x_1,\ldots,x_{j-1},x_{j+1},
\ldots,x_d) \in \R^{d-1} $ в $ L_p(\R) $ выполняется равенство
\begin{equation*} \tag{1.3.1}
(\mathcal T^j f)(x_1,\ldots, x_{j-1},\cdot,x_{j+1},
\ldots,x_d) = (T(f(x_1,\ldots,x_{j-1},\cdot,x_{j+1},
\ldots,x_d)))(\cdot),
\end{equation*}

2) при этом, для каждого $ j=1,\ldots,d $ отображение $ V_j^{L_p}, $ которое
каждому оператору $ T \in \mathcal B(L_p(\R), L_p(\R)) $ ставит в соответствие
оператор $ V_j^{L_p}(T) = \mathcal T^j \in \mathcal B(L_p(\R^d), L_p(\R^d)), $
удовлетворяющий (1.3.1), является непрерывным гомоморфизмом банаховой алгебры
$ \mathcal B(L_p(\R), L_p(\R)) $ в банахову алгебру $ \mathcal B(L_p(\R^d),
L_p(\R^d)), $
и
\begin{equation*} \tag{1.3.2}
\| V_j^{L_p}(T) \| _{\mathcal B(L_p(\R^d), L_p(\R^d))} \le
\| T \|_{ \mathcal B(L_p(\R), L_p(\R))},
\end{equation*}

3) причём, для любых операторов $ S,T \in \mathcal B(L_p(\R), L_p(\R)) $
при любых $ i,j =1,\ldots,d: i\ne j, $ соблюдается равенство
\begin{equation*} \tag{1.3.3}
(V_i^{L_p}(S) V_j^{L_p}(T))f = (V_j^{L_p}(T) V_i^{L_p}(S))f, f \in L_p(\R^d).
\end{equation*}

Доказательство.

Фиксируем $ j=1,\ldots,d $ и рассмотрим линейное подпространство
$ L_p^j(\R^d), $ состоящее из тех функций $ f \in L_p(\R^d), $ для
которых при некотором $ n \in \N $ существуют наборы функций $ \{ \phi_i
\in L_p(\R^{d-1}), i=1,\ldots,n \} $ и $ \{\psi_i \in L_p(\R), i=1,\ldots,n \} $
такие, что в $ L_p(\R^d) $ выполняется равенство
\begin{equation*} \tag{1.3.4}
f(x_1,\ldots,x_{j-1},x_j,x_{j+1},\ldots,x_d) = \sum_{i=1}^n \phi_i(x_1,\ldots,
x_{j-1},x_{j+1},\ldots,x_d) \psi_i(x_j),
\end{equation*}
и, следовательно, ввиду теоремы Фубини, почти для всех $ (x_1,\ldots,x_{j-1},
x_{j+1},\ldots,x_d) \in \R^{d-1} $ в $ L_p(\R) $ соблюдается равенство
\begin{equation*} \tag{1.3.5}
f(x_1,\ldots,x_{j-1},\cdot,x_{j+1},\ldots,x_d) = \sum_{i=1}^n \phi_i(x_1,\ldots,
x_{j-1},x_{j+1},\ldots,x_d) \psi_i(\cdot).
\end{equation*}

Пусть $ T: L_p(\R) \mapsto L_p(\R) $ -- непрерывный линейный оператор. Тогда
каждой функции $ f \in L_p^j(\R^d), $ для которой в $ L_p(\R^d) $ имеет место
(1.3.4), сопоставим функцию $ \mathcal T^j f \in L_p(\R^d), $ определяемую в
$ L_p(\R^d) $ равенством
$$
(\mathcal T^j f)(x_1,\ldots,x_{j-1},x_j,x_{j+1},\ldots,x_d) = \sum_{i=1}^n
\phi_i(x_1,\ldots,x_{j-1},x_{j+1},\ldots,x_d) (T \psi_i)(x_j).
$$

Из этого равенства, принимая во внимание теорему Фубини, линейность оператора
$ T $ и (1.3.5), получаем, что почти для всех  \linebreak $ (x_1,\ldots,x_{j-1},x_{j+1},
\ldots,x_d) \in \R^{d-1} $ в $ L_p(\R) $ справедливо равенство
\begin{multline*}
(\mathcal T^j f)(x_1,\ldots,x_{j-1},\cdot,x_{j+1},\ldots,x_d) = \\
\sum_{i=1}^n \phi_i(x_1,\ldots,x_{j-1},x_{j+1},\ldots,x_d) (T \psi_i)(\cdot) =\\
(T (\sum_{i=1}^n \phi_i(x_1,\ldots,x_{j-1},x_{j+1},\ldots,x_d) \psi_i(\cdot)))
(\cdot) =\\
 (T f(x_1,\ldots,x_{j-1},\cdot,x_{j+1},\ldots,x_d))(\cdot),
\end{multline*}
что совпадает с (1.3.1).

Из последнего равенства, в частности, вытекает, что приведенное выше
определение функции $ \mathcal T^j f $ для функции $ f \in L_p^j(\R^d) $ не
зависит от выбора представления (1.3.4).

Ясно, что построенное отображение $ \mathcal T^j: L_p^j(\R^d) \mapsto L_p(\R^d) $ --
линейно, и для $ f \in L_p^j(\R^d), $ благодаря теореме Фубини, соотношению
(1.3.1) и непрерывности оператора $ T, $ имеем
\begin{multline*}
\| \mathcal T^j f \|_{L_p(\R^d)}^p \\ = \int_{ \R^{d-1}}
\|(\mathcal T^j f)(x_1,\ldots,x_{j-1},\cdot,x_{j+1},\ldots,x_d)\|_{L_p(\R)}^p dx_1
\ldots dx_{j-1} dx_{j+1} \ldots dx_d \\
= \int_{ \R^{d-1}} \| T f(x_1,\ldots,x_{j-1},
\cdot,x_{j+1},\ldots,x_d)\|_{L_p(\R)}^p dx_1 \ldots dx_{j-1} dx_{j+1} \ldots dx_d
\\ \le \int_{ \R^{d-1}} (\|T\| \cdot \|f(x_1,\ldots,x_{j-1},\cdot,x_{j+1},\ldots,
x_d)\|_{L_p(\R)})^p dx_1 \ldots dx_{j-1} dx_{j+1} \ldots dx_d \\
=\|T\|^p \|f\|_{L_p(\R^d)}^p
\end{multline*}
или
\begin{equation*} \tag{1.3.6}
\| \mathcal T^j f\|_{L_p(\R^d)} \le \|T\| \cdot \|f\|_{L_p(\R^d)}.
\end{equation*}

Далее, заметим, что поскольку $ L_p^j(\R^d) $ содержит всюду плотное в
$ L_p(\R^d) $ подпространство, являющееся линейной оболочкой
функций вида $ \prod_{k=1}d (x_k^{\lambda_k} \chi_k(x_k)), $ где
$ \chi_k $ -- характеристическая функция некоторого интервала
$ (-\delta_k, \delta_k) \subset \R, \delta_k >0, \lambda_k \in \Z_+,
k = 1, \ldots, d, $
то подпространство $ L_p^j(\R^d) $ --
всюду плотно в $ L_p(\R^d). $

Принимая во внимание это обстоятельство, для каждой функции $ f \in L_p(\R^d) $
выберем последовательность функций $ \{ f_n \in L_p^j(\R^d), n \in \N\} $,
сходящуюся в $ L_p(\R^d) $ к $ f, $ и, пользуясь тем, что в силу (1.3.6) при
$ n,m \in \N $ имеет место оценка
$$
\|\mathcal T^j f_n -\mathcal T^j f_m\|_{L_p(\R^d)}
\le \|T\| \cdot \|f_n -f_m\|_{L_p(\R^d)},
$$
а, следовательно, последовательность $ \{ \mathcal T^j f_n, n \in \N \} $ --
фундаментальна в банаховом пространстве $ L_p(\R^d), $ положим
\begin{equation*} \tag{1.3.7}
\mathcal T^j f = \lim_{n\to \infty} \mathcal T^j f_n \text{в} L_p(\R^d).
\end{equation*}

Из (1.3.6) вытекает, что определяемое посредством (1.3.7) значение
$ \mathcal T^j f, $ не зависит от выбора последовательности $ \{ f_n \in
L_p^j(\R^d), n \in \N\}, $ сходящейся к $ f $ в $ L_p(\R^d). $

Заданное таким образом отображение $ \mathcal T^j: L_p(\R^d) \mapsto L_p(\R^d), $
в силу линейности $ \mathcal T^j $ на $ L_p^j(\R^d) $ и соотношения (1.3.7), --
также линейно.

Из соблюдения неравенства (1.3.6) для $ f \in L_p^j(\R^d) $ и равенства
(1.3.7) следует, что (1.3.6) имеет место для $ f \in L_p(\R^d). $ Тем
самым, $ \mathcal T^j: L_p(\R^d) \mapsto L_p(\R^d) $ -- непрерывный оператор.

Проверим справедливость (1.3.1) для $ f \in L_p(\R^d). $ Для этого, выбирая
для $ f \in L_p(\R^d) $ последовательность $ \{ f_n \in L_p^j(\R^d),
n \in \N\}, $ сходящуюся к $ f $ в $ L_p(\R^d), $ ввиду (1.3.1) для каждого
$ n \in \N $ почти для всех $ (x_1,\ldots,x_{j-1},x_{j+1},\ldots,x_d) \in
\R^{d-1} $ в $ L_p(\R) $ имеем равенство
\begin{equation*} \tag{1.3.8}
(\mathcal T^j f_n)(x_1,\ldots,x_{j-1},\cdot,x_{j+1},\ldots,x_d) = (T f_n(x_1,\ldots,
x_{j-1},\cdot,x_{j+1},\ldots,x_d))(\cdot),
\end{equation*}
откуда получаем, что почти для всех $ (x_1,\ldots,x_{j-1},x_{j+1},\ldots,x_d)
\in \R^{d-1} $ для всех $ n \in \N $ в $ L_p(\R) $ имеет место (1.3.8).

Учитывая (1.3.7), переходя на основании леммы 1.3.1 к соответствующей
последовательности $ \{n_k \in \N, k \in \N\}, $ заключаем, что почти для
всех $ (x_1,\ldots,x_{j-1},x_{j+1},\ldots,x_d) \in \R^{d-1} $ в $ L_p(\R) $
выполняется равенство
\begin{equation*} \tag{1.3.9}
(\mathcal T^j f)(x_1,\ldots,x_{j-1},\cdot,x_{j+1},\ldots,x_d) = \lim_{ k \to \infty}
 (\mathcal T^j f_{n_k})(x_1,\ldots,x_{j-1},\cdot,x_{j+1},\ldots,x_d),
\end{equation*}
а также почти для всех $ (x_1,\ldots,x_{j-1},x_{j+1},\ldots,x_d) \in \R^{d-1} $
подпоследовательность $ \{ f_{n_k}(x_1,\ldots,x_{j-1},\cdot, x_{j+1},\ldots,x_d),
k \in \N\} $ сходится в $ L_p(\R) $ к \linebreak $ f(x_1,\ldots,x_{j-1},\cdot,x_{j+1},
\ldots,x_d), $ и, следовательно, в силу непрерывности оператора $ T: L_p(\R)
\mapsto L_p(\R) $ почти для всех $ (x_1,\ldots,x_{j-1},x_{j+1},\ldots,x_d)  \in
\R^{d-1} $ в $ L_p(\R) $ справедливо равенство
\begin{equation*} \tag{1.3.10}
(T f(x_1,\ldots,x_{j-1},\cdot,x_{j+1},\ldots,x_d))(\cdot) = \lim_{ k \to \infty}
(T f_{n_k}(x_1,\ldots,x_{j-1},\cdot,x_{j+1},\ldots,x_d))(\cdot).
\end{equation*}

Соединяя (1.3.8), (1.3.9) и (1.3.10), приходим к выводу, что для $ f \in
L_p(\R^d) $ почти для всех $ (x_1,\ldots,x_{j-1},x_{j+1},\ldots,x_d)  \in
\R^{d-1} $ в $ L_p(\R) $ имеет место (1.3.1).

При $ j=1,\ldots,d $ для операторов $ S,T \in \mathcal B(L_p(\R),
L_p(\R)), $ функции $ f \in L_p(\R^d) $ почти для всех
$ (x_1,\ldots,x_{j-1},x_{j+1},\ldots,x_d) \in \R^{d-1} $  в $ L_p(\R)$, согласно
(1.3.1)  имеем
\begin{multline*}
(V_j(S+T)f)(x_1,\ldots,x_{j-1},\cdot, x_{j+1},\ldots, x_d) =\\
((S+T)f(x_1,\ldots,x_{j-1},\cdot, x_{j+1},\ldots, x_d))(\cdot) = \\
(S f(x_1,\ldots,x_{j-1},\cdot, x_{j+1},\ldots, x_d))(\cdot) +
(T f(x_1,\ldots,x_{j-1},\cdot, x_{j+1},\ldots, x_d))(\cdot) = \\
(V_j(S)f)(x_1,\ldots,x_{j-1},\cdot, x_{j+1},\ldots, x_d) +
(V_j(T)f)(x_1,\ldots,x_{j-1},\cdot, x_{j+1},\ldots, x_d),
\end{multline*}
а поэтому
$$
V_j(S+T)f = V_j(S)f +V_j(T)f
$$
и
 \begin{multline*}
(V_j(ST)f)(x_1,\ldots,x_{j-1},\cdot, x_{j+1},\ldots, x_d) =\\
((ST)f(x_1,\ldots,x_{j-1},\cdot, x_{j+1},\ldots, x_d))(\cdot) = \\
(S(Tf(x_1,\ldots,x_{j-1},\cdot, x_{j+1},\ldots, x_d)))(\cdot) = \\
(S(V_j(T)f)(x_1,\ldots,x_{j-1},\cdot, x_{j+1},\ldots, x_d))(\cdot) = \\
(V_j(S)(V_j(T)f))(x_1,\ldots,x_{j-1},\cdot, x_{j+1},\ldots, x_d) = \\
((V_j(S)V_j(T))f)(x_1,\ldots,x_{j-1},\cdot, x_{j+1},\ldots, x_d).
\end{multline*}

Неравенство (1.3.2) следует из (1.3.6).

Переходя к доказательству (1.3.3), заметим, что при $ j =1,\ldots,d $
для любой функции $ f $ вида
\begin{equation*} \tag{1.3.11}
f(x) = \prod_{k=1}^d f_k(x_k),
f_k \in L_p(\R), k =1, \ldots, d,
\end{equation*}
почти для всех  $ (x_1,\ldots,x_{j-1}, x_{j+1},\ldots, x_d) \in \R^{d-1} $
в  $ L_p(\R) $ выполняется равенство
\begin{multline*}
(V_j(T) f)(x_1,\ldots,x_{j-1},\cdot, x_{j+1},\ldots, x_d) =\\
(T(f_j(\cdot) \prod_{\{k=1,\ldots,d: k \ne j\}} f_k(x_k)))(\cdot) =
(T(f_j))(\cdot) \prod_{\{k=1,\ldots,d: k \ne j\}} f_k(x_k),
\end{multline*}
и, следовательно,
$$
(V_j(T) f)(x) = (T(f_j))(x_j) \prod_{\{k=1,\ldots,d:
k \ne j\}} f_k(x_k).
$$

В силу тех же соображений, что и выше, при $ i, j =1,\ldots,d: i \ne j, $
почти для всех $ (x_1,\ldots,x_{i-1},
x_{i+1},\ldots, x_d) \in \R^{d-1} $ в  $ L_p(\R) $ справедливо равенство
\begin{multline*}
((V_i(S)V_j(T)) f)(x_1,\ldots,x_{i-1},\cdot, x_{i+1},\ldots, x_d) =\\
(S(V_j(T) f) (x_1,\ldots,x_{i-1},\cdot, x_{i+1},\ldots, x_d))(\cdot) = \\
(S(f_i(\cdot) (T f_j)(x_j) \prod_{\{k =1,\ldots,d: k \ne i,j\}}
f_k(x_k)))(\cdot) =
(S(f_i))(\cdot) (T(f_j))(x_j)
\prod_{\{k =1,\ldots,d: k \ne i,j\}}
f_k(x_k)
\end{multline*}
и, значит,
$$
((V_i(S)V_j(T)) f)(x) = (S(f_i))(x_i) (T(f_j))
(x_j) \prod_{\{k =1,\ldots,d: k \ne i,j\}}
f_k(x_k).
$$
Точно так же
$$
((V_j(T)V_i(S)) f)(x) = (S(f_i))(x_i) (T(f_j))
(x_j) \prod_{\{k =1,\ldots,d: k \ne i,j\}}
f_k(x_k).
$$

Таким образом, равенство (1.3.3) справедливо для люлой функции
$ f $ вида (1.3.11).
Ввиду линейности операторов $ V_i(S), V_j(T) $
отсюда следует, что равенство (1.3.3) соблюдается для любой линейной комбинации таких функцй.
На основании сказанного, учитывая плотность в $ L_p(\R^d) $ линейной
оболочки всех функций вида (1.3.11) и непрерывность в $ L_p(\R^d) $
операторов $V_i(S) $ и $ V_j(T), $ заключаем, что (1.3.3) выполняется
для любой функции $ f \in L_p(\R^d). $
Плотность в $ L_p(\R^d) $ линейной оболочки всех функций вида (1.3.11)
следует из того, что она содержит все функции вида
$$
\sum_{ \lambda \in \Z_+^d(l)} a_\lambda x^\lambda \chi(x),
$$
где $ a_\lambda \in \C, \lambda \in \Z_+^d(l), l \in \Z_+^d,
x^\lambda = x_1^{\lambda_1} \ldots x_d^{\lambda_d}, \chi $ --
характеристическая функция некоторого куба$ \delta B^d, \delta \in \R_+. $
$ \square $

Замечание.

Если при $ d \in \N, 1 \le p, q < \infty, $ оператор $ T \in
\mathcal B(L_p(\R), L_p(\R)) \cap \mathcal B(L_q(\R), L_q(\R)), $ то
при $ j=1,\ldots,d $ для $ f \in L_p(\R^d) \cap L_q(\R^d)$ справедливо равенство
\newline
$ (V_j^{L_p} T)f = (V_j^{L_q} T)f. $ Поэтому символы $ L_p, L_q $
в качестве индексов у $ V_j $ можно опускать.
\bigskip

1.4. В этом пункте дается описание свойств операторов проектирования
на подпространства всплесков, соответствующие кратно-масштабному
анализу, порожденному тензорным произведением ограниченных финитных
функций, и других проекторов, которые используются в п. 2.2. при
доказательстве основных результатов работы.

Предложение 1.4.1

Пусть $ d \in \N $ и при $ j =1, \ldots, d $ функции $ \phi_j, \phi_j^* $
удовлетворяют условиям предложения 1.2.2, а
$ \phi(x) = \prod_{j =1}^d \phi_j(x_j),
\phi^*(x) = \prod_{j =1}^d \phi_j^*(x_j), x \in \R^d. $
Тогда при $ 1 \le p < \infty, \kappa \in \Z_+^d $ справедливы следующие
утверждения:

1) операторы
$$
E_\kappa^p = E_\kappa^{\phi, \phi^*} \mid_{L_p(\R^d)},
E_{\kappa_j}^{j, p} =
E_{\kappa_j}^{\phi_j, \phi_j^*}
\mid_{L_p(\R)}, j =1,\ldots,d,
$$
связаны равенством
\begin{equation*} \tag{1.4.1}
E_\kappa^p = \prod_{j=1}^d V_j (E_{\kappa_j}^{j, p});
\end{equation*}

2) для линейных операторов
\begin{eqnarray*}
\mathcal E_\kappa &= \mathcal E_\kappa^{\phi, \phi^*}: L_1^{\loc(\R^d)}
\mapsto L_1^{\loc(\R^d)}, \\
\mathcal E_\kappa^{p}&: L_p(\R^d) \mapsto L_p(\R^d), \\
\mathcal E_{\kappa_j}^{j, p}&: L_p(\R) \mapsto L_p(\R), j =1,\ldots,d,
\end{eqnarray*}
определяемых соотношениями
\begin{eqnarray*}
\mathcal E_\kappa &= \sum_{\epsilon \in \Upsilon^d:
\s(\epsilon) \subset \s(\kappa)} (-\e)^\epsilon E_{\kappa -\epsilon}, \\
\mathcal E_\kappa^{p} &= \sum_{\epsilon \in \Upsilon^d:
\s(\epsilon) \subset \s(\kappa)} (-\e)^\epsilon E_{\kappa -\epsilon}^{p} =
\mathcal E_\kappa^{\phi, \phi^*} \mid_{L_p(\R^d)}, \\
\mathcal E_{\kappa_j}^{j, p} &=
\mathcal E_{\kappa_j}^{\phi_j, \phi_j^*}
\mid_{L_p(\R)}, j =1,\ldots,d,
\end{eqnarray*}
имеет место равенство
\begin{equation*} \tag{1.4.2}
\mathcal E_\kappa^{p} = \prod_{j=1}^d V_j(\mathcal E_{\kappa_j}^{j, p}).
\end{equation*}

Доказательство.

Начнем с доказательства первого утверждения. Для этого в условиях
предложения сначала покажем, что при $ 1 \le p < \infty, \kappa \in \Z_+^d $
для любого непутого множества $ J \subset \{1,\ldots,d\}, $
для любой функции $ f \in L_p(\R^d) $ такой, что
существует $ \delta \in \R_+^d, $ для которого $ f =0 $ почти всюду на
множестве  $ \{ x \in \R^d: x^J \notin (\delta B^d)^J \}, $ имеет место
равенство
\begin{multline*} \tag{1.4.3}
((\prod_{j \in J} V_j (E_{\kappa_j}^{j, p})) f)(x) =
\sum_{\nu^J \in (\Z^d)^J: (\supp \phi_j^*(2^{\kappa_j} \cdot -\nu_j))
\cap (\delta_j B^1) \ne \emptyset, j \in J} 2^{(\kappa^J, \e^J)} \cdot \\
(\int_{(\delta B^d)^J \cap \prod_{j \in J}
(\supp \phi_j^*(2^{\kappa_j} \cdot -\nu_j))} f(x \chi_{J^\prime} +y \chi_J) \cdot \\
(\prod_{j \in J} \overline{\phi_j^*(2^{\kappa_j} y_j -\nu_j)}) dy^J)
(\prod_{j \in J} \phi_j(2^{\kappa_j} x_j -\nu_j)),
\end{multline*}
где $ J^\prime = \{1,\ldots,d\} \setminus J. $

Доказательство (1.4.3) проведем по индукции относительно $ \card J. $
При $ \card J = 1, $ т.е. для $ J = \{j\},
j \in \{1,\ldots,d\}, $ для $ f \in L_p(\R^d) $ такой, что при некотором
$ \delta \in \R_+^d $ значение функции $ f(x)  =0 $ почти для всех
$ x \in \R^d: x_j \notin \delta_j B^1, $ ввиду (1.3.1), (1.2.19) почти для всех
$ (x_1, \ldots, x_{j -1}, x_{j +1}, \ldots, x_d) \in \R^{d -1} $ имеем
\begin{multline*}
((V_j (E_{\kappa_j}^{j, p})) f)(x_1, \ldots, x_{j -1}, \cdot,
x_{j +1}, \ldots, x_d) = \\
(E_{\kappa_j}^{j, p} f(x_1, \ldots, x_{j -1}, \cdot,
x_{j +1}, \ldots, x_d))(\cdot) = \\
\sum_{\nu_j \in \Z} 2^{\kappa_j} \biggl(\int_{\R} f(x_1, \ldots, x_{j -1}, y_j,
x_{j +1}, \ldots, x_d) \overline {\phi_j^*(2^{\kappa_j} y_j -\nu_j)} dy_j\biggr)
\phi_j(2^{\kappa_j} \cdot -\nu_j) = \\
\sum_{\nu_j \in \Z} 2^{\kappa_j} \biggl(\int_{(\delta_j B^1) \cap
\supp \phi_j^*(2^{\kappa_j} \cdot -\nu_j)} f(x_1, \ldots, x_{j -1}, y_j,
x_{j +1}, \ldots, x_d)\cdot\\  \overline {\phi_j^*(2^{\kappa_j} y_j -\nu_j)} dy_j\biggr)
\phi_j(2^{\kappa_j} \cdot -\nu_j) = \\
\sum_{\substack{\nu_j \in \Z: (\delta_j B^1) \cap \\
(\supp \phi_j^*(2^{\kappa_j} \cdot -\nu_j)) \ne \emptyset}} 2^{\kappa_j}
\biggl(\int_{(\delta_j B^1) \cap \supp \phi_j^*(2^{\kappa_j} \cdot -\nu_j)}
f(x_1, \ldots, x_{j -1}, y_j, x_{j +1}, \ldots, x_d) \cdot\\
\overline {\phi_j^*(2^{\kappa_j} y_j -\nu_j)} dy_j\biggr)
\phi_j(2^{\kappa_j} \cdot -\nu_j),
\end{multline*}
откуда
\begin{multline*} \tag{1.4.4}
((V_j (E_{\kappa_j}^{j, p})) f)(x) =\\
\sum_{\substack{\nu_j \in \Z: (\delta_j B^1) \cap \\
(\supp \phi_j^*(2^{\kappa_j} \cdot -\nu_j)) \ne \emptyset}} 2^{\kappa_j}
\biggl(\int_{(\delta_j B^1) \cap (\supp \phi_j^*(2^{\kappa_j} \cdot -\nu_j)) }
f(x_1, \ldots, x_{j -1}, y_j, x_{j +1}, \ldots, x_d) \cdot\\
\overline {\phi_j^*(2^{\kappa_j} y_j -\nu_j)} dy_j\biggr)
\phi_j(2^{\kappa_j} x_j -\nu_j),
\end{multline*}
что совпадает с (1.4.3) в  случае $ \card J =1. $

Предположим, что равенство (1.4.3) справедливо для любого множества
$ J \subset \{1, \ldots, d\}: \card J = m \le d -1. $ Проверим, что тогда оно
имеет место для любого множества $ \mathcal J \subset \{1, \ldots, d\}:
\card \mathcal J = m +1. $
Представляя такое множество в виде $ \mathcal J = \{i\} \cup J, $ где
$ \card J = m, i \notin J, $ для  функции $ f \in L_p(\R^d) $ такой, что существует
$ \delta \in \R_+^d, $ для которого $ f =0 $ почти всюду на множестве
$ \{ x \in \R^d: x^{\mathcal J} \notin (\delta B^d)^{\mathcal J} \}, $
на основании (1.3.3) и предположения индукции, учитывая, что
$ \{ x \in \R^d: x^J \notin (\delta B^d)^J \} \subset \{ x \in \R^d:
x^{\mathcal J} \notin (\delta B^d)^{\mathcal J} \}, $ получаем
\begin{multline*} \tag{1.4.5}
((\prod_{j \in \mathcal J} V_j (E_{\kappa_j}^{j, p})) f)(x) =
((V_i (E_{\kappa_i}^{i, p}))
((\prod_{j \in J} V_j (E_{\kappa_j}^{j, p})) f))(x) = \\
((V_i (E_{\kappa_i}^{i, p}))
(\sum_{\substack{\nu^J \in (\Z^d)^J: (\supp \phi_j^*(2^{\kappa_j} \cdot -\nu_j))
\cap \\ (\delta_j B^1) \ne \emptyset, j \in J}} 2^{(\kappa^J, \e^J)}
(\int_{(\delta B^d)^J \cap \prod_{j \in J}
(\supp \phi_j^*(2^{\kappa_j} \cdot -\nu_j))} f(u \chi_{J^\prime} +y \chi_J) \cdot \\
(\prod_{j \in J} \overline{\phi_j^*(2^{\kappa_j} y_j -\nu_j)} dy^J)
(\prod_{j \in J} \phi_j(2^{\kappa_j} u_j -\nu_j))))(x) = \\
\sum_{\nu^J \in (\Z^d)^J: (\supp \phi_j^*(2^{\kappa_j} \cdot -\nu_j))
\cap (\delta_j B^1) \ne \emptyset, j \in J} 2^{(\kappa^J, \e^J)} \cdot \\
((V_i (E_{\kappa_i}^{i, p}))
((\int_{(\delta B^d)^J \cap \prod_{j \in J}
(\supp \phi_j^*(2^{\kappa_j} \cdot -\nu_j))} f(u \chi_{J^\prime} +y \chi_J) \cdot \\
(\prod_{j \in J} \overline{\phi_j^*(2^{\kappa_j} y_j -\nu_j)} dy^J)
(\prod_{j \in J} \phi_j(2^{\kappa_j} u_j -\nu_j))))(x).
\end{multline*}
Теперь, учитывая, что
$$
\{u \in \R^d: u_i \notin (\delta_i B^1)\} \subset
\{u \in \R^d: u^{\mathcal J} \notin (\delta B^d)^{\mathcal j} \},
$$
видим, что $ f(u) =0 $ почти для всех $ u \in \R^d: u_i \notin (\delta_i B^1), $
а, следовательно, для $ \nu^J \in (\Z^d)^J: (2^{-\kappa_j} \nu_j +
2^{-\kappa_j} \supp \phi_j^*) \cap (\delta_j B^1) \ne \emptyset, j \in J, $
почти для всех $ u \in \R^d: u_i \notin (\delta_i B^1), $ выполняется равенство
\begin{multline*}
\biggl(\int_{(\delta B^d)^J \cap \prod_{j \in J}
(\supp \phi_j^*(2^{\kappa_j} \cdot -\nu_j))} f(u \chi_{J^\prime} +y \chi_J)\cdot \\
(\prod_{j \in J} \overline{\phi_j^*(2^{\kappa_j} y_j -\nu_j)}) dy^J\biggr)
(\prod_{j \in J} \phi_j(2^{\kappa_j} u_j -\nu_j)) =0.
\end{multline*}
Поэтому, обозначая $ \mathcal J^{\prime} = \{1,\ldots,d\} \setminus \mathcal J, $
согласно (1.4.4) с учетом теоремы Фубини находим, что для
$ \nu^J \in (\Z^d)^J: (\supp \phi_j^*(2^{\kappa_j} \cdot -\nu_j))
\cap (\delta_j B^1) \ne \emptyset, j \in J, $ справедливо равенство
\begin{multline*} \tag{1.4.6}
((V_i (E_{\kappa_i}^{i, p}))
((\int_{(\delta B^d)^J \cap \prod_{j \in J}
(\supp \phi_j^*(2^{\kappa_j} \cdot -\nu_j))} f(u \chi_{J^\prime} +y \chi_J) \cdot \\
(\prod_{j \in J} \overline{\phi_j^*(2^{\kappa_j} y_j -\nu_j)}) dy^J)
(\prod_{j \in J} \phi_j(2^{\kappa_j} u_j -\nu_j))))(x) = \\
\sum_{\nu_i \in \Z: (\delta_i B^1) \cap
(\supp \phi_i^*(2^{\kappa_i} \cdot -\nu_i)) \ne \emptyset} 2^{\kappa_i} \cdot \\
(\int_{(\delta_i B^1) \cap
(\supp \phi_i^*(2^{\kappa_i} \cdot -\nu_i)) }
((\int_{(\delta B^d)^J \cap \prod_{j \in J}
(\supp \phi_j^*(2^{\kappa_j} \cdot -\nu_j))}
f(x \chi_{\mathcal J^\prime} +y_i e_i +y \chi_J) \cdot \\
(\prod_{j \in J} \overline{\phi_j^*(2^{\kappa_j} y_j -\nu_j)}) dy^J)
(\prod_{j \in J} \phi_j(2^{\kappa_j} x_j -\nu_j)))
\overline {\phi_i^*(2^{\kappa_i} y_i -\nu_i)} dy_i)
\phi_i(2^{\kappa_i} x_i -\nu_i) = \\
\sum_{\nu_i \in \Z: (\delta_i B^1) \cap
(\supp \phi_i^*(2^{\kappa_i} \cdot -\nu_i)) \ne \emptyset} 2^{\kappa_i} \cdot \\
(\int_{(\delta_i B^1) \times (\delta B^d)^J \cap
(\supp \phi_i^*(2^{\kappa_i} \cdot -\nu_i)) \times
\prod_{j \in J} (\supp \phi_j^*(2^{\kappa_j} \cdot -\nu_j))}
f(x \chi_{\mathcal J^\prime} +y_i e_i +y \chi_J) \cdot \\
(\prod_{j \in J} \overline{\phi_j^*(2^{\kappa_j} y_j -\nu_j)})
\overline {\phi_i^*(2^{\kappa_i} y_i -\nu_i)} dy^J dy_i)
(\prod_{j \in J} \phi_j(2^{\kappa_j} x_j -\nu_j)) \phi_i(2^{\kappa_i} x_i -\nu_i),
\end{multline*}
где $ e_i = (0, \ldots, 0, 1_i, 0, \ldots, 0). $

Подставляя (1.4.6) в (1.4.5), приходим к (1.4.3) с $ \mathcal J $ вместо $ J: $
\begin{multline*}
((\prod_{j \in \mathcal J} V_j (E_{\kappa_j}^{j, p})) f)(x) =
\sum_{\nu^J \in (\Z^d)^J: (\supp \phi_j^*(2^{\kappa_j} \cdot -\nu_j))
\cap (\delta_j B^1) \ne \emptyset, j \in J} 2^{(\kappa^J, \e^J)} \cdot \\
(\sum_{\substack{\nu_i \in \Z: (\delta_i B^1) \cap \\
(\supp \phi_i^*(2^{\kappa_i} \cdot -\nu_i)) \ne \emptyset}} 2^{\kappa_i}
\biggl(\int_{\substack{(\delta_i B^1) \times (\delta B^d)^J \cap
(\supp \phi_i^*(2^{\kappa_i} \cdot -\nu_i)) \times\\
\prod_{j \in J} (\supp \phi_j^*(2^{\kappa_j} \cdot -\nu_j))}}
f(x \chi_{\mathcal J^\prime} +y_i e_i +y \chi_J) \cdot \\
(\prod_{j \in J} \overline{\phi_j^*(2^{\kappa_j} y_j -\nu_j)})
\overline {\phi_i^*(2^{\kappa_i} y_i -\nu_i)} dy^J dy_i)
(\prod_{j \in J} \phi_j(2^{\kappa_j} x_j -\nu_j)) \phi_i(2^{\kappa_i} x_i -\nu_i)) = \\
\sum_{\substack{\nu^J \in (\Z^d)^J, \nu_i \in \Z:  (\supp \phi_j^*(2^{\kappa_j} \cdot -\nu_j))
\cap (\delta_j B^1) \ne \emptyset, \\ j \in J; (\delta_i B^1) \cap
(\supp \phi_i^*(2^{\kappa_i} \cdot -\nu_i)) \ne \emptyset}} 2^{(\kappa^J, \e^J)} 2^{\kappa_i} \cdot\\
\biggl(\int_{\substack{(\delta_i B^1) \times (\delta B^d)^J \cap
(\supp \phi_i^*(2^{\kappa_i} \cdot -\nu_i)) \times\\
\prod_{j \in J} (\supp \phi_j^*(2^{\kappa_j} \cdot -\nu_j))}}
f(x \chi_{\mathcal J^\prime} +y_i e_i +y \chi_J) \cdot \\
(\prod_{j \in J} \overline{\phi_j^*(2^{\kappa_j} y_j -\nu_j)})
\overline {\phi_i^*(2^{\kappa_i} y_i -\nu_i)} dy^J dy_i\biggr)  \\
(\prod_{j \in J} \phi_j(2^{\kappa_j} x_j -\nu_j)) \phi_i(2^{\kappa_i} x_i -\nu_i) =
\sum_{\substack{\nu^{\mathcal J} \in (\Z^d)^{\mathcal J}:
(\supp \phi_j^*(2^{\kappa_j} \cdot -\nu_j)) \cap \\ (\delta_j B^1)
\ne \emptyset, j \in \mathcal J}} 2^{(\kappa^{\mathcal J}, \e^{\mathcal J})} \cdot \\
\biggl(\int_{(\delta B^d)^{\mathcal J} \cap \prod_{j \in \mathcal J}
(\supp \phi_j^*(2^{\kappa_j} \cdot -\nu_j))}
f(x \chi_{\mathcal J^\prime} +y \chi_{\mathcal J}) \cdot \\
(\prod_{j \in \mathcal J} \overline{\phi_j^*(2^{\kappa_j} y_j -\nu_j)}) dy^{\mathcal J}\biggr)
(\prod_{j \in \mathcal J} \phi_j(2^{\kappa_j} x_j -\nu_j)).
\end{multline*}

Тем самым справедливость (1.4.3) для произвольного непустого $ J \subset
\{1,\ldots,d\} $ установлена. В частности, при $ J = \{1,\ldots,d\} $
из (1.4.3) получаем, что для любой функции $ f \in L_p(\R^d) $ такой,
что при некотором $ \delta \in \R_+^d $ функция $ f $ равна $ 0 $ почти всюду
на $ \R^d \setminus (\delta B^d),$ соблюдается равенство
\begin{multline*}
((\prod_{j =1}^d V_j (E_{\kappa_j}^{j, p})) f)(x) =
\sum_{\nu \in \Z^d: (\supp \phi_j^*(2^{\kappa_j} \cdot -\nu_j))
\cap (\delta_j B^1) \ne \emptyset, j =1,\ldots,d} 2^{(\kappa, \e)} \cdot \\
(\int_{(\delta B^d) \cap \prod_{j =1}^d
(\supp \phi_j^*(2^{\kappa_j} \cdot -\nu_j))} f(y)
(\prod_{j =1}^d \overline{\phi_j^*(2^{\kappa_j} y_j -\nu_j)}) dy)
(\prod_{j =1}^d \phi_j(2^{\kappa_j} x_j -\nu_j)) = \\
\sum_{\nu \in \Z^d} 2^{(\kappa, \e)}
(\int_{(\delta B^d) \cap \prod_{j =1}^d
(\supp \phi_j^*(2^{\kappa_j} \cdot -\nu_j))} f(y)
(\prod_{j =1}^d \overline{\phi_j^*(2^{\kappa_j} y_j -\nu_j)}) dy)
(\prod_{j =1}^d \phi_j(2^{\kappa_j} x_j -\nu_j)) = \\
\sum_{\nu \in \Z^d} 2^{(\kappa, \e)}
(\int_{\R^d} f(y)
(\prod_{j =1}^d \overline{\phi_j^*(2^{\kappa_j} y_j -\nu_j)}) dy)
(\prod_{j =1}^d \phi_j(2^{\kappa_j} x_j -\nu_j)) = \\
\sum_{\nu \in \Z^d} 2^{(\kappa, \e)}
(\int_{\R^d} f(y) \overline{\phi^*(2^\kappa y -\nu)} dy)
\phi(2^\kappa x -\nu) =
(E_\kappa^{p} f)(x) \text{ (см. (1.2.19)).}
\end{multline*}
Из последнего равенства в силу непрерывности в $ L_p(\R^d) $ входящих
в него операторов и плотности в $ L_p(\R^d) $ множества функций
$ f \in L_p(\R^d), $ для которых $ f = f \chi_{\delta B^d} $ при некотором
$ \delta \in \R_+^d, $ вытекает (1.4.1).

Теперь проверим второе утверждение. Используя (1.4.1),
(1.3.3) и п. 2 леммы 1.3.2, имеем
\begin{multline*}
\mathcal E_\kappa^{p} = \sum_{\epsilon \in \Upsilon^d:
\s(\epsilon) \subset \s(\kappa)} (-\e)^\epsilon E_{\kappa
-\epsilon}^{p} =
\sum_{\epsilon \in \Upsilon^d: \s(\epsilon)
\subset \s(\kappa)}
(\prod_{j=1}^d (-1)^{\epsilon_j}) (\prod_{j=1}^d V_j(E_{\kappa_j -\epsilon_j}^{j, p})) =\\
=\sum_{\epsilon \in \Upsilon^d: \s(\epsilon) \subset \s(\kappa)}
(\prod_{j=1,\ldots,d: j \notin \s(\kappa)} (-1)^{\epsilon_j}
V_j(E_{\kappa_j -\epsilon_j}^{j, p})) (\prod_{j \in \s(\kappa)}
(-1)^{\epsilon_j} V_j(E_{\kappa_j -\epsilon_j}^{j, p})) =\\
=\sum_{\epsilon \in \Upsilon^d: \s(\epsilon) \subset \s(\kappa)}
(\prod_{j=1,\ldots,d: j \notin \s(\kappa)} V_j(E_0^{j, p}))
(\prod_{j \in \s(\kappa)} (-1)^{\epsilon_j} V_j(E_{\kappa_j -\epsilon_j}^{j, p})) =\\
=(\prod_{j=1,\ldots,d: j \notin \s(\kappa)} V_j(E_0^{j, p}))
\sum_{\epsilon \in \Upsilon^d: \s(\epsilon) \subset \s(\kappa)}
(\prod_{j \in \s(\kappa)} (-1)^{\epsilon_j} V_j(E_{\kappa_j -\epsilon_j}^{j, p})) =\\
=(\prod_{j=1,\ldots,d: j \notin \s(\kappa)} V_j(E_0^{j, p}))
\sum_{\epsilon^{\s(\kappa)} \in (\Upsilon^d)^{\s(\kappa)}}
(\prod_{j \in \s(\kappa)} (-1)^{\epsilon_j} V_j(E_{\kappa_j -\epsilon_j}^{j, p})) =\\
=(\prod_{j=1,\ldots,d: j \notin \s(\kappa)} V_j(E_0^{j, p}))
(\prod_{j \in \s(\kappa)} (\sum_{\epsilon_j =0,1}
(-1)^{\epsilon_j}
V_j(E_{\kappa_j -\epsilon_j}^{j, p}))) =\\
=(\prod_{j=1,\ldots,d: j \notin \s(\kappa)} V_j(E_0^{j, p}))
(\prod_{j \in \s(\kappa)} (V_j(E_{\kappa_j}^{j, p}) - V_j(E_{\kappa_j -1}^{j, p}))) =\\
=(\prod_{j=1,\ldots,d: j \notin \s(\kappa)} V_j(E_0^{j, p}))
(\prod_{j \in \s(\kappa)} V_j(E_{\kappa_j}^{j, p} - E_{\kappa_j
-1}^{j, p})). \square
\end{multline*}

Предложение 1.4.2

Пусть $ d \in \N $ и при $ j = 1,\ldots,d $ для функций $ \phi_j, \phi_j^* $
соблюдаются условия предложения 1.2.4, а функции $ \phi, \phi^* $
задаются равенствами
$$
\phi(x) = \prod_{j =1}^d \phi_j(x_j),
\phi^*(x) = \prod_{j =1}^d \phi_j^*(x_j).
$$
Тогда при $ \kappa \in \Z_+^d, 1 < p < \infty $ имеют место следующие
утверждения:

1) оператор $ E_\kappa^p = E_\kappa^{\phi, \phi^*} \mid_{L_p(\R^d)} $
является проектором в пространстве $ L_p(\R^d); $

2) для любой функции $ f \in L_p(\R^d) $ справедливо соотношение
\begin{equation*} \tag{1.4.7}
\| f -E_\kappa^p f \|_{L_p(\R^d)} \to 0 \text{при } \mn(\kappa) \to \infty;
\end{equation*}

3) для $ f \in L_p(\R^d), g \in L_{p^\prime}(\R^d) $
выполняются равенства
\begin{equation*} \tag{1.4.8}
\int_{\R^d} (E_\kappa f) \cdot \overline g dx = \int_{\R^d} f \cdot
\overline {(E_\kappa g)} dx;
\end{equation*}
и
\begin{equation*} \tag{1.4.9}
\int_{\R^d} (\mathcal E_\kappa f) \cdot \overline g dx =
\int_{\R^d} f \cdot \overline {(\mathcal E_\kappa g)} dx;
\end{equation*}

4) для $ \kappa^\prime \in \Z_+^d $ соблюдаются равенства
\begin{equation*} \tag{1.4.10}
\mathcal E_\kappa^{p} \mathcal E_{\kappa^\prime}^{p} =
\begin{cases} \mathcal E_\kappa^{p},
  \text{при $ \kappa = \kappa^\prime $}; \\
       0,  \text{при $ \kappa \ne \kappa^\prime $}.
\end{cases}
\end{equation*}

Доказательство.

Проверяя первое утверждение предложения, на основании (1.4.1), (1.3.3),
п. 2 лемммы 1.3.2 и (1.2.34) имеем
\begin{multline*}
E_\kappa^{p} E_\kappa^{p} = (\prod_{j=1}^d V_j (E_{\kappa_j}^{j, p}))
(\prod_{j=1}^d V_j (E_{\kappa_j}^{j, p})) =
\prod_{j=1}^d (V_j (E_{\kappa_j}^{j, p})  V_j (E_{\kappa_j}^{j, p})) =\\
\prod_{j=1}^d V_j (E_{\kappa_j}^{j, p}  E_{\kappa_j}^{j, p}) =
\prod_{j=1}^d V_j (E_{\kappa_j}^{j, p}) = E_\kappa^{p},
\end{multline*}
т.е. первое утверждение верно.

Для вывода (1.4.7) понадобятся некоторые факты. Прежде всего заметим, что
\begin{equation*} \tag{1.4.11}
\span \{ \phi(2^\kappa \cdot -\nu), \nu \in \Z^d \} \subset
\Im E_\kappa^{p}.
\end{equation*}

В самом деле, пользуясь (1.2.19), при $ \nu \in \Z^d $ получаем
\begin{multline*}
(E_\kappa^{p} (\phi(2^\kappa \cdot -\nu)))(x) =
\sum_{\nu^\prime \in \Z^d} 2^{(\kappa, \e)}
(\int_{\R^d} \phi(2^\kappa y -\nu) \overline {\phi^*(2^\kappa y -
\nu^\prime)} dy) \phi(2^\kappa x -\nu^\prime).
\end{multline*}
А поскольку для $ \nu, \nu^\prime \in \Z^d $ в силу (1.2.8) и теоремы Фубини
\begin{multline*}
\int_{\R^d} \phi(2^\kappa y -\nu) \overline {\phi^*(2^\kappa y -
\nu^\prime)} dy =
2^{-(\kappa, \e)} \int_{\R^d} \phi(y -\nu) \overline {\phi^*(y -
\nu^\prime)} dy = \\
2^{-(\kappa, \e)} \int_{\R^d} (\prod_{j =1}^d
\phi_j(y_j -\nu_j)) \overline { (\prod_{j =1}^d \phi_j^*(y_j -
\nu_j^\prime))} dy = \\
2^{-(\kappa, \e)} \prod_{j =1}^d
(\int_{\R} \phi_j(y_j -\nu_j) \overline { \phi_j^*(y_j -\nu_j^\prime)} dy_j) =
\begin{cases} 2^{-(\kappa, \e)}, \text{ при } \nu^\prime = \nu; \\
0, \text{ при } \nu^\prime \ne \nu,
\end{cases}
\end{multline*}
то
$$
(E_\kappa^{p} (\phi(2^\kappa \cdot -\nu)))(x) = \phi(2^\kappa x -\nu),
\nu \in \Z^d.
$$
Отсюда в силу линейности оператора $ E_\kappa^{p} $ вытекает (1.4.11).

Далее, покажем, что при $ \kappa,
\kappa^\prime \in \Z_+^d: \kappa^\prime \le \kappa, $ имеет место
включение
\begin{equation*} \tag{1.4.12}
\span \{ \phi(2^{\kappa^\prime} \cdot -\nu^\prime), \nu^\prime \in \Z^d \}
\subset \span \{ \phi(2^{\kappa} \cdot -\nu), \nu \in \Z^d \}.
\end{equation*}

Действительно, ввиду (1.2.33) при $ \kappa^\prime \le \kappa,
\nu^\prime \in \Z^d $ имеем
\begin{multline*}
\phi(2^{\kappa^\prime} x -\nu^\prime) = \prod_{j =1}^d
\phi_j(2^{\kappa_j^\prime} x_j -\nu_j^\prime) =
\prod_{j =1}^d
(\sum_{ \nu_j \in \Nu_j} c_{\nu_j}^j \phi_j(2^{\kappa_j} x_j -\nu_j)) = \\
\sum_{\nu \in \Z^d: \nu_j \in \Nu_j, j=1,\ldots,d}
(\prod_{j =1}^d (c_{\nu_j}^j \phi_j(2^{\kappa_j} x_j -\nu_j))) =
\sum_{\nu \in \Z^d: \nu_j \in \Nu_j, j=1,\ldots,d} c_\nu \phi(2^\kappa x -\nu)
\in \\ \span \{ \phi(2^{\kappa} \cdot -\nu), \nu \in \Z^d \}.
\end{multline*}
Отсюда следует (1.4.12).

Проверим еще, что
\begin{equation*} \tag{1.4.13}
\cup_{ \kappa \in \Z_+^d} \span \{ \phi(2^{\kappa} \cdot -\nu),
\nu \in \Z^d \} = \span \{ \phi(2^{\kappa} \cdot -\nu), \nu \in \Z^d,
\kappa \in \Z_+^d \}.
\end{equation*}

Ясно, что
$$
\cup_{ \kappa \in \Z_+^d} \span \{ \phi(2^{\kappa} \cdot -\nu),
\nu \in \Z^d \} \subset \span \{ \phi(2^{\kappa} \cdot -\nu), \nu \in \Z^d,
\kappa \in \Z_+^d \}.
$$
Установим обратное включение. Пусть $ f \in \span \{ \phi(2^{\kappa} \cdot -\nu),
\nu \in \Z^d, \kappa \in \Z_+^d \}, $ т.е.
$ f = \sum_{i=1}^m c^i \phi(2^{\kappa^i} \cdot -\nu^i), $ где $ c^i \in \C,
\kappa^i \in \Z_+^d, \nu^i \in \Z^d, i =1,\ldots,m. $
Положим $ \kappa_j = \max_{i=1,\ldots,m} \kappa_j^i \in \Z_+, j =1,\ldots,d. $
Тогда, учитывая, что благодаря (1.4.12), для $ i =1,\ldots,m $ справедливо
включение
$$
\phi(2^{\kappa^i} \cdot -\nu^i) \in
\span \{ \phi(2^{\kappa^i} \cdot -\nu), \nu \in \Z^d \}
\subset \span \{ \phi(2^{\kappa} \cdot -\nu), \nu \in \Z^d \},
$$
видим, что
$ f \in \span \{ \phi(2^{\kappa} \cdot -\nu), \nu \in \Z^d \}, $
а, следовательно, выполняется (1.4.13).

Теперь покажем по индукции относительно $ d, $ что
\begin{equation*} \tag{1.4.14}
\close_{L_p(\R^d)} \span \{ \phi(2^{\kappa} \cdot -\nu), \nu \in \Z^d,
\kappa \in \Z_+^d \} = L_p(\R^d).
\end{equation*}

При $ d =1 $ в силу (1.2.40), (1.2.27) равенство (1.4.14) выполняется.
Предположим, что (1.4.14) справедливо при $ d -1 $ вместо $ d. $
И предположим, что
$$
\close_{L_p(\R^d)} \span \{ \phi(2^{\kappa} \cdot -\nu), \nu \in \Z^d,
\kappa \in \Z_+^d \} \ne L_p(\R^d).
$$
Тогда на основании одного из следствий из теоремы Хана-Банаха существует
ненулевой непрерывный линейный функционал на $ L_p(\R^d), $
аннулирующий
$ \close_{L_p(\R^d)} \span \{ \phi(2^{\kappa} \cdot -\nu), \nu \in \Z^d,
\kappa \in \Z_+^d \}. $ Учитывая, что $ (L_p(\R^d))^* = L_{p^\prime}(\R^d), $
заключаем, что существует функция $ g \in L_{p^\prime}(\R^d), $ отличная
от $ 0 $ на множестве положительной меры, такая, что при
$ \kappa \in \Z_+^d, \nu \in \Z^d $ верно равенство
$$
\int_{ \R^d} \phi(2^{\kappa} \cdot -\nu)  g dx =0.
$$
Поэтому по теореме Фубини при каждом $ (\kappa_1, \ldots, \kappa_{d -1})
\in \Z_+^{d -1}, (\nu_1, \ldots, \nu_{d -1}) \in \Z^{d -1} $
для любых $ \kappa_d \in \Z_+ $ и $ \nu_d \in \Z $ соблюдается равенство
\begin{multline*}
0 = \int_{ \R^d} \phi(2^{\kappa} \cdot -\nu) g dx =
\int_{ \R^d} (\prod_{j=1}^d \phi_j(2^{\kappa_j} x_j -\nu_j))
 g(x) dx = \\
\int_\R \int_{ \R^{d -1}} (\prod_{j=1}^d \phi_j(2^{\kappa_j} x_j -\nu_j))
g(x_1, \ldots,x_{d -1}, x_d) dx_1 \ldots dx_{d -1} dx_d = \\
\int_\R \phi_d(2^{\kappa_d} x_d -\nu_d)
( \int_{ \R^{d -1}} (\prod_{j=1}^{d -1} \phi_j(2^{\kappa_j} x_j -\nu_j))
g(x_1, \ldots,x_{d -1}, x_d) dx_1 \ldots dx_{d -1} ) dx_d.
\end{multline*}
Из этого равенства, принимая во внимание соблюдение
(1.4.14) при $ d =1, $ а также включение
$$
( \int_{ \R^{d -1}} (\prod_{j=1}^{d -1} \phi_j(2^{\kappa_j} x_j -\nu_j))
g(x_1, \ldots,x_{d -1}, x_d) dx_1 \ldots dx_{d -1} ) \in
L_{p^\prime}(\R),
$$
заключаем, что при каждом $ (\kappa_1, \ldots, \kappa_{d -1})
\in \Z_+^{d -1}, (\nu_1, \ldots, \nu_{d -1}) \in \Z^{d -1} $
для любой функции $ f \in L_p(\R) $ имеет место равенство
$$
\int_\R f(x_d) ( \int_{ \R^{d -1}} (\prod_{j=1}^{d -1} \phi_j(2^{\kappa_j}
x_j -\nu_j)) g(x_1, \ldots,x_{d -1}, x_d) dx_1 \ldots
dx_{d -1} ) dx_d =0,
$$
и, значит, при каждом $ (\kappa_1, \ldots, \kappa_{d -1})
\in \Z_+^{d -1}, (\nu_1, \ldots, \nu_{d -1}) \in \Z^{d -1} $ почти для всех
$ x_d \in \R $ соблюдается равенство
$$
\int_{ \R^{d -1}} (\prod_{j=1}^{d -1} \phi_j(2^{\kappa_j}
x_j -\nu_j)) g(x_1, \ldots,x_{d -1}, x_d) dx_1 \ldots
dx_{d -1} =0.
$$
Из сказанного вытекает, что почти для всех
$ x_d \in \R $ при любых $ (\kappa_1, \ldots, \kappa_{d -1})
\in \Z_+^{d -1}, (\nu_1, \ldots, \nu_{d -1}) \in \Z^{d -1} $ выполняется
равенство
$$
\int_{ \R^{d -1}} (\prod_{j=1}^{d -1} \phi_j(2^{\kappa_j}
x_j -\nu_j)) g(x_1, \ldots,x_{d -1}, x_d) dx_1 \ldots
dx_{d -1} =0.
$$

Учитывая, что в соответствии с предположением индукции соблюдается
(1.4.14) с $ d -1 $ вместо $ d, $ как и выше, приходим к выводу, что почти
для всех $ x_d \in \R $ для любой функции $ f \in L_p(\R^{d -1}) $
верно равенство
$$
\int_{ \R^{d -1}} f(x_1, \ldots, x_{d -1})
g(x_1, \ldots,x_{d -1}, x_d) dx_1 \ldots
dx_{d -1} =0,
$$
а, следовательно, почти для всех $ x_d \in \R $ почти всюду в $ \R^{d -1} $
имеет место равенство $ g(x_1, \ldots, x_{d -1}, x_d) =0, $
или $ g(x) =0 $ почти всюду в $ \R^d, $
что противоречит выбору $ g. $
Тем самым, установлена справедливость (1.4.14) при любом $ d. $

Переходя к выводу (1.4.7), для $ f \in L_p(\R^d) $ и произвольного
$ \epsilon > 0 $ в силу (1.4.14) выберем функцию
$$
g \in \span \{ \phi(2^{\kappa} \cdot -\nu), \nu \in \Z^d,
\kappa \in \Z_+^d\},
$$
для которой
$$
\| f - g \|_{L_p(\R^d)} < \epsilon.
$$
Согласно (1.4.13) существует $ \kappa^0 \in \Z_+^d $ такое, что
$$
g \in \span \{ \phi(2^{\kappa^0} \cdot -\nu), \nu \in \Z^d \}.
$$
Тогда, замечая, что для $ \kappa \in \Z_+^d: \kappa \ge \kappa^0, $
вследствие (1.4.12) имеет место включение
$ g \in \span \{ \phi(2^{\kappa} \cdot -\nu), \nu \in \Z^d \}, $
и принимая во внимание (1.4.11), тот факт, что $ E_\kappa^{p} $ --
проектор, а также (1.2.20), для $ \kappa \in \Z_+^d: \mn(\kappa) >
\mx(\kappa^0), $ получаем
\begin{multline*}
\| f -E_\kappa^{p} f \|_{L_p(\R^d)} =
\| f -g +g -E_\kappa^{p} f \|_{L_p(\R^d)} = \\
\| f -g +E_\kappa^{p} g -E_\kappa^{p} f \|_{L_p(\R^d)} \le
\| f -g \|_{L_p(\R^d)} +\| E_\kappa^{p} (g -f) \|_{L_p(\R^d)} \le \\
\| f -g \|_{L_p(\R^d)} +c_1 \| g -f \|_{L_p(\R^d)} =
(1 +c_1) \| f -g \|_{L_p(\R^d)} < c_2 \epsilon,
\end{multline*}
что влечет (1.4.7).

Для доказательства (1.4.8) сначала установим, что при $ \kappa \in \Z_+^d $
для $ f \in L_p(\R^d), g \in L_{p^\prime}(\R^d) $ для любого непустого
множества $ J \subset \{1,\ldots,d\} $ выполняется равенство
\begin{multline*} \tag{1.4.15}
\int_{\R^d} ((\prod_{j \in J} V_j(E_{\kappa_j}^{j, p}))f) \overline g dx =
\int_{\R^d} f \cdot \overline {((\prod_{j \in J}
V_j(E_{\kappa_j}^{j, p^\prime}))g)} dx.
\end{multline*}

Вывод (1.4.15) проведем по индукции относительно $ \card J. $
При $ \card J =1, $ т.е. при $ J = \{j\}, j =1, \ldots,d, $ в силу теоремы
Фубини, (1.3.1), (1.2.29) имеем
\begin{multline*}
 \int_{\R^d} ((V_j(E_{\kappa_j}^{j, p}))f)(x)
\overline {g(x)} dx =
 \int_{\R^{d -1}} \int_\R  ((V_j(E_{\kappa_j}^{j, p}))f)(x_1,\ldots, x_{j -1},
x_j, x_{j +1}, \ldots, x_d)\cdot\\ \overline {g(x_1,\ldots, x_{j -1}, x_j,
x_{j +1}, \ldots, x_d)} dx_j dx_1 \ldots dx_{j -1} dx_{j +1} \ldots dx_d = \\
 \int_{\R^{d -1}} \int_\R  ( E_{\kappa_j}^{j, p} f(x_1,\ldots, x_{j -1},
\cdot, x_{j +1}, \ldots, x_d))(x_j)\cdot\\ \overline {g(x_1,\ldots, x_{j -1}, x_j,
x_{j +1}, \ldots, x_d)} dx_j dx_1 \ldots dx_{j -1} dx_{j +1} \ldots dx_d = \\
 \int_{\R^{d -1}} \int_\R f(x_1,\ldots, x_{j -1}, x_j, x_{j +1}, \ldots, x_d)\cdot\\
\overline {( E_{\kappa_j}^{j, p^\prime} g(x_1,\ldots, x_{j -1},\cdot, x_{j +1},
\ldots, x_d))(x_j) } dx_j dx_1 \ldots dx_{j -1} dx_{j +1} \ldots dx_d = \\
 \int_{\R^{d -1}} \int_\R f(x_1,\ldots, x_{j -1}, x_j, x_{j +1}, \ldots, x_d)\cdot\\
\overline {( V_j( E_{\kappa_j}^{j, p^\prime}) g)(x_1,\ldots, x_{j -1}, x_j,
x_{j +1}, \ldots, x_d) } dx_j dx_1 \ldots dx_{j -1} dx_{j +1} \ldots dx_d = \\
 \int_{\R^d} f(x) \overline{(V_j( E_{\kappa_j}^{j, p^\prime}) g)(x)} dx,
\end{multline*}
что совпадает с (1.4.15) при $ J = \{j\}, j =1,\ldots,d. $

Предположим, что (1.4.15) справедливо для любого множества $ J \subset
\{1,\ldots,d\}, $  для которого $ \card J \le m \le d -1. $ Покажем, что тогда
(1.4.15) имеет место для любого множества $ \mathcal J \subset
\{1,\ldots,d\}, $ для которого $ \card \mathcal J = m +1. $
Представляя $ \mathcal J $  в виде $ \mathcal J = J \cup \{i\}, $
где $ i \notin J, $ и пользуясь (1.3.3), предположением индукции, получаем
\begin{multline*}
\int_{\R^d} ((\prod_{j \in \mathcal J} V_j(E_{\kappa_j}^{j, p}))f) \overline g dx =
\int_{\R^d} (V_i(E_{\kappa_i}^{i, p}))
((\prod_{j \in J} V_j(E_{\kappa_j}^{j, p}))f) \overline g dx = \\
\int_{\R^d} ((\prod_{j \in J} V_j(E_{\kappa_j}^{j, p}))f)
\overline {((V_i(E_{\kappa_i}^{i, p^\prime})) g) } dx =
\int_{\R^d} f \cdot \overline {((\prod_{j \in J}
V_j(E_{\kappa_j}^{j, p^\prime})) ((V_i(E_{\kappa_i}^{i, p^\prime})) g)) } dx = \\
\int_{\R^d} f \cdot \overline {((\prod_{j \in \mathcal J}
V_j(E_{\kappa_j}^{j, p^\prime}))g)} dx,
\end{multline*}
что завершает вывод (1.4.15).

Полагая в (1.4.15) $ J = \{1,\ldots,d\}, $ и учитывая (1.4.1),
приходим к (1.4.8).

Проверяя (1.4.9), для $ f \in L_p(\R^d), g \in L_{p^\prime}(\R^d) $
ввиду (1.4.8) получаем
\begin{multline*}
\int_{\R^d} (\mathcal E_\kappa f) \overline {g} dx = \!\int_{\R^d}
(\sum_{\epsilon \in \Upsilon^d: \s(\epsilon) \subset \s(\kappa)}
\!(-\e)^\epsilon E_{\kappa -\epsilon} f) \overline {g} dx =\\
=\!\int_{\R^d} \sum_{\epsilon \in \Upsilon^d:\s(\epsilon) \subset
\s(\kappa)} \!(-\e)^\epsilon (E_{\kappa -\epsilon} f) \overline {g} dx =
\sum_{\epsilon \in \Upsilon^d: \s(\epsilon) \subset \s(\kappa)}
\!(-\e)^\epsilon \!\int_{\R^d} (E_{\kappa -\epsilon} f) \overline {g} dx =\\
=\sum_{\epsilon \in \Upsilon^d: \s(\epsilon) \subset \s(\kappa)}
\!(-\e)^\epsilon \!\int_{\R^d} f \cdot \overline
{(E_{\kappa -\epsilon} g)} dx
= \!\int_{\R^d} \sum_{\epsilon \in \Upsilon^d: \s(\epsilon) \subset
\s(\kappa)}
\!(-\e)^\epsilon f \cdot \overline {(E_{\kappa -\epsilon} g)} dx =\\
=\!\int_{\R^d} f (\sum_{\epsilon \in \Upsilon^d:\s(\epsilon)
\subset \s(\kappa)} \!(-\e)^\epsilon \overline
{(E_{\kappa -\epsilon} g)}) dx = \\
=\!\int_{\R^d} f \overline {(\sum_{\epsilon \in \Upsilon^d: \s(\epsilon)
\subset \s(\kappa)} \!(-\e)^\epsilon
(E_{\kappa -\epsilon} g))} dx =
\int_{\R^d} f \overline {(\mathcal E_\kappa g)} dx.
\end{multline*}

Наконец, убедимся в справедливости (1.4.10). При $ d =1 $ соотношение
(1.4.10) совпадает с (1.2.35).
Установим соблюдение (1.4.10) при произвольном $ d \in \N. $
Согласно (1.4.2), (1.3.3), п. 2 леммы 1.3.2 имеем
$$
\mathcal E_\kappa^{p} \mathcal E_{\kappa^\prime}^{p} \!=\!
(\prod_{j=1}^d \!V_j(\mathcal E_{\kappa_j}^{j, p}))
(\prod_{j=1}^d \!V_j(\mathcal E_{\kappa_j^\prime}^{j, p})) \!=\!
\prod_{j=1}^d \!(V_j(\mathcal E_{\kappa_j}^{j, p}) V_j(\mathcal
E_{\kappa_j^\prime}^{j, p})) \!=\! \prod_{j=1}^d \!V_j(\mathcal
E_{\kappa_j}^{j, p} \mathcal E_{\kappa_j^\prime}^{j, p}).
$$
Отсюда в силу равенства (1.2.35)
с учетом (1.4.2)
вытекает справедливость (1.4.10) при произвольном $ d \in \N. $
$\square$
\bigskip

\centerline{ \S 2. Теорема типа Литтлвуда-Пэли}
\centerline{для операторов $ \mathcal E_\kappa $}
\bigskip

2.1. В этом пункте приведены вспомогательные утверждения,
на которые опирается доказательство занимающей центральное место в работе
леммы 2.2.1. Все эти утверждения можно извлечь из [11, гл. I].

Будем обозначать $ |x| = (\sum_{j=1}^d |x_j|^2)^{1/2}, x \in \R^d.$

Для замкнутого множества $ F \subset \R^d $ и $ x \in \R^d $
положим
$$
\rho(x) = \rho(x,F) = \min_{y \in F} | x -y|.
$$

Лемма 2.1.1

Пусть $ d \in \N. $ Тогда существует константа $ c_1(d) >0 $
такая, что для любого замкнутого множества $ F \subset \R^d, $
для которого $ \mes (\R^d \setminus F) < \infty, $ функция $ M(x), $
определяемая для $ x \in F $ равенством
$$
M(x) = \int_{\R^d} \rho(u) |x -u|^{-(d+1)} du,
$$
суммируема на $ F $ и соблюдается  неравенство
\begin{equation*} \tag{2.1.1}
\int_F M(x) dx \le c_1 \mes (\R^d \setminus F).
\end{equation*}
(см. \S 2 гл. I из [11])

Предложение 2.1.2

Пусть $ d \in \N. $ Тогда существуют константы $ c_2(d) >0, c_3(d) >0,
c_4(d) >0, c_5(d) >0 $ такие, что
для любой (вещественной) функции $ f \in L_1(\R^d) $ при любом $ \alpha \in \R_+ $
можно построить замкнутое множесво $ F \subset \R^d $ и семейство открытых кубов
$ \{ Q_r, r \in \N \} $ со следующими свойствами:

1) почти для всех $ x \in F $  выполняется неравенство
\begin{equation*} \tag{2.1.2}
| f(x) | \le \alpha;
\end{equation*}

2) для $ W = \R^d \setminus F $ справедливо неравенство
\begin{equation*} \tag{2.1.3}
\mes W \le (c_2 / \alpha) \int_{\R^d} |f| dx;
\end{equation*}

3) \begin{equation*} \tag{2.1.4}
Q_r \cap Q_s = \emptyset, \text{ для}   r,s \in \N: r \ne s;
\end{equation*}

4) \begin{equation*} \tag{2.1.5}
W = \cup_{ r \in \N} \overline Q_r,
\end{equation*}
где  $ \overline Q_r $ -- замыкание куба $ Q_r, r \in \N; $

5) при $ r \in \N $ справедливы неравенства
\begin{equation*} \tag{2.1.6}
c_3 \diam Q_r < \inf_{ x \in Q_r} \rho(x, F) \le
c_4 \diam Q_r,
\end{equation*}
где  $ \diam Q = \sup_{x,y \in Q} |x -y|, Q \subset \R^d; $

а также

6) при  $ r \in \N $ имеет место оценка
\begin{equation*} \tag{2.1.7}
(1 / \mes Q_r) \int_{ Q_r} | f(x) | dx \le c_5 \alpha.
\end{equation*}
(см. \S 3 гл. I из [11])

При $ d \in \N $ через $ L(\R^d) $ обозначим пространство вещественных
измеримых по Лебегу функций в $ \R^d. $ Как обычно, при $ 1 \le p_1, p_2 \le
\infty $ под суммой $ L_{p_1}(\R^d) +L_{p_2}(\R^d) $ понимается
подпространство в $ L(\R^d), $ состоящее из всех функций $ f \in L(\R^d), $
для которых существуют функции $ f_1 \in L_{p_1}(\R^d) $ и
$ f_2 \in L_{p_2}(\R^d) $ такие, что $ f = f_1 +f_2.$

Напомним (см., например, [11]), что при $ 1 \le p_1 \le p \le p_2 < \infty $ справедливо включение
$ L_p(\R^d) \subset L_{p_1}(\R^d) +L_{p_2}(\R^d). $

Теорема 2.1.3

Пусть $ d \in \N, 1 < q < \infty, C_0 \in \R_+, C_1 \in \R_+ $ и
$ T: (L_1(\R^d) +L_q(\R^d)) \mapsto L(\R^d) $ --
отображение, удовлетворяющее следующим условиям:

1) для любых $ f, g \in (L_1(\R^d) +L_q(\R^d)) $ почти для всех
$ x \in \R^d $ выполняется неравенство
  \begin{equation*} \tag{2.1.8}
| (T(f +g))(x) | \le | (T f)(x) | +| (T g)(x) |;
\end{equation*}

2) для $ f \in L_1(\R^d) $ при $ \alpha >0 $ соблюдается неравенство
\begin{equation*} \tag{2.1.9}
\mes \{ x \in \R^d: | (T f)(x) | > \alpha \} \le
(C_0 / \alpha) \| f \|_{L_1(\R^d)};
\end{equation*}

3) для $ f \in L_q(\R^d) $ при $ \alpha >0 $ выполняется неравенство
\begin{equation*} \tag{2.1.10}
\mes \{ x \in \R^d: | (T f)(x) | > \alpha \} \le
((C_1 / \alpha) \| f \|_{L_q(\R^d)})^q.
\end{equation*}
Тогда при $ 1< p < q $ существует константа $ c_6(p, q, C_0, C_1)
>0 $ такая, что для любого отображения $ T: (L_1(\R^d) +L_q(\R^d))
\mapsto L(\R^d), $ подчиненного условиям (2.1.8) -- (2.1.10), для
$ f \in L_p(\R^d) $ имеет место неравенство
\begin{equation*} \tag{2.1.11}
\| T f \|_{L_p(\R^d)} \le c_6 \| f \|_{L_p(\R^d)}.
\end{equation*}
(см. \S 4 гл. I из [11])
\bigskip

2.2. В этом пункте устанавливается аналог теоремы Литтлвуда-Пэли для
ортопроекторов на подпространства всплесков
$ \{\mathcal E_\kappa, \kappa \in \Z_+^d\} $ (см. теорему 2.2.3),
соответствующие кратно-масштабному анализу, порожденному тензорным
произведением гладких финитных функций. При этом при
доказательстве теоремы 2.2.3 будем придерживаться того же подхода,
что в [3, п. 1.5.2] в случае теоремы Литтлвуда-Пэли для кратных рядов Фурье.
Убедимся, что имеет место

Лемма  2.2.1

Пусть функции $ \phi, \phi^* $ подчинены условиям предложения 1.2.4,
$ \phi^* \in C_1(\R) $ и $ 1 < p < \infty. $ Тогда существует константа
$ c_1(\phi, \phi^*, p) >0 $ такая, что при любом $ k \in \Z_+ $ для любого
набора чисел $ \sigma = \{ \sigma_\kappa \in \{-1, 1\}: \kappa =0, \ldots, k \}, $
для $ f \in L_p(\R) $ справедливо неравенство
\begin{equation*} \tag{2.2.1}
\| \sum_{\kappa =0}^k \sigma_\kappa \cdot (\mathcal E_\kappa f)
\|_{L_p(\R)} \le c_1 \| f \|_{L_p(\R)}.
\end{equation*}

Отметим, что если $ \phi_N^D $ -- масштабирующая функция Добеши
порядка $ N $ (см., например, п. 7.3 из [6] или [7]),  то при достаточно
больших $ N \in \N $ пара функций $ \phi = \phi_N^D, \phi^* = \phi_N^D $
удовлетворяет условиям леммы 2.2.1.

Отметим также, что доказательство леммы 2.2.1 проводится по схеме, использованной
в [11] при доказательстве теоремы 1 из гл. II.

Доказательство.

Сначала установим справедливость (2.2.1) при $ 1 < p \le 2. $
Неравенство (2.2.1) достаточно получить для функций $ f, $
принимающих вещественные значения.
В самом деле, если (2.2.1) верно для вещественнозначных функий $ f, $
то для комплекснозначной функции $ f $  имеем
\begin{multline*}
\| \sum_{\kappa =0}^k \sigma_\kappa \cdot (\mathcal E_\kappa f)
\|_{L_p(\R)} =
\| \sum_{\kappa =0}^k \sigma_\kappa \cdot (\mathcal E_\kappa ( \Re f
+i \Im f)) \|_{L_p(\R)} = \\
\| \sum_{\kappa =0}^k \sigma_\kappa \cdot (\mathcal E_\kappa ( \Re f))
+i \sum_{\kappa =0}^k \sigma_\kappa \cdot (\mathcal E_\kappa
( \Im f)) \|_{L_p(\R)}  \le \\
\| \sum_{\kappa =0}^k \sigma_\kappa \cdot (\mathcal E_\kappa
( \Re f)) \|_{L_p(\R)}
+\| \sum_{\kappa =0}^k \sigma_\kappa \cdot (\mathcal E_\kappa
( \Im f)) \|_{L_p(\R)} \le \\
c_1 \| \Re f \|_{L_p(\R)} +c_1 \| \Im f \|_{L_p(\R)}
\le 2 c_1 \| f \|_{L_p(\R)}.
\end{multline*}

Определим при $ k \in \Z_+, \sigma = \{ \sigma_\kappa \in \{-1, 1\}:
\kappa =0, \ldots, k \} $ отображение
$ T = T_{k,\sigma}: L_1(\R) +L_2(\R) \mapsto L(\R), $
(где $ L_1(\R), L_2(\R), L(\R) $ -- пространства
вещественнозначных функций), полагая для
$ f \in (L_1(\R) +L_2(\R)) \subset L_1^{\loc(\R)} $ значение
$$
( T f)(x) = | \sum_{\kappa =0}^k \sigma_\kappa (\mathcal E_\kappa ( f ))(x) |.
$$

При $ k \in \Z_+, \sigma = \{ \sigma_\kappa \in \{-1, 1\}:
\kappa =0, \ldots, k \} $
для $ f \in (L_1(\R) +L_2(\R)), g \in (L_1(\R) +L_2(\R)) $ имеем
\begin{multline*} \tag{2.2.2}
| (T (f +g))(x) | =
| \sum_{\kappa =0}^k \sigma_\kappa (\mathcal E_\kappa (f +g))(x) | = \\
| \sum_{\kappa =0}^k \sigma_\kappa (\mathcal E_\kappa (f))(x) +
\sum_{\kappa =0}^k \sigma_\kappa (\mathcal E_\kappa ( g))(x) | \le \\
| \sum_{\kappa =0}^k \sigma_\kappa (\mathcal E_\kappa (f))(x) | +
| \sum_{\kappa =0}^k \sigma_\kappa (\mathcal E_\kappa ( g))(x) | = \\
| (T f)(x) | +| (T g)(x) |,
\end{multline*}
т.е. выполняется (2.1.8).

Далее, покажем, что при $ k \in \Z_+, \sigma = \{ \sigma_\kappa \in
\{-1, 1\}: \kappa =0, \ldots, k \}, $ для $ f \in L_2(\R) $ и $ \alpha >0 $
соблюдается неравенство
\begin{equation*} \tag{2.2.3}
\mes \{ x \in \R: | (T f)(x) | > \alpha \} \le
((1 / \alpha) \| f \|_{L_2(\R)})^2.
\end{equation*}

В самом деле, при $ k \in \Z_+, \sigma = \{ \sigma_\kappa \in \{-1, 1\}:
\kappa =0, \ldots, k \}, $ для $ f \in L_2(\R) $ и $ \alpha >0 $ в виду
(1.2.36), (1.2.38) имеем
\begin{multline*}
\int_\R | (T f)(x)|^2 dx = \int_\R | \sum_{\kappa =0}^k \sigma_\kappa
(\mathcal E_\kappa(f))(x)|^2 dx = \\
\int_\R (\sum_{\kappa =0}^k
\sigma_\kappa (\mathcal E_\kappa( f ))(x))
\overline {(\sum_{\kappa^\prime =0}^k
\sigma_{\kappa^\prime} (\mathcal E_{\kappa^\prime}( f))(x))} dx = \\
\int_\R \sum_{\kappa =0}^k \sum_{\kappa^\prime =0}^k
\sigma_\kappa (\mathcal E_\kappa( f ))(x)
\sigma_{\kappa^\prime} \overline {(\mathcal E_{\kappa^\prime}( f))(x)} dx = \\
\sum_{\kappa =0}^k \sum_{\kappa^\prime =0}^k
\sigma_\kappa \sigma_{\kappa^\prime}
\int_\R (\mathcal E_\kappa( f))(x)
\overline {(\mathcal E_{\kappa^\prime} ( f))(x)} dx \\
= \sum_{\kappa =0}^k \int_\R
| (\mathcal E_\kappa( f ))(x)|^2 dx
= \sum_{ \kappa =0}^k \| \mathcal E_\kappa f \|_{L_2(\R)}^2
\le \sum_{\kappa \in \Z_+} \| \mathcal E_\kappa f \|_{L_2(\R)}^2
= \| f \|_{L_2(\R)}^2,
\end{multline*}
откуда, как обычно, получаем
\begin{multline*}
\alpha^2 \mes \{ x \in \R: | (T f)(x) | > \alpha \} = \int_{ \{
x \in \R: | (T f)(x) | > \alpha \} } \alpha^2 dx \\
\le \int_{ \{
x \in \R: | (T f)(x) | > \alpha \} } | (T f )(x)|^2 dx \le
\int_\R | (T f)(x)|^2 dx \le \| f \|_{L_2(\R)}^2,
\end{multline*} и, значит, верно (2.2.3).

Теперь установим, что существует константа $ C_0 >0 $ такая,
что при $ k \in \Z_+, \sigma = \{ \sigma_\kappa \in \{-1, 1\}:
\kappa =0, \ldots, k \} $ для $ f \in L_1(\R) $ и $ \alpha >0 $
соблюдается неравенство
\begin{equation*} \tag{2.2.4}
\mes \{ x \in \R: | (T f)(x) | > \alpha \} \le
(C_0 / \alpha) \| f \|_{L_1(\R)}.
\end{equation*}

Пусть $ k \in \Z_+, \sigma = \{ \sigma_\kappa \in \{-1, 1\}:
\kappa =0, \ldots, k \}, f \in L_1(\R) $ и $ \alpha >0. $
Для функции $ f $ и числа $ \alpha $ построим замкнутое множество
$ F,$ множество $ W = \R \setminus F $ и семейство интервалов
$ \{ Q_r, r \in \N \}, $ для которых соблюдаются условия
(2.1.2) -- (2.1.7) при $ d =1. $

Определим функции $ g \in L_1(\R) \cap
L_2(\R) $ и $ h \in L_1(\R), $ полагая
$$
g(x) = f(x) \chi_F(x) +\sum_{r=1}^\infty (1/ \mes Q_r) (\int_{Q_r} f(y) dy)
\chi_{Q_r}(x), x \in \R,
$$
и $ h = f -g. $

Из (2.1.5), (2.1.4) с учетом того, что $ \mes (W \setminus (\cup_{
r \in \N } Q_r)) =0 $ (ибо $ (W \setminus (\cup_{ r \in \N } Q_r))
\subset \cup_{ r \in \N }( \overline Q_r \setminus Q_r)), $
вытекает, что почти для всех $ x \in \R $ имеет место равенство $
\chi_W (x) = \sum_{r =1}^\infty \chi_{Q_r} (x). $ Поэтому почти
для всех $ x \in \R $ получаем
\begin{multline*}
h(x) = f(x) -g(x) = f(x) \chi_F (x) +f(x) \chi_W (x) -g(x) \\
=f(x) \chi_F (x) +f(x) (\sum_{r =1}^\infty \chi_{Q_r} (x)) -g(x)
\\ = f(x) (\sum_{r =1}^\infty \chi_{Q_r} (x)) -\sum_{r=1}^\infty (1/
\mes Q_r) (\int_{Q_r} f(y) dy) \chi_{Q_r}(x)\\
 = \sum_{r=1}^\infty
( f(x) -(1/ \mes Q_r) \int_{Q_r} f(y) dy) \chi_{Q_r}(x) =
\sum_{r=1}^\infty h_r (x), \end{multline*}
где $ h_r (x) = ( f(x)
-(1/ \mes Q_r) \int_{Q_r} f(y) dy) \chi_{Q_r}(x). $

 Учитывая (2.1.4), (2.1.5), на основании (2.1.2) и (2.1.7) заключаем,
что почти для всех $ x \in \R $ выполняется неравенство
$ | g(x) | \le c_2 \alpha, $ из которого вытекает оценка
\begin{multline*} \tag{2.2.5}
\| g \|_{L_2(\R)}^2 = \int_\R | g(x) |^2 dx \le \int_\R c_2 \alpha
| g(x) | dx = \\
c_2 \alpha \int_\R | f(x) \chi_F(x) +\sum_{r=1}^\infty
(1/ \mes Q_r) (\int_{Q_r} f(y) dy) \chi_{Q_r}(x) | dx \le \\
c_2 \alpha \int_\R | f(x) | \chi_F(x) +\sum_{r=1}^\infty
(1/ \mes Q_r) | \int_{Q_r} f(y) dy | \chi_{Q_r}(x)  dx = \\
c_2 \alpha (\int_\R | f(x) | \chi_F(x) dx +\sum_{r=1}^\infty \int_\R
(1/ \mes Q_r) | \int_{Q_r} f(y) dy | \chi_{Q_r}(x)  dx ) = \\
c_2 \alpha (\int_F | f(x) | dx +\sum_{r=1}^\infty
(1/ \mes Q_r) | \int_{Q_r} f(y) dy | \int_\R \chi_{Q_r}(x) dx ) \le \\
c_2 \alpha (\int_F | f(x) | dx +\sum_{r=1}^\infty
(1/ \mes Q_r) (\int_{Q_r} | f(y) | dy) \mes Q_r ) = \\
c_2 \alpha (\int_F | f(x) | dx +\sum_{r=1}^\infty
\int_{Q_r} | f(x) | dx ) = \\
c_2 \alpha \int_{F \cup (\cup_{r=1}^\infty Q_r)} | f(x) | dx  =
c_2 \alpha \int_\R | f(x) | dx = c_2 \alpha \| f \|_{L_1(\R)}.
\end{multline*}

Для получения (2.2.4), фиксируя множество $ A \subset \R:
\mes A =0 $ и для $ x \in \R \setminus A $ ввиду (2.2.2)
имеет место неравенство
$$
| (T f)(x) | = | (T (g +h))(x) | \le | (T g)(x) | +| (T h)(x) |,
$$
видим, что
\begin{multline*}
(\{ x \in \R: | (T f)(x) | > \alpha \} \setminus A) \subset \{ x
\in \R: | (T g)(x)| +| (T h)(x) | > \alpha \} \\
\subset \{ x \in \R: | (T g)(x) | > \alpha /2 \} \cup \{ x \in
\R: | (T h)(x) |
> \alpha /2 \},
\end{multline*} и, значит,
\begin{multline*} \tag{2.2.6}
\mes \{ x \in \R: | (T f)(x) | > \alpha \} = \mes (\{ x \in \R:
| (T f)(x) | > \alpha \} \setminus A ) \\
\le \mes \{ x \in \R: | (T g)(x) | > \alpha /2 \} + \mes \{ x
\in \R: | (T h)(x) | > \alpha /2 \}. \end{multline*}

Из (2.2.3) и (2.2.5) выводим
\begin{multline*} \tag{2.2.7}
\mes \{ x \in \R: | (T g)(x) | > \alpha /2 \} \le ((2 / \alpha)
\| g \|_{L_2(\R)})^2 = (2 / \alpha)^2 \| g \|_{L_2(\R)}^2 \\
\le (2 / \alpha)^2 c_2 \alpha \| f \|_{L_1(\R)} = (c_3 / \alpha)
\| f \|_{L_1(\R)}.
\end{multline*}

Для оценки второго слагаемого в правой части (2.2.6) имеем
\begin{multline*}
\{ x \in \R: | (T h)(x) | > \alpha /2 \} \\
= \{ x \in F: | (T
h)(x) | > \alpha /2 \} \cup \{ x \in W: | (T h)(x) | > \alpha /2
\},
\end{multline*}
а, следовательно,
\begin{multline*} \tag{2.2.8}
\mes \{ x \in \R: | (T h)(x) | > \alpha /2 \} \\
= \mes \{ x \in
F: | (T h)(x) | > \alpha /2 \} + \mes \{ x \in W: | (T h)(x) |
> \alpha /2 \}.
\end{multline*}

Второе слагаемое в правой части (2.2.8) в силу (2.1.3) удовлетворяет
неравенству
\begin{equation*} \tag{2.2.9}
\mes \{ x \in W: | (T h)(x) | > \alpha /2 \} \le
\mes W \le (c_4 / \alpha) \int_\R |f| dx = (c_4 / \alpha) \| f \|_{L_1(\R)}.
\end{equation*}

С другой стороны
\begin{equation*} \tag{2.2.10}
\mes \{ x \in F: | (T h)(x) | > \alpha /2 \} \\
\le (2 / \alpha) \| T h \|_{L_1(F)}.
\end{equation*}

Для проведения оценки правой части (2.2.10) определим при $ m \in \N $
функцию $ h_m^\prime $ равенством
$$
h_m^\prime = h -\sum_{r=1}^m h_r
$$
и заметим, что вследствие (2.2.2) при $ m \in \N $ почти для всех
$ x \in F$ справедливо неравенство
$$
| (T h)(x) | = | (T (\sum_{r=1}^m h_r +h_m^\prime))(x) | \le
\sum_{r=1}^m | (T h_r)(x) | +| (T h_m^\prime)(x) |,
$$
которое влечет оценку
\begin{multline*} \tag{2.2.11}
\| T h \|_{L_1(F)} = \int_{F} | (T h)(x) | dx
\le \int_{F} (\sum_{r=1}^m | (T h_r)(x) |) +| (T
h_m^\prime)(x) | dx\\
= \sum_{r=1}^m \int_{F} | (T h_r)(x) | dx +\int_{F} | (T h_m^\prime)(x) | dx.
\end{multline*}

 При $ r \in \N $ оценим сверху значения $ | (T h_r)(x) | $
для $ x \in F. $

Отметим, что
$$
\R = ( \cup_{ n \in \Z} Q_{k, n})  \cup \{ 2^{-k} \nu:  \nu \in \Z \},
$$
где $ Q_{k, n} = 2^{-k} n +2^{-k} I, n \in \Z, $
и, значит,
\begin{multline*}
F = F \cap \R = ( \cup_{ n \in \Z} ( F \cap Q_{k, n}))
\cup ( F \cap \{ 2^{-k} \nu:  \nu \in \Z \}) \\
= (\cup_{ n \in \Z: F \cap Q_{k, n} \ne
\emptyset } ( F \cap Q_{k, n}))  \cup ( F \cap \{ 2^{-k} \nu:
\nu \in \Z \}).
\end{multline*}

Учитывая это замечание, проведем оценку сверху $ | (T h_r)(x) | $
при $ r \in \N, n \in \Z:
F \cap Q_{k, n} \ne \emptyset $ для $ x \in F \cap Q_{k, n}. $
Принимая во внимание (1.2.19), при $ r \in \N, n \in \Z: F \cap Q_{k, n} \ne
\emptyset, $ почти для всех $ x \in F \cap
Q_{k, n} $ имеем
\begin{multline*}
| (T h_r)(x) | = \biggl| \sum_{\kappa =0}^k \sigma_\kappa
(\mathcal E_\kappa ( H_r))(x)\biggr|  =
\biggl| \sum_{\kappa =0}^k \sigma_\kappa
(( E_\kappa -E_{\kappa -1})( H_r))(x)\biggr| \\ \le
\sum_{\kappa =0}^k \biggl| \sigma_\kappa
(( E_\kappa(h_r))(x) -(E_{\kappa -1}( H_r))(x)) \biggr|  \le
\sum_{\kappa =0}^k ( | ( E_\kappa(h_r))(x) | +
| (E_{\kappa -1}( H_r))(x) |) \\ \le
2 \sum_{\kappa =0}^k | ( E_\kappa(h_r))(x) | =
2 \sum_{\kappa =0}^{k} \biggl| \sum_{ \rho_\kappa \in \Z} 2^\kappa
\times\biggl(\int_{ \supp \phi_{\kappa, \rho_\kappa}^* } h_r(y)
\overline { \phi_{\kappa, \rho_\kappa}^*(y) } dy \biggr)
\phi_{\kappa, \rho_\kappa}(x) \biggr| \\
= 2 \sum_{\kappa =0, \ldots, k: \exists \rho_\kappa \in \Z:
\supp \phi_{\kappa, \rho_\kappa} \cap Q_{k, n} \ne \emptyset}
\biggl|  \sum_{ \rho_\kappa \in \Z: \supp \phi_{\kappa, \rho_\kappa} \cap
Q_{k, n} \ne \emptyset} 2^\kappa \\
\times\biggl(\int_{ \supp \phi_{\kappa , \rho_\kappa}^* } h_r(y)
\overline { \phi_{\kappa, \rho_\kappa}^*(y) } dy \biggr)
\phi_{\kappa, \rho_\kappa}(x) \biggr| \le \\
2 \sum_{\substack{\kappa =0, \ldots, k: \exists \rho_\kappa \in \Z:
\supp \phi_{\kappa, \rho_\kappa} \cap Q_{k, n} \ne \emptyset, \\
\supp \phi_{\kappa, \rho_\kappa}^* \cap Q_r \ne \emptyset}}
\sum_{\substack{ \rho_\kappa \in \Z: \supp \phi_{\kappa, \rho_\kappa} \cap
Q_{k, n} \ne \emptyset, \\ \supp \phi_{\kappa, \rho_\kappa}^* \cap
Q_r \ne \emptyset}} 2^\kappa \\
\times \biggl| \int_{ \supp \phi_{\kappa, \rho_\kappa}^* \cap Q_r} h_r(y)
\overline { \phi_{\kappa, \rho_\kappa}^*(y) } dy \biggr|
\times  | \phi_{\kappa, \rho_\kappa}(x)| =
\end{multline*}
\begin{multline*} \tag{2.2.12}
2 \sum_{\substack{\kappa =0, \ldots, k: \exists \rho_\kappa \in \Z:
\supp \phi_{\kappa, \rho_\kappa} \cap Q_{k, n} \ne \emptyset, \\
\supp \phi_{\kappa, \rho_\kappa}^* \cap Q_r \ne \emptyset}}
\sum_{\substack{ \rho_\kappa \in \Z: \supp \phi_{\kappa, \rho_\kappa} \cap
Q_{k, n} \ne \emptyset, \\ \supp \phi_{\kappa, \rho_\kappa}^* \cap
Q_r \ne \emptyset}} 2^\kappa \\
\times \biggl| \int_{ \supp \phi_{\kappa, \rho_\kappa}^* \cap Q_r} h_r(y)
\overline { \phi^*(2^\kappa y -\rho_\kappa) } dy \biggr|
\times | \phi(2^\kappa x -\rho_\kappa) | \le \\
2 \sum_{\substack{\kappa =0, \ldots, k: \exists \rho_\kappa \in \Z:
\supp \phi_{\kappa, \rho_\kappa} \cap Q_{k, n} \ne \emptyset, \\
\supp \phi_{\kappa, \rho_\kappa}^* \cap Q_r \ne \emptyset}}
\sum_{\substack{ \rho_\kappa \in \Z: \supp \phi_{\kappa, \rho_\kappa} \cap
Q_{k, n} \ne \emptyset,\\ \supp \phi_{\kappa, \rho_\kappa}^* \cap
Q_r \ne \emptyset}} 2^\kappa \\
\times \biggl| \int_{ \supp \phi_{\kappa, \rho_\kappa}^* \cap Q_r} h_r(y)
\overline { \phi^*(2^\kappa y -\rho_\kappa) } dy \biggr|
\times \| \phi \|_{L_\infty(\R)} = \\
c_5 \sum_{\substack{\kappa =0, \ldots, k: \exists \rho_\kappa \in \Z:
\supp \phi_{\kappa, \rho_\kappa} \cap Q_{k, n} \ne \emptyset, \\
\supp \phi_{\kappa, \rho_\kappa}^* \cap Q_r \ne \emptyset}}
\sum_{\substack{ \rho_\kappa \in \Z: \supp \phi_{\kappa, \rho_\kappa} \cap
Q_{k, n} \ne \emptyset,\\ \supp \phi_{\kappa, \rho_\kappa}^* \cap
Q_r \ne \emptyset}} 2^\kappa \\
\times \biggl| \int_{ \supp \phi_{\kappa, \rho_\kappa}^* \cap Q_r} h_r(y)
\overline { \phi^*(2^\kappa y -\rho_\kappa) } dy \biggr| = \\
c_5 \sum_{\substack{\kappa =0, \ldots, k: \exists \rho_\kappa \in \Z:
\supp \phi_{\kappa, \rho_\kappa} \cap Q_{k, n} \ne \emptyset, \\
\supp \phi_{\kappa, \rho_\kappa}^* \cap Q_r \ne \emptyset}}
\sum_{\substack{ \rho_\kappa \in \Z: \supp \phi_{\kappa, \rho_\kappa} \cap
Q_{k, n} \ne \emptyset,\\ \supp \phi_{\kappa, \rho_\kappa}^* \cap
Q_r \ne \emptyset}} 2^\kappa \\
\times \biggl| \int_{ Q_r} h_r(y)
\overline { \phi^*(2^\kappa y -\rho_\kappa) } dy \biggr| = \\
c_5 \sum_{\substack{\kappa =0, \ldots, \iota: \exists \rho_\kappa \in \Z:
\supp \phi_{\kappa, \rho_\kappa} \cap Q_{k, n} \ne \emptyset, \\
\supp \phi_{\kappa, \rho_\kappa}^* \cap Q_r \ne \emptyset}}
\sum_{\substack{ \rho_\kappa \in \Z: \supp \phi_{\kappa, \rho_\kappa} \cap
Q_{k, n} \ne \emptyset,\\ \supp \phi_{\kappa, \rho_\kappa}^* \cap
Q_r \ne \emptyset}} 2^\kappa \\
\times \biggl| \int_{ Q_r} h_r(y)
\overline { \phi^*(2^\kappa y -\rho_\kappa) } dy \biggr|,
\end{multline*}
где
$
\iota = \iota(r,n) = \max \{ \kappa =0, \ldots, k| \exists \rho_\kappa \in \Z:
\supp \phi_{\kappa, \rho_\kappa} \cap Q_{k, n} \ne \emptyset,
\supp \phi_{\kappa, \rho_\kappa}^* \cap Q_r \ne \emptyset \};
r \in \N, n \in \Z: F \cap Q_{k, n} \ne \emptyset.
$

Для оценки правой части (2.2.12) отметим, что при $ r \in \N $
выполняется равенство
\begin{multline*} \tag{2.2.13}
\int_{Q_r} h_r (x) dx = \int_{Q_r} ( f(x) -(1/ \mes Q_r)
\int_{Q_r} f(y) dy) \chi_{Q_r}(x) dx\\ = \int_{Q_r} ( f(x) -(1/
\mes Q_r) \int_{Q_r} f(y) dy) dx \\
= \int_{Q_r} f(x) dx -(1/ \mes Q_r) (\int_{Q_r} f(y) dy)
\int_{Q_r} dx \\
= \int_{Q_r} f(x) dx -\int_{Q_r} f(y) dy =0.
\end{multline*}

Фиксируя для каждого $ r \in \N $ точку $ y_r \in Q_r, $
при $ r \in \N, n \in \Z: F \cap
Q_{k, n} \ne \emptyset, \kappa =0, \ldots, \iota: \exists \rho_\kappa
\in \Z: \supp \phi_{\kappa, \rho_\kappa} \cap Q_{k, n} \ne \emptyset,
\supp \phi_{\kappa, \rho_\kappa}^* \cap Q_r \ne \emptyset; \rho_\kappa
\in \Z: \supp \phi_{\kappa, \rho_\kappa} \cap Q_{k, n} \ne \emptyset,
\supp \phi_{\kappa, \rho_\kappa}^* \cap Q_r \ne \emptyset $
с учетом (2.2.13) выводим
\begin{multline*} \tag{2.2.14}
| \int_{ Q_r} h_r(y) \overline { \phi^*( 2^\kappa y -\rho_\kappa) } dy | = \\
| \int_{ Q_r} h_r(y) (\overline { \phi^*( 2^\kappa y -\rho_\kappa) } -
\overline { \phi^*( 2^\kappa y_r -\rho_\kappa) }) +h_r(y)
\overline { \phi^*( 2^\kappa y_r -\rho_\kappa) } dy | = \\
| \int_{ Q_r} h_r(y) (\overline { \phi^*( 2^\kappa y -\rho_\kappa) } -
\overline { \phi^*( 2^\kappa y_r -\rho_\kappa) }) dy +
\overline { \phi^*( 2^\kappa y_r -\rho_\kappa) } \int_{ Q_r} h_r(y) dy | = \\
| \int_{ Q_r} h_r(y) (\overline { \phi^*( 2^\kappa y -\rho_\kappa) } -
\overline { \phi^*( 2^\kappa y_r -\rho_\kappa) }) dy | \le \\
\int_{ Q_r} | h_r(y)| \cdot | \overline { \phi^*( 2^\kappa y -\rho_\kappa) } -
\overline { \phi^*( 2^\kappa y_r -\rho_\kappa) } | dy.
\end{multline*}

Оценивая правую часть (2.2.14), при $ r \in \N, n \in \Z:
F \cap Q_{k, n} \ne \emptyset, \kappa =0, \ldots, \iota: \exists \rho_\kappa
\in \Z: \supp \phi_{\kappa, \rho_\kappa} \cap Q_{k, n} \ne \emptyset,
\supp \phi_{\kappa, \rho_\kappa}^* \cap Q_r \ne \emptyset; \rho_\kappa
\in \Z: \supp \phi_{\kappa, \rho_\kappa} \cap Q_{k, n} \ne \emptyset,
\supp \phi_{\kappa, \rho_\kappa}^* \cap Q_r \ne \emptyset, $
для $ y \in Q_r $ имеем
\begin{multline*} \tag{2.2.15}
| \overline { \phi^*( 2^\kappa y -\rho_\kappa) } -
\overline { \phi^*( 2^\kappa y_r -\rho_\kappa) } | =
| \phi^*( 2^\kappa y -\rho_\kappa) -\phi^*( 2^\kappa y_r -\rho_\kappa) | = \\
| \int_0^1 \frac{d} {dt} ( \phi^*( 2^\kappa (y_r +t (y -y_r)) -
\rho_\kappa) ) dt | = \\
| \int_0^1 2^\kappa (y -y_r) \frac{d \phi^*} {du} ( 2^\kappa (y_r +t (y -y_r))
-\rho_\kappa) dt | = \\
2^\kappa | y -y_r | \cdot | \int_0^1 \frac{d \phi^*} {du} ( 2^\kappa (y_r +
t (y -y_r)) -\rho_\kappa) dt | \le \\
2^\kappa | y -y_r | \cdot \int_0^1 | \frac{d \phi^*} {du} ( 2^\kappa (y_r +
t (y -y_r)) -\rho_\kappa) | dt \le \\
2^\kappa | y -y_r | \cdot \| \frac{d \phi^*} {du} \|_{L_\infty(\R)} \le
2^\kappa ( \diam Q_r ) \| \frac{d \phi^*} {du} \|_{L_\infty(\R)} =
c_6(\phi^*) 2^\kappa ( \diam Q_r ).
\end{multline*}

Соединяя (2.2.14) и (2.2.15), находим, что при $ r \in \N, n \in \Z:
F \cap Q_{k, n} \ne \emptyset, \kappa =0, \ldots, \iota: \exists \rho_\kappa
\in \Z: \supp \phi_{\kappa, \rho_\kappa} \cap Q_{k, n} \ne \emptyset,
\supp \phi_{\kappa, \rho_\kappa}^* \cap Q_r \ne \emptyset; \rho_\kappa
\in \Z: \supp \phi_{\kappa, \rho_\kappa} \cap Q_{k, n} \ne \emptyset,
\supp \phi_{\kappa, \rho_\kappa}^* \cap Q_r \ne \emptyset, $ выполняется
неравенство
\begin{multline*}
| \int_{ Q_r} h_r(y) \overline { \phi^*( 2^\kappa y -\rho_\kappa) } dy | \le \\
\int_{ Q_r} c_6 2^\kappa ( \diam Q_r ) | h_r(y) | dy = c_6
2^\kappa ( \diam Q_r ) \int_{ Q_r} | h_r(y) | dy.
\end{multline*}

Подставляя последнюю оценку в (2.2.12), получаем, что
при $ r \in \N, n \in \Z: F \cap Q_{k, n} \ne \emptyset, $ почти для всех
$ x \in F \cap Q_{k, n} $ соблюдается неравенство
\begin{multline*} \tag{2.2.16}
| (T h_r)(x) | \le
c_5 \sum_{\substack{\kappa =0, \ldots, \iota: \exists \rho_\kappa \in \Z:
\supp \phi_{\kappa, \rho_\kappa} \cap Q_{k, n} \ne \emptyset,\\
\supp \phi_{\kappa, \rho_\kappa}^* \cap Q_r \ne \emptyset}}
\sum_{\substack{ \rho_\kappa \in \Z: \supp \phi_{\kappa, \rho_\kappa} \cap
Q_{k, n} \ne \emptyset,\\ \supp \phi_{\kappa, \rho_\kappa}^* \cap
Q_r \ne \emptyset}} 2^\kappa \\
c_6 2^\kappa ( \diam Q_r )
\int_{ Q_r} | h_r(y) | dy =
c_7 ( \diam Q_r ) (\int_{ Q_r} | h_r(y) | dy) \cdot \\
\sum_{\substack{\kappa =0, \ldots, \iota: \exists \rho_\kappa \in \Z:
\supp \phi_{\kappa, \rho_\kappa} \cap Q_{k, n} \ne \emptyset, \\
\supp \phi_{\kappa, \rho_\kappa}^* \cap Q_r \ne \emptyset}}
\card \{ \rho_\kappa \in \Z: \supp \phi_{\kappa, \rho_\kappa} \cap
Q_{k, n} \ne \emptyset, \\
\supp \phi_{\kappa, \rho_\kappa}^* \cap
Q_r \ne \emptyset \} 2^{2 \kappa} \le
c_7 ( \diam Q_r ) (\int_{ Q_r} | h_r(y) | dy) \cdot \\
\sum_{\substack{\kappa =0, \ldots, \iota: \exists \rho_\kappa \in \Z:
\supp \phi_{\kappa, \rho_\kappa} \cap Q_{k, n} \ne \emptyset, \\
\supp \phi_{\kappa, \rho_\kappa}^* \cap Q_r \ne \emptyset}}
\card \{ \rho_\kappa \in \Z: \supp \phi_{\kappa, \rho_\kappa} \cap
Q_{k, n} \ne \emptyset \} 2^{2 \kappa} \le \\
c_7 ( \diam Q_r ) (\int_{ Q_r} | h_r(y) | dy)
\sum_{\kappa =0}^\iota
\card \{ \rho_\kappa \in \Z: \supp \phi_{\kappa, \rho_\kappa} \cap
Q_{k, n} \ne \emptyset \} 2^{2 \kappa}.
\end{multline*}

При $ \kappa =0, \ldots, \iota, $ для $  \rho_\kappa \in \Z:
\supp \phi_{\kappa, \rho_\kappa} \cap Q_{k, n} \ne \emptyset, $
выбирая
$ \xi \in I $ и $ \eta \in \supp \phi, $ для которых
$ 2^{-k} n +2^{-k} \xi = 2^{-\kappa} \rho_\kappa +2^{-\kappa} \eta, $
находим, что
$ \rho_\kappa = 2^{\kappa -k} n +2^{\kappa -k} \xi -\eta, $
откуда
$$
2^{\kappa -k} n  -\eta \le \rho_\kappa \le 2^{\kappa -k} n +1 -\eta,
$$
или
$$
2^{\kappa -k} n  -\sup_{u \in \supp \phi} u \le \rho_\kappa \le
2^{\kappa -k} n +1 -\inf_{u \in \supp \phi} u,
$$
а, следоаательно,
$$
\card \{ \rho_\kappa \in \Z: \supp \phi_{\kappa, \rho_\kappa} \cap
Q_{k, n} \ne \emptyset \} \le c_8(\phi).
$$

Применяя эту оценку в (2.2.16), находим, что
при $ r \in \N, n \in \Z: F \cap Q_{k, n} \ne \emptyset, $ почти для всех
$ x \in F \cap Q_{k, n} $ справедливо неравенство
\begin{multline*} \tag{2.2.17}
| (T h_r)(x) | \le
c_7 ( \diam Q_r ) (\int_{ Q_r} | h_r(y) | dy)
\sum_{\kappa =0}^\iota c_8 2^{2 \kappa} \le \\
c_9 ( \diam Q_r ) (\int_{ Q_r} | h_r(y) | dy) 2^{2 \iota}
\sum_{\kappa =0}^\iota 2^{-2(\iota -\kappa)} \le \\
c_9 ( \diam Q_r ) (\int_{ Q_r} | h_r(y) | dy) 2^{2 \iota}
\sum_{s=0}^\infty 2^{-2s} =  c_{10} ( \diam Q_r ) 2^{2 \iota} \int_{
Q_r} | h_r(y) | dy.
\end{multline*}

Заметим, что при $ r \in \N $ в силу (2.1.7) верно неравенство
\begin{multline*} \tag{2.2.18}
\int_{Q_r} | h_r (y) | dy = \int_{Q_r} | ( f(y) -(1/ \mes Q_r)
\int_{Q_r} f(z) dz) \chi_{Q_r}(y) | dy = \\
\int_{Q_r} | f(y) -(1/ \mes Q_r)
\int_{Q_r} f(z) dz | dy \le \\
\int_{Q_r} | f(y) | dy +(1/ \mes Q_r)
| \int_{Q_r} f(z) dz | \int_{Q_r} dy
\le \\
\int_{Q_r} | f(y) | dy +\int_{Q_r} | f(z) | dz = 2 \int_{Q_r} | f(y) | dy
\le c_{11} \alpha \mes Q_r.
\end{multline*}

Кроме того, благодаря (2.1.6), при $ r \in \N, $ для $ y \in Q_r $
справедливо неравенство
\begin{equation*} \tag{2.2.19}
\diam Q_r \le c_{12} \rho(y, F).
\end{equation*}

А также при $ r \in \N,
n \in \Z: F \cap Q_{k, n} \ne \emptyset, $ для
$ x \in F \cap Q_{k, n} $ и $ y \in Q_r, $ выбирая $ \rho_\iota \in \Z:
\supp \phi_{\iota, \rho_\iota} \cap Q_{k, n} \ne \emptyset,
\supp \phi_{\iota, \rho_\iota}^* \cap Q_r \ne \emptyset $ и
$ \xi \in \supp \phi_{\iota, \rho_\iota} \cap Q_{k, n},
\eta \in \supp \phi_{\iota, \rho_\iota}^* \cap Q_r, $
имеем
\begin{multline*}
| x -y| \le | x -\xi| +| \xi -\eta| +| \eta -y| \le \\
\diam Q_{k, n} +
\sup_{u \in \supp \phi_{\iota, \rho_\iota}, v \in
\supp \phi_{\iota, \rho_\iota}^* } | u -v | +\diam Q_r,
\end{multline*}
причем,
\begin{multline*}
\diam Q_{k, n} = 2^{-k} \le 2^{-\iota}, \\
\sup_{u \in \supp \phi_{\iota, \rho_\iota},
v \in \supp \phi_{\iota, \rho_\iota}^* } | u -v | =
\sup_{u \in (2^{-\iota} \rho_\iota +2^{-\iota} \supp \phi),
v \in (2^{-\iota} \rho_\iota +2^{-\iota} \supp \phi^*)} | u -v | = \\
\sup_{u \in \supp \phi, v \in \supp \phi^*}
| (2^{-\iota} \rho_\iota +2^{-\iota} u)
-(2^{-\iota} \rho_\iota +2^{-\iota} v) | = \\
\sup_{u \in \supp \phi, v \in \supp \phi^*} 2^{-\iota} | u -v | =
2^{-\iota} \sup_{u \in \supp \phi, v \in \supp \phi^*} | u -v | =
c_{13} 2^{-\iota}, \\
\text{и с учетом (2.1.6)} \\
\diam Q_r \le c_{14} \inf_{u \in F, v \in Q_r} | u -v| \le
c_{14} | x -\eta| \le c_{14} (| x -\xi| +| \xi -\eta|) \\
\le c_{14} (\diam Q_{k, n} +\sup_{u \in \supp \phi_{\iota, \rho_\iota},
v \in \supp \phi_{\iota, \rho_\iota}^* } | u -v | ),
\end{multline*}
а, следовательно, соблюдается неравенство
\begin{equation*} \tag{2.2.20}
2^\iota \le c_{15} | x -y |^{-1}.
\end{equation*}

Подставляя (2.2.18) -- (2.2.20) в (2.2.17), заключаем, что
при $ r \in \N, n \in \Z: F \cap Q_{k, n} \ne \emptyset, $ почти для всех
$ x \in F \cap Q_{k, n} $ имеет место неравенство
\begin{multline*} \tag{2.2.21}
| (T h_r)(x) | \le c_{16} ( \diam Q_r ) 2^{2 \iota} \alpha \mes
Q_r = c_{16} \alpha ( \diam Q_r ) 2^{2 \iota} \int_{ Q_r } dy \\
= c_{16} \alpha \int_{ Q_r } ( \diam Q_r ) 2^{2 \iota} dy \le
c_{17} \alpha \int_{ Q_r } \rho(y, F) | x -y |^{-2} dy.
\end{multline*}

Из сказанного ясно, что при $ r \in \N $
неравенство (2.2.21) выполняется почти для всех $ x \in F, $
а, следовательно,
\begin{multline*} \tag{2.2.22}
\int_{ F } | (T h_r)(x) | dx \le \int_{ F } c_{17}
\alpha (\int_{ Q_r } \rho(y, F) | x -y |^{-2} dy ) dx\\
= c_{17} \alpha \int_F \int_{ Q_r } \rho(y, F) | x -y |^{-2} dy
dx.
\end{multline*}

Подставляя (2.2.22) в (2.2.11) и с учетом (2.1.4), (2.1.5)
применяя (2.1.1), а затем (2.1.3), приходим к выводу, что при $ m
\in \N $ справедливо неравенство
\begin{multline*} \tag{2.2.23}
\| T h \|_{L_1(F)} \le
c_{17} \alpha \sum_{r=1}^m \int_F \int_{ Q_r } \rho(y, F) | x -y |^{-2} dy dx
+ \\ +\int_{F} | (T h_m^\prime)(x) | dx \le
c_{17} \alpha \int_F (\sum_{r=1}^m \int_{ Q_r } \rho(y, F) | x -y |^{-2} dy) dx
 + \\ +\int_{\R} | (T h_m^\prime)(x) | dx =
c_{17} \alpha \int_F \int_{ \cup_{r=1}^m Q_r } \rho(y, F) | x -y |^{-2} dy dx
+ \\ +\int_{\R} | (T h_m^\prime)(x) | dx \le
c_{17} \alpha \int_F \int_W \rho(y, F) | x -y |^{-2} dy dx + \\
+\int_{\R} | (T h_m^\prime)(x) | dx \le
c_{18} \alpha \mes (\R \setminus F)
+\int_{\R} | (T h_m^\prime)(x) | dx = \\
c_{18} \alpha \mes W +\int_{\R} | (T h_m^\prime)(x) | dx
\le c_{19} \int_{\R} |f(x)| dx
+\int_{\R} | (T h_m^\prime)(x) | dx = \\
c_{19}  \| f \|_{L_1(\R)}
+\int_\R | (T h_m^\prime)(x) | dx.
\end{multline*}

Принимая во внимание, что при $ m \in \N $ в силу
оценок (1.2.20) и (2.2.18), условия (2.1.4) и включения
$ f \in L_1(\R) $ имеет место соотношение
\begin{multline*}
\int_\R | (T h_m^\prime)(x) | dx \\
= \int_\R | \sum_{\kappa =0}^k \sigma_\kappa \cdot (\mathcal E_\kappa
( h_m^\prime ))(x) | dx \le \int_\R \sum_{\kappa =0}^k |
(\mathcal E_\kappa( h_m^\prime ))(x) | dx \\=
\sum_{\kappa =0}^k \int_\R | (\mathcal E_\kappa( h_m^\prime
))(x) | dx = \sum_{\kappa =0}^k \| \mathcal E_\kappa(h_m^\prime)
\|_{L_1(\R)} \\ \le
\sum_{\kappa =0}^k (\| E_\kappa(h_m^\prime) \|_{L_1(\R)} +
\| E_{\kappa -1}(h_m^\prime) \|_{L_1(\R)}) \\
\le \sum_{\kappa =0}^k c_{20} \| h_m^\prime \|_{L_1(\R)} = c_{20} (k +1) \| h_m^\prime \|_{L_1(\R)} = \\
c_{20} (k +1) \| \sum_{r=m+1}^\infty h_r \|_{L_1(\R)} \le
c_{20} (k +1) \sum_{r=m+1}^\infty \| h_r \|_{L_1(\R)} = \\
c_{20} (k +1) \sum_{r=m+1}^\infty \int_\R | h_r (x) | dx =
c_{20} (k +1) \sum_{r=m+1}^\infty \int_{ Q_r} | h_r (x) | dx \le \\
c_{20} (k +1) \sum_{r=m+1}^\infty 2 \int_{ Q_r} | f(x) | dx \to 0
\text{ при } m \to \infty,
\end{multline*}
из (2.2.23) следует, что выполняется неравенство
\begin{equation*} \tag{2.2.24}
\| T h \|_{L_1(F)} \le c_{19}  \| f \|_{L_1(\R)}.
\end{equation*}

Подставляя (2.2.24) в (2.2.10), выводим неравенство
$$
\mes \{ x \in F: | (T h)(x) | > \alpha /2 \} \le
(c_{21} / \alpha) \| f \|_{L_1(\R)},
$$
которое в соединении с (2.2.8), (2.2.9) дает оценку
\begin{equation*} \tag{2.2.25}
\mes \{ x \in \R: | (T h)(x) | > \alpha /2 \} \le
(c_{22} / \alpha) \| f \|_{L_1(\R)}.
\end{equation*}

Объединяя (2.2.6), (2.2.7) и (2.2.25), приходим к (2.2.4).
Сопоставляя (2.2.2) -- (2.2.4) с (2.1.8) -- (2.1.10), заключаем, что
при $ k \in \N, \sigma = \{ \sigma_\kappa \in \{-1, 1\}:
\kappa =0, \ldots, k \}, 1 < p < 2 $ для вещественнозначной $ f \in L_p(\R) $
согласно (2.1.11) верно неравенство
$$
\| T f \|_{L_p(\R)} \le c_1 \| f \|_{L_p(\R)},
$$
из которого следует (2.2.1) при $ 1 < p < 2 $ для комплекснозначных
$ f \in L_p(\R). $
Справедливость (2.2.1) при $ p =2 $ установлена при выводе (2.2.3).

Теперь проверим соблюдение (2.2.1)
при $ 2 < p < \infty. $
В самом деле, при $ k \in \N, \sigma = \{ \sigma_\kappa \in \{-1, 1\}:
\kappa =0, \ldots, k \}, 2 < p < \infty, p^\prime = p/(p-1) \in (1,2), $
в силу (1.2.30) для $ f \in L_p(\R), g \in L_{p^\prime}(\R) $
справедливо равенство
$$
\int_{\R} (\sum_{\kappa =0}^k \sigma_\kappa \cdot \mathcal E_\kappa f)
\cdot \overline g dx =
\int_{\R} f \cdot \overline {(\sum_{\kappa =0}^k \sigma_\kappa \cdot
\mathcal E_\kappa g )} dx,
$$
из которого ввиду того, что $ L_p(\R) = (L_{p^\prime}(\R))^*, $
и отображение $ S: L_{p^\prime}(\R) \ni g \mapsto S g = \overline g \in
L_{p^\prime}(\R) $ является изометрией $ L_{p^\prime}(\R) $ на себя,
следует, что оператор
$$
(\sum_{\kappa =0}^k \sigma_\kappa \cdot \mathcal E_\kappa^{p}) =
(\sum_{\kappa =0}^k \sigma_\kappa \cdot \mathcal E_\kappa) \mid_{L_p(\R)}
$$
является сопряженным к оператору
$$
S (\sum_{\kappa =0}^k \sigma_\kappa \cdot \mathcal E_\kappa^{p^\prime}) S
$$
а, значит,
\begin{multline*}
\| (\sum_{\kappa =0}^k \sigma_\kappa \cdot \mathcal E_\kappa^{p})
\|_{\mathcal B(L_p(\R), L_p(\R))} =
\| S (\sum_{\kappa =0}^k \sigma_\kappa \cdot \mathcal
E_\kappa^{p^\prime}) S \|_{ \mathcal
B(L_{p^\prime}(\R), L_{p^\prime}(\R))} = \\
\| (\sum_{\kappa =0}^k \sigma_\kappa \cdot \mathcal
E_\kappa^{p^\prime}) \|_{ \mathcal
B(L_{p^\prime}(\R), L_{p^\prime}(\R))}.
\end{multline*}
Отсюда, учитывая выполнение (2.2.1) при $ p^\prime $ вместо $ p, $
заключаем, что (2.2.1) соблюдается при $ 2 < p < \infty. \square $

При помощи леммы 2.2.1 устанавливается

Теорема 2.2.2

Пусть $ d \in \N, 1 < p < \infty $ и для $ j = 1,\ldots,d $ функции
$ \phi_j, \phi_j^* $ удовлетворяют условиям леммы 2.2.1, а функции
$ \phi, \phi^* $ задаются равенствами
$$
\phi(x) = \prod_{j =1}^d \phi_j(x_j),
\phi^*(x) = \prod_{j =1}^d \phi_j^*(x_j), x \in \R^d.
$$
Тогда существует константа$ c_{23}(d, \phi, \phi^*,p) >0 $ такая, что для
любого семейства чисел $ \{ \sigma_\kappa: \kappa \in \Z_+^d \} $ вида
$ \sigma_\kappa = \prod_{j=1}^d
\sigma^j_{\kappa_j}, $
где $ \sigma^j_{\kappa_j} \in \{-1, 1\}, \kappa_j \in \Z_+, j =1,\ldots,d, $
для $ f \in L_p(\R^d) $ справедливо неравенство
\begin{equation*} \tag{2.2.26}
\| \sum_{\kappa \in \Z_+^d} \sigma_\kappa \cdot
(\mathcal E_\kappa f ) \|_{L_p(\R^d)}
\le c_{23} \| f \|_{L_p(\R^d)},
\end{equation*}
где $ \mathcal E_\kappa = \mathcal E_\kappa^{\phi, \phi^*}, \kappa \in \Z_+^d. $

Доказательство.

Сначала покажем, что в условиях теоремы для любого непустого
множества $ J \subset \{1, \ldots, d \} $ при любом $ k^J \in (\Z_+^d)^J $
и любых наборах чисел $ \{ \sigma^j_{\kappa_j} \in \{-1, 1\}, \kappa_j
= 0,\ldots, k_j \}, j \in J, $
для $ f \in L_p(\R^d) $ имеет место неравенство
\begin{multline*} \tag{2.2.27}
\| \sum_{\kappa^J \in \Z_+^m(k^J)}  (\prod_{j \in J} (\sigma^j_{\kappa_j}
V_j (\mathcal E_{\kappa_j}^{j, p})))f \|_{L_p(\R^d)}
\le (\prod_{j \in J} c_{1}(\phi_j, \phi^*_j,p)) \cdot
\| f \|_{L_p(\R^d)},
\end{multline*}
где $ m = \card J, $ а $ \mathcal E_{\kappa_j}^{j, p} =
\mathcal E_{\kappa_j}^{\phi_j, \phi_j^*} \mid_{L_p(\R)}, \kappa_j \in \Z_+,
j =1,\ldots,d. $

Доказательство (2.2.27) проведем по индукции относительно $ m. $
При $ m =1, $ т.е. для $ j =1, \ldots, d, $ используя п. 2) леммы
1.3.2, теорему Фубини, (1.3.1), (2.2.1), имеем
\begin{multline*}
\| \sum_{\kappa_j =0}^{k_j} \sigma^j_{\kappa_j} \cdot (V_j
(\mathcal E_{\kappa_j}^{j, p})) f \|_{L_p(\R^d)}^p =
\| (V_j(\sum_{\kappa_j =0}^{k_j} \sigma^j_{\kappa_j} \mathcal
E_{\kappa_j}^{j, p})) f \|_{L_p(\R^d)}^p \\
= \int_{\R^d} | (V_j (\sum_{\kappa_j =0}^{k_j} \sigma^j_{\kappa_j}
\mathcal E_{\kappa_j}^{j, p})) f |^p dx\\ = \int_{\R^{d -1}}
\int_\R | ((V_j (\sum_{\kappa_j =0}^{k_j} \sigma^j_{\kappa_j}
\mathcal E_{\kappa_j}^{j, p})) f) (x_1,\ldots, x_{j -1}, x_j,
x_{j +1}, \ldots, x_d)|^p\\
\times dx_j dx_1 \ldots dx_{j -1} dx_{j +1} \ldots dx_d
\\
= \int_{\R^{d -1}} \int_\R | ( (\sum_{\kappa_j =0}^{k_j}
\sigma^j_{\kappa_j} \mathcal E_{\kappa_j}^{j, p}) f (x_1,\ldots,
x_{j -1}, \cdot, x_{j +1}, \ldots, x_d))(x_j) |^p\\
\times dx_j dx_1 \ldots
dx_{j -1} dx_{j +1} \ldots dx_d\\
 = \int_{\R^{d -1}} \|
(\sum_{\kappa_j =0}^{k_j} \sigma^j_{\kappa_j} \mathcal
E_{\kappa_j}^{j, p}) f (x_1,\ldots, x_{j -1}, \cdot, x_{j +1},
\ldots, x_d) \|_{L_p(\R)}^p \\
\times dx_1 \ldots dx_{j -1} dx_{j +1} \ldots
dx_d \\
\le \int_{\R^{d -1}} ( c_{1}(\phi_j, \phi_j^*,p) \| f (x_1,\ldots, x_{j -1},
\cdot, x_{j +1}, \ldots, x_d) \|_{L_p(\R)})^p\\
\times dx_1 \ldots dx_{j -1} dx_{j +1} \ldots dx_d\\ =
(c_{1}(\phi_j, \phi_j^*, p))^p \int_{\R^{d -1}} \int_\R | f (x_1,\ldots, x_{j -1},
x_j, x_{j +1}, \ldots, x_d) |^p\\
\times dx_j
dx_1 \ldots dx_{j -1} dx_{j +1} \ldots dx_d \\
= ( c_{1}(\phi_j, \phi_j^*,p))^p
\int_{\R^d} | f(x)|^p dx = ( c_{1}(\phi_j, \phi_j^*,p) \| f \|_{L_p(\R^d)})^p,
\end{multline*}
откуда
\begin{equation*} \tag{2.2.28}
\| \sum_{\kappa_j =0}^{k_j} \sigma^j_{\kappa_j} (V_j (\mathcal
E_{\kappa_j}^{j, p})) f \|_{L_p(\R^d)} \le c_{1}(\phi_j, \phi_j^*, p)
\| f \|_{L_p(\R^d)},
\end{equation*}
что совпадает с (2.2.27) при $ m =1, J = \{j\}. $

Предположим теперь, что при некотором $ m: 1 \le m \le d -1, $
оценка (2.2.27) имеет место для любого множества $ J \subset \{1,
\ldots, d \}: \card J = m, $ при любом $ k^J \in (\Z_+^d)^J, $
любых наборах чисел $ \{ \sigma^j_{\kappa_j} \in \{-1, 1\},
\kappa_j = 0, \ldots, k_j \}, j \in J, $ и любой функции $ f \in
L_p(\R^d).$ Покажем, что тогда неравенство (2.2.27) справедливо при
$ m+1 $ вместо $ m $ для любого множества $ \mathcal J \subset
\{1, \ldots, d \} $ вместо $ J, $ у которого $ \card \mathcal J =
m +1, $ при любом $ k^{ \mathcal J} \in (\Z_+^d)^{ \mathcal J}, $
любых наборах чисел $ \{ \sigma^j_{\kappa_j} \in \{-1, 1\},
\kappa_j = 0, \ldots, k_j \}, j \in \mathcal J, $ и любой функции
$ f \in L_p(\R^d). $ Представляя множество $ \mathcal J \subset
\{1, \ldots, d \}: \card \mathcal J = m +1, $ в виде $ \mathcal J
= J \cup \{i\}, i \notin J, $ с учетом (1.3.3) в силу (2.2.28) и
предположения индукции получаем
\begin{multline*}
\| \sum_{\kappa^{\mathcal J} \in \Z_+^{m+1}(k^{\mathcal J})}
(\prod_{j \in \mathcal J} (\sigma^j_{\kappa_j} V_j (\mathcal
E_{\kappa_j}^{j, p})))f \|_{L_p(\R^d)} \\
= \| \sum_{(\kappa_i, \kappa^J): \kappa_i =0, \ldots, k_i,
\kappa^J \in \Z_+^m(k^J)} \sigma^i_{\kappa_i} (V_i (\mathcal
E_{\kappa_i}^{i, p})) ((\prod_{j \in J} (\sigma^j_{\kappa_j} V_j
(\mathcal E_{\kappa_j}^{j, p})) )f) \|_{L_p(\R^d)}\\ =
\| \sum_{\kappa_i =0}^{k_i} \sigma^i_{\kappa_i} (V_i (\mathcal
E_{\kappa_i}^{i, p})) (\sum_{\kappa^J \in \Z_+^m(k^J)}  (\prod_{j
\in J} (\sigma^j_{\kappa_j} V_j (\mathcal E_{\kappa_j}^{j, p}))
)f) \|_{L_p(\R^d)} \\
\le c_{1}(\phi_i, \phi_i^*, p) \| \sum_{\kappa^J \in \Z_+^m(k^J)}  (\prod_{j
\in J} (\sigma^j_{\kappa_j} V_j (\mathcal E_{\kappa_j}^{j, p})
))f \|_{L_p(\R^d)} \\ \le c_{1}(l\phi_i, \phi_i^*,p) (\prod_{j \in J}
c_{1}(\phi_j, \phi_j^*, p)) \cdot \| f \|_{L_p(\R^d)} = (\prod_{j \in \mathcal
J} c_{1}(\phi_j, \phi_j^*, p)) \cdot \| f \|_{L_p(\R^d)},
\end{multline*}
что завершает вывод (2.2.27).

В частности, из (2.2.27) при $ m = d $ ввиду (1.4.2) получаем, что
в условиях теоремы при любом $ k \in \Z_+^d $ соблюдается неравенство
\begin{multline*} \tag{2.2.29}
\| \sum_{\kappa \in \Z_+^d(k)} \sigma_\kappa \cdot (\mathcal
E_\kappa f ) \|_{L_p(\R^d)} \le c_{23} \| f \|_{L_p(\R^d)},
\sigma_\kappa = \prod_{j=1}^d \sigma^j_{\kappa_j}, \\ \text{ где }
\sigma^j_{\kappa_j} \in \{-1, 1\}, \kappa_j
\in \Z_+, j =1,\ldots,d,  f \in L_p(\R^d),
c_{23} = \prod_{j=1}^d c_{1}(\phi_j, \phi_j^*, p).
\end{multline*}

Теперь убедимся в справедливости (2.2.26).
Для произвольного семейства чисел $ \{ \sigma_\kappa = \prod_{j=1}^d
\sigma^j_{\kappa_j}: \sigma^j_{\kappa_j} \in \{-1, 1\}, j =1,\ldots,d,
\kappa \in \Z_+^d \}, $ функции $ f \in L_p(\R^d) $ рассмотрим последовательность
$$
\{ (\sum_{\kappa \in \Z_+^d(k \e)} \sigma_\kappa \cdot
(\mathcal E_\kappa f )) \in L_p(\R^d), k \in \Z_+ \}
$$
и, принимая во внимание (2.2.29), секвенциальную компактность
шара $ B(L_p(\R^d)) $ относительно $*$-слабой топологии в пространсте
$ L_p(\R^d) = (L_{p^\prime}(\R^d))^*, $ выберем
подпоследовательность
$$
\{ \sum_{\kappa \in \Z_+^d(k_n \e)} \sigma_\kappa \cdot
(\mathcal E_\kappa f ): k_n < k_{n+1}, n \in \N \}
$$
и функцию $ F \in L_p(\R^d), $ обладающие тем свойством, что для
любой функции $ g \in L_{p^\prime}(\R^d) $ выполняется равенство
\begin{equation*} \tag{2.2.30}
\lim_{ n \to \infty} \int_{\R^d} (\sum_{\kappa \in \Z_+^d(k_n \e)}
\sigma_\kappa (\mathcal E_\kappa f )) g dx =
\int_{\R^d} F g dx.
\end{equation*}
Заметим, что при любом $ \kappa \in \Z_+^d, $ благодаря (1.4.9),
(2.2.30), (1.4.10), для  $ g \in L_{p^\prime}(\R^d) $ имеет место равенство
\begin{multline*}
\int_{\R^d} (\mathcal E_\kappa F) \cdot \overline g dx = \int_{\R^d} F
\cdot \overline {(\mathcal E_\kappa g)} dx = \lim_{ n \to \infty}
\int_{\R^d} (\sum_{\kappa^\prime \in \Z_+^d(k_n \e)}
\sigma_{\kappa^\prime} \cdot (\mathcal E_{\kappa^\prime} f))
\overline{(\mathcal E_\kappa g)} dx \\
= \lim_{ n \to \infty} \int_{\R^d} \mathcal E_\kappa
(\sum_{\kappa^\prime \in \Z_+^d(k_n \e)}
\sigma_{\kappa^\prime} \cdot (\mathcal E_{\kappa^\prime} f)) \cdot
\overline g dx \\ = \lim_{ n \to \infty} \int_{\R^d}
(\sum_{\kappa^\prime \in \Z_+^d(k_n \e)} \sigma_{\kappa^\prime}
\cdot \mathcal E_\kappa^{p} (\mathcal E_{\kappa^\prime}^{p} f))
\overline g dx  = \int_{\R^d} \sigma_{\kappa} \cdot (\mathcal
E_\kappa f) \overline g dx,
\end{multline*}
и, значит,
\begin{equation*} \tag{2.2.31}
\mathcal E_\kappa F =
\sigma_{\kappa} \cdot (\mathcal E_\kappa f).
\end{equation*}
Учитывая (2.2.31), (1.1.1), (1.4.7), заключаем, что
\begin{equation*}
\sum_{\kappa \in \Z_+^d(k)} \sigma_\kappa \cdot (\mathcal
E_\kappa f ) = \sum_{\kappa \in \Z_+^d(k)} (\mathcal
E_\kappa F) = E_k^{p} F
\end{equation*}
сходится к $ F $ в $ L_p(\R^d) $ при $ \mn(k) \to \infty. $
Поэтому, переходя к пределу при $ \mn(k) \to \infty $ в
неравенстве (2.2.29), приходим к (2.2.26). $ \square $

Следствие

В условиях теоремы 2.2.2 для любого семейства чисел $ \{ \sigma_\kappa:
\kappa \in \Z_+^d \} $ вида $ \sigma_\kappa = \prod_{j=1}^d
\sigma^j_{\kappa_j}, $ где $ \sigma^j_{\kappa_j} \in \{-1, 1\}, \kappa_j
\in \Z_+, j =1,\ldots,d, $ и любой функции $ f \in L_p(\R^d) $
соблюдается неравенство
\begin{equation*} \tag{2.2.32}
\| f \|_{L_p(\R^d)} \le c_{23}
\| \sum_{\kappa \in \Z_+^d} \sigma_\kappa \cdot
(\mathcal E_\kappa f ) \|_{L_p(\R^d)}.
\end{equation*}

Доказательство.

Сначала покажем, что при любом $ k \in \Z_+^d, $ любом наборе чисел
$ \{ \sigma_\kappa = \prod_{j=1}^d \sigma^j_{\kappa_j}: \sigma^j_{\kappa_j} \in \{-1, 1\}, j =1,\ldots,d,
\kappa \in \Z_+^d(k) \}, $ для $ f \in L_p(\R^d) $
справедливо неравенство
\begin{equation*} \tag{2.2.33}
\| E_k f \|_{L_p(\R^d)} \le c_{23}
\| \sum_{\kappa \in \Z_+^d(k)}
\sigma_\kappa \cdot (\mathcal E_\kappa f ) \|_{L_p(\R^d)}.
\end{equation*}

В самом деле, ввиду (1.4.10), (1.1.1) имеем
\begin{multline*}
(\sum_{\kappa \in \Z_+^d(k)} \sigma_\kappa \cdot \mathcal
E_\kappa^{p})^2 f = (\sum_{\kappa, \kappa^\prime \in \Z_+^d(k)}
\sigma_\kappa \sigma_{\kappa^\prime} \cdot \mathcal E_\kappa^{p}
\mathcal E_{\kappa^\prime}^p) f \\
= (\sum_{\kappa \in \Z_+^d(k)} \sigma_\kappa^2 \cdot \mathcal E_\kappa^{p}) f =
\sum_{\kappa \in \Z_+^d(k)} \mathcal E_\kappa f = E_k f.
\end{multline*}
Откуда, применяя (2.2.29), выводим
\begin{multline*}
\| E_k f \|_{L_p(\R^d)} = \| (\sum_{\kappa \in \Z_+^d(k)}
\sigma_\kappa \cdot \mathcal E_\kappa)^2 f \|_{L_p(\R^d)}\\
 =\| \sum_{\kappa \in \Z_+^d(k)} \sigma_\kappa \cdot \mathcal
E_\kappa (\sum_{\kappa^\prime \in \Z_+^d(k)}
\sigma_{\kappa^\prime} \cdot (\mathcal E_{\kappa^\prime} f))
\|_{L_p(\R^d)} \\
\le c_{23} \| \sum_{\kappa^\prime \in \Z_+^d(k)}
\sigma_{\kappa^\prime} \cdot (\mathcal E_{\kappa^\prime} f)
\|_{L_p(\R^d)}.
\end{multline*}

Как видно из вывода (2.2.26) и (1.4.7), в неравенстве (2.2.33) можно
перейти к пределу при $ \mn(k) \to \infty, $ в результате чего
получим (2.2.32). $ \square $

С помощью теоремы 2.2.2 и следствия из нее, опираясь на схему доказательства
теоремы Литтлвуда-Пэли, изложенную в [3] для
операторов взятия частных сумм кратных рядов Фурье, устанавливается

Теорема 2.2.3

В условиях теоремы 2.2.2 существуют константы $ c_{24}(d,\phi,\phi^*,p) >0,\\
c_{25}(d,\phi,\phi^*,p) >0 $ такие, что для любой функции $ f \in L_p(\R^d) $
выполняются неравенства
\begin{equation*} \tag{2.2.34}
c_{24} \| f \|_{L_p(\R^d)} \le
( \int_{\R^d} (\sum_{\kappa \in \Z_+^d} |(\mathcal E_\kappa f)(x)|^2 )^{p/2}
dx)^{1/p} \le c_{25} \| f \|_{L_p(\R^d)}.
\end{equation*}

Доказательство.

Рассмотрим систему Радемахера, состоящую из функций
$$
\omega_\kappa(t) = \operatorname{sign} \sin( 2^{\kappa +1} \pi t),
t \in I, \kappa \in \Z_+,
$$
и определим семейство функций $ \omega_\kappa^d, \kappa \in \Z_+^d, $
полагая
$$
\omega_\kappa^d(t) = \prod_{j=1}^d \omega_{\kappa_j}(t_j), t \in I^d.
$$
Как известно (см., например, [3, п. 1.5.2]), существуют константы $ c_{26}(d,p) >0,
c_{27}(d,p) >0 $ такие, что при любом $ k \in \Z_+^d $ для любого набора
чисел $ \{ a_\kappa \in \C, \kappa \in \Z_+^d(k) \} $
имеет место неравенство
\begin{equation*} \tag{2.2.35}
c_{26} ( \sum_{ \kappa \in \Z_+^d(k)} | a_\kappa|^2 )^{1/2}
\le
\| \sum_{ \kappa \in \Z_+^d(k)} a_\kappa \omega_\kappa^d(\cdot) \|_{L_p(I^d)}
\le
c_{27} ( \sum_{ \kappa \in \Z_+^d(k)} | a_\kappa|^2 )^{1/2}.
\end{equation*}

Для $ f \in L_p(\R^d) $ при $ k \in \Z_+^d, $ используя (2.2.33), теорему
Фубини, (2.2.35), выводим
\begin{multline*} \tag{2.2.36}
\| E_k f \|_{L_p(\R^d)}^p = \int_{I^d} \| E_k f \|_{L_p(\R^d)}^p dt \le
\int_{I^d} (c_{23} \| \sum_{\kappa \in
\Z_+^d(k)} \omega_\kappa^d(t) \cdot (\mathcal E_\kappa f )
\|_{L_p(\R^d)})^p dt\\
 = (c_{23})^p \int_{I^d} \int_{\R^d} |
\sum_{\kappa \in \Z_+^d(k)} \omega_\kappa^d(t) \cdot (\mathcal
E_\kappa f )(x) |^p dx dt \\
= (c_{23})^p \int_{\R^d}
\int_{I^d} | \sum_{\kappa \in \Z_+^d(k)} \omega_\kappa^d(t) \cdot
(\mathcal
E_\kappa f )(x) |^p dt dx\\
 = (c_{23})^p \int_{\R^d} \|
\sum_{\kappa \in \Z_+^d(k)} (\mathcal E_\kappa f )(x) \cdot
\omega_\kappa^d(\cdot) \|_{L_p(I^d)}^p dx \\
\le (c_{23})^p \int_{\R^d} (c_{27} (\sum_{\kappa \in \Z_+^d(k)}
|(\mathcal E_\kappa f )(x)|^2 )^{1/2})^p dx\\
 = (c_{28})^p \int_{\R^d}
(\sum_{\kappa \in \Z_+^d(k)} |(\mathcal E_\kappa f )(x)|^2
)^{p/2} dx,
\end{multline*}
и, пользуясь (2.2.35), теоремой Фубини, (2.2.29), получаем
\begin{multline*}
\int_{\R^d} (\sum_{\kappa \in \Z_+^d(k)} |(\mathcal E_\kappa f )(x)|^2 )^{p/2} dx
\le \int_{\R^d} (c_{29} \| \sum_{\kappa \in \Z_+^d(k)}
(\mathcal E_\kappa f )(x) \cdot \omega_\kappa^d(\cdot) \|_{L_p(I^d)})^p dx\\
= (c_{29})^p \int_{\R^d} \int_{I^d} | \sum_{\kappa \in \Z_+^d(k)}
\omega_\kappa^d(t) \cdot (\mathcal E_\kappa f )(x) |^p dt dx\\
= (c_{29})^p \int_{I^d} \int_{\R^d} | \sum_{\kappa \in
\Z_+^d(k)} \omega_\kappa^d(t) \cdot (\mathcal E_\kappa f)(x) |^p dx dt\\
= (c_{29})^p \int_{I^d} \| \sum_{\kappa \in \Z_+^d(k)}
\omega_\kappa^d(t) \cdot (\mathcal E_\kappa f )
\|_{L_p(\R^d)}^p dt \\
\le (c_{29})^p \int_{i^d} (c_{23} \| f \|_{L_p(\R^d)})^p dt =
(c_{25})^p \| f \|_{L_p(\R^d)}^p,
\end{multline*}
откуда, в частности, имеем
\begin{equation*} \tag{2.2.37}
\int_{\R^d} (\sum_{\kappa \in \Z_+^d(n \e)}
|(\mathcal E_\kappa f )(x)|^2 )^{p/2} dx \le (c_{25})^p
\| f\|_{L_p(\R^d)}^p, n \in \Z_+.
\end{equation*}

Для получения второго неравенства в (2.2.34) достаточно применить
теорему Леви о предельном переходе под знаком интеграла к монотонно
возрастающей последовательности функций
$ \{ (\sum_{\kappa \Z_+^d(n \e)} |(\mathcal E_\kappa f)(x)|^2 )^{p/2},
n \in \Z_+\}, $ принимая во внимание (2.2.37) и учитывая, что
почти для всех $ x \in \R^d $ предел
\begin{multline*}
\lim_{n \to \infty} \biggl(\sum_{\kappa \in \Z_+^d(n \e)} |(\mathcal
E_\kappa f)(x)|^2\biggr )^{p/2} =
\biggl(\lim_{n \to \infty} \sum_{\kappa \in \Z_+^d(n \e)}
|(\mathcal E_\kappa f)(x)|^2\biggr )^{p/2} \\ =
\biggl(\sum_{\kappa
\in \Z_+^d} |(\mathcal E_\kappa f)(x)|^2 \biggr)^{p/2}.
\end{multline*}
Последнее равенство справедливо в силу того, что при $ k \in
\Z_+^d $ имеет место соотношение
\begin{multline*}
\sum_{\kappa \in \Z_+^d(\mn(k) \e)} |(\mathcal E_\kappa f)(x)|^2 \le
\sum_{\kappa \in \Z_+^d(k)} |(\mathcal E_\kappa f)(x)|^2 \le
\sum_{\kappa \in \Z_+^d(\mx(k) \e)} |(\mathcal E_\kappa f)(x)|^2, x \in \R^d.
\end{multline*}

Ввиду сказанного и на основании (1.4.7) переходя к пределу в (2.2.36) при
$ \mn(k) \to \infty, $ приходим
к первому неравенству в (2.2.34). $ \square $

\end{document}